\documentclass[12pt,reqno]{amsart}

\numberwithin{equation}{section}

\usepackage{mdwlist}

\usepackage{nameref}

\usepackage[final, 
color]{showkeys}
\definecolor{refkey}{gray}{.85}
\definecolor{labelkey}{gray}{.85}

\makeatletter
\let\orgdescriptionlabel\descriptionlabel
\renewcommand*{\descriptionlabel}[1]{%
  \let\orglabel\label
  \let\label\@gobble
  \phantomsection
  \edef\@currentlabel{#1}%
  \let\label\orglabel
  \orgdescriptionlabel{#1}%
}
\makeatother

\usepackage{AKstyle}

\usepackage{caption}
\usepackage[labelfont=rm]{subcaption}

\usepackage[pagebackref=true, colorlinks]{hyperref}

\hypersetup{pdffitwindow=true,linkcolor=blue,citecolor=blue,urlcolor=blue,menucolor=blue}

\usepackage{comment}

\theoremstyle{plain}

\newtheorem*{conja}{Conjecture A}
\newtheorem*{thma}{Theorem A}

\newtheorem*{conjz}{Conjecture Z}
\newtheorem*{thmz}{Theorem Z}

\newtheorem*{thmp}{Theorem P}

\begin{document}

\author{Alex Kontorovich}
\thanks{Partially supported by NSF grants  DMS-1209373, DMS-1064214 and DMS-1001252.}
\email{alex.kontorovich@yale.edu}
\address{Department of Mathematics, Yale University, New Haven, CT}

\title[From Apollonius to Zaremba]{From Apollonius To Zaremba: Local-Global Phenomena in Thin Orbits}

\begin{abstract}
%We survey
%In this survey, 
%w
We %describe
discuss
 a number of naturally arising problems in arithmetic, % which are 
 culled from completely unrelated sources, %and 
% yet
which
  turn out to have a common formulation involving ``thin'' orbits. 
These include the
 %Strong Density Conjecture
 local-%to-
 global problem for integral 
 Apollonian gaskets and Zaremba's Conjecture on finite continued fractions with %uniformly 
absolutely
bounded partial quotients.
Though
these 
%such
problems  could have been posed by the ancient Greeks, 
recent
progress
 %had to wait for 
%uses modern
 %tools 
comes from a pleasant synthesis of
modern
techniques
 from
a variety of
fields, including 
% more % some % of % the latest
%modern
% technology, %which only recently became available,
%specifically 
%in
%infinite volume spectral, ergodic and representation theory, additive combinatorics and expanders,
%and some modern variants of the circle method. 
%in 
harmonic 
analysis,
algebra, 
geometry,
combinatorics,
 and dynamics.
We describe the problems, partial progress, and some of the
 %techniques
 tools alluded to above.
\end{abstract}
\date{\today}

\subjclass[2010]{11F41, 11J70, 11P55, 20H10, 22E40}
\maketitle
\tableofcontents

\newpage

\section{Introduction}

In this article we will discuss recent developments on several
seemingly unrelated arithmetic problems, which each boil down to the
same issue of proving a ``local-global principle for thin orbits''. In
each of these problems, we study the\emph{ orbit} 
$$
\cO=\G\cdot\bv_{0},
$$ 
of some given
vector $\bv_{0}\in\Z^d$, under the action of some  given % set%ubgroup
group or semigroup%
, $\Gamma$,
(under multiplication) of %the 
$d$-by-$d$ integer matrices. % ; that is
It will turn out that
the %se naturally-arising 
orbits
arising naturally
in our problems
 %will be 
 are
 ``thin'';
%in a sense which we will formally define later. 
%, meaning r
%R
roughly 
speaking, this means that
%in the sense that they 
%the 
each
orbit
is ``degenerate'' in
%their
its
 algebro-geometric closure,
containing %very
relatively
very few points. % in $\Z^{d}$. %, relative to

% the ambient set in which it naturally sits. % contains the orbit.

%The 
Each of the
 problems %each
  then takes 
%ing 
another vector $\bw_{0}\in\Z^{d}$, %we 
and
for the standard inner product $\<\cdot,\cdot\>$ on $\R^{d}$,
forms
%ing 
the set 
$$
\sS:=\<\bw_{0},
%\G\cdot\bv_{0}
\cO
\>
\subset\Z
$$ 
of integers, asking
% whether %it 
what numbers are in
$\sS$.
%satisfies a ``local-global principle''.
\begin{comment}

\newpage

The goal of these notes is to survey recent developments on several {\it a priori} unrelated problems about whole numbers. It will turn out that the problems fall under a common umbrella, in that each asks for what we will refer to as a ``local-global principle for thin orbits.'' Let us explain.
\\

%First recall
We first clarify % explain % the notion of 
what we mean by
an orbit. Suppose some set $\G$ {\it acts} on a set $X$, that is, each $\g\in\G$ is a map $\g:X\to X$. Then, fixing some base point $x_{0}\in X$, one can study the {\it orbit} $\cO$ of $x_{0}$ under $\G$, namely,
$$
\cO:=\G x_{0}=\{\g(x_{0}):\g\in\G\}\subset X.
$$
We shall explain in due course
% explain 
what we mean by a ``thin'' orbit, but for now one can have in mind that such an orbit is somehow ``degenerate,'' 
in that 
%the vast technology developed to understand orbits in great generality does not apply to the particular orbits which will arise in our applications.
%having asymptotic
 %``thin'' orbit is one which
 it has asymptotically many fewer points than $X$.
\\

Next %recall
we describe what one means by
 %the notion of 
 a ``local-global principle.'' Say we have a set $\sS\subset\Z$ of %positive 
 integers. 
 We would like to know, for a particular integer $n$, whether $n\in \sS$.
\end{comment}
For an integer $q\ge1$, the projection map 
$$
\Z\to\Z/q\Z
$$ 
%gives
can give an obvious obstruction to membership.
 Let $\sS(\mod q)$ be the image of this
projection, % map,
$$
\sS(\mod q):=\{s(\mod q):s\in \sS\}\subset \Z/q\Z. %,
$$ 
 %in $\sS$. 
%Then 
%So
%f
For example, 
suppose that any number in $\sS$ leaves a remainder of $1,2$ or $3$ when divided by $4$, that is,
 $\sS(\mod 4)=\{1,2,3\}$.
% that is,
%if every number in $\sS$ is odd, 
Then one can % immediately %answer whether 
 conclude, without any further consideration, that
 $10^{10^{10}}\notin\sS$, since $10^{10^{10}}\equiv0(\mod 4)$.
This is called a {\it local} obstruction. 
%  Such an 
%  An obvious obstruction to membership in $\sS$ is if $n$ fails to be in the 
 %reduction 
 %For an integer $q\ge1$, l
%for some integer $q\ge1$. 
%We c
Call $n$ {\it admissible} if it %passes 
avoids
all local obstructions,
$$
n\in \sS(\mod q),\qquad\text{ for all $q\ge1$.}
$$
%Note that
%It %may not be int
%is not 
%is by no means
%obvious, but i
In many applications,
 the set $\sS(\mod q)$ is %often %infinitely 
significantly 
easier to analyze than the set $\sS$ itself. But a local to global phenomenon predicts that, if $n$ is admissible, then in fact $n\in \sS$, thereby reducing the seemingly %much 
more difficult problem to the easier one. 
\\

It is the combination of these
 %ideas
 concepts, $(i)$ thin orbits, and $(ii)$ local-global phenomena, which will turn out to be %at 
 the %heart
 ``beef''
  of the problems we intend to discuss.
\\

\subsection{Outline}\

We %proceed % in reverse alphabetical order, starting
begin
 in \S\ref{sec:Z} with Zaremba's Conjecture. We will explain how this problem arose naturally in the study of ``good lattice points'' for quasi-Monte Carlo methods in multi-dimensional numerical integration, and how it also has applications to the linear congruential method for pseudo-random number generators. But the assertion of the conjecture  is a statement about continued fraction expansions of rational numbers, and as such is so elementary that Euclid himself could have posed it. We will discuss recent progress by Bourgain and the author, proving 
 %$100\%$ 
a density version
 of the conjecture. %(but falling far short of settling it outright).

\begin{figure}
\includegraphics[width=2.5in]{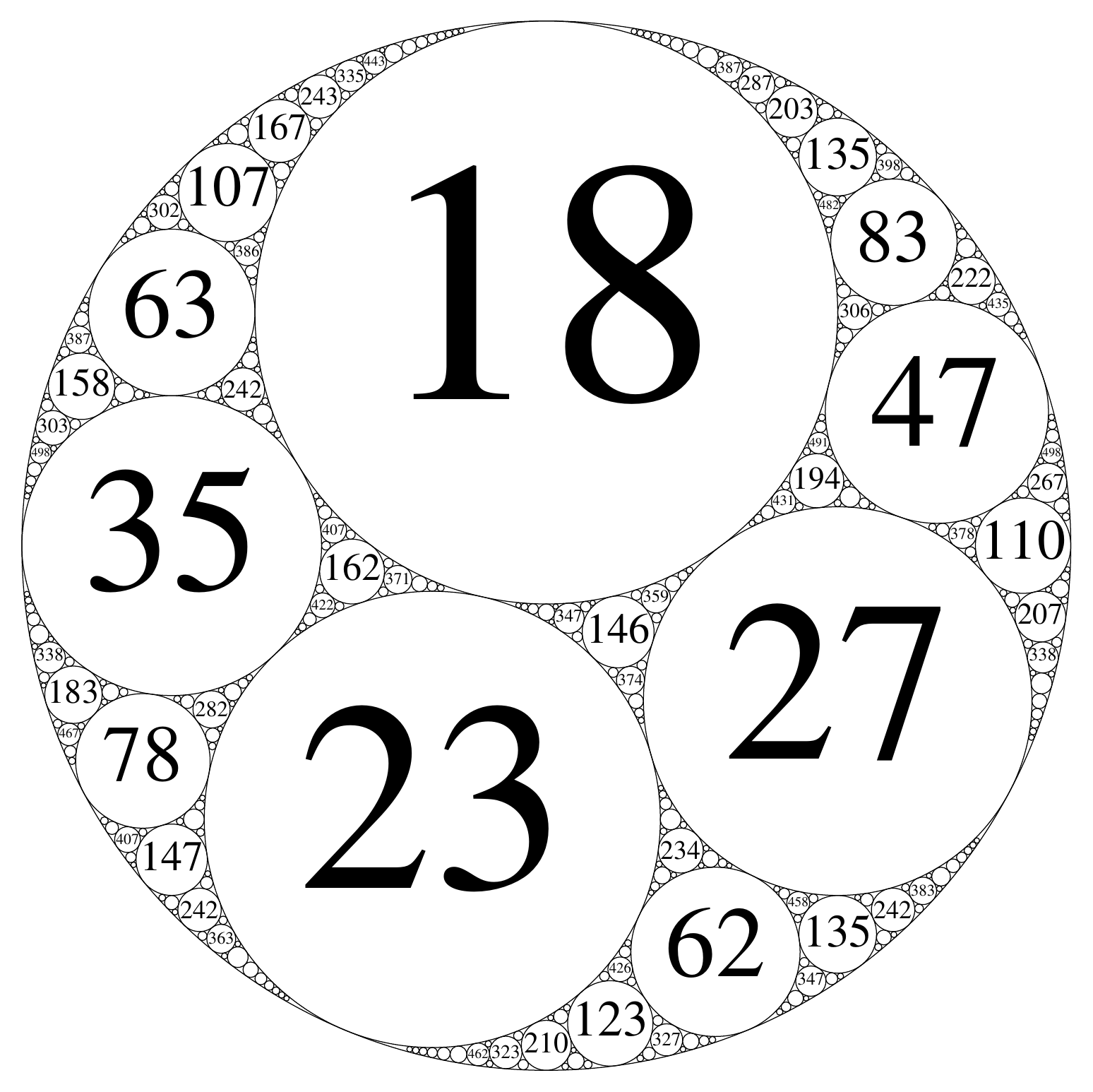}
\caption{An integral Apollonian gasket.}
\label{fig:Apoll}
\end{figure}

We % move from one ancient Greek to another, changing 
change
our focus in \S\ref{sec:A} to the %great 
ancient % Greek 
geometer
Apollonius of Perga.
As we will explain, his straight-edge and compass construction of tangent 
%(``kissing'') 
circles, when iterated {\it ad infinitum}, gives rise to a beautiful fractal circle packing in the plane,
% see 
such as that shown in
Figure \ref{fig:Apoll}. Recall that the curvature of a circle is just one over its radius. For special configurations, all the curvatures of circles in the given packing turn out to be integers; these are the numbers shown in Figure \ref{fig:Apoll}. We will present in \S\ref{sec:A} 
%recent 
progress on the problem: which integers appear?
%As we will explain, i
It was recently proved by Bourgain and the author that %$100\%$ of the 
almost every
admissible number appears.

In \S\ref{sec:P}, also stemming from Greek mathematics, we describe a local-global problem for a thin orbit of Pythagorean triples, as will be defined there. 
This problem is a %refinement 
variant
of the so-called Affine Sieve, recently introduced by Bourgain, Gamburd, and Sarnak. We will 
%briefly discuss the Affine and some more classical sieves, before 
explain %ing
 %the solution, 
 % of %$100\%$ 
an % density 
``almost''
%version
%of the
 local-global
  %problem 
  theorem
  in %that 
  this
  context
  due to Bourgain and the author.

Finally, these three problems are reformulated to the aforementioned common umbrella in \S\ref{sec:C}, where some of the ingredients of the proofs are sketched. 
The %se 
problems do not naturally fit in an established area of research, having no $L$-functions or Hecke theory (though they are unquestionably problems about whole numbers), being not part of the Langlands Program (though involving automorphic forms and representations), nor falling under the purview of the classical circle method or sieve, which
attempt to solve equations or produce primes in polynomials (here it is not polynomials that %produce 
generate
points, but the %se %affine linear 
%linear
aforementioned
matrix
actions%, as we will explain
). 
Instead %they tend to 
one must
borrow bits and pieces from %a variety of 
these 
fields and others.
The 
major 
tools which we aim to highlight
%there and 
throughout
% are
include
analysis
(the circle method, exponential sum bounds, %and
infinite volume spectral theory),
algebra (%
%infinite-dimensional representations of  Lie groups,
strong approximation, 
Zariski density,
spin and orthogonal groups associated to quadratic forms, %, 
representation theory%
), 
geometry (hyperbolic manifolds, circle packings, %and 
diophantine approximation), 
combinatorics (sum-product, expander graphs, spectral gaps),
and dynamics (ergodic theory, mixing rates, %,
 %expander graphs, 
 the thermodynamic formalism%
 ).
%\\

%players are the circle method; bounds for exponential sums and bilinear forms;
%infinite volume spectral, representation, and ergodic  theory; the thermodynamical formalism and transfer operators; %additive combinatorics; and expander graphs.

\subsection{Notation}\

We use the following standard notation. A quantity is defined 
%by 
via
the symbol ``$:=$'', and a concept being defined is italicized. Write $f\sim g$ for $f/g\to1$, $f=o(g)$ for $f/g\to0$, and $f=O(g)$ or $f\ll g$ for $f\le C g$. Here $C>0$ is called an implied constant, and 
is absolute unless otherwise specified. Moreover, $f\asymp g$ means $f\ll g \ll f$. We use $e(x)=e^{2\pi i x}$.
 %and $e_{q}(x)=e(x/q)$. 
 The cardinality of a finite set $S$ is written as $|S|$ or $\#S$.
The transpose of a vector $\bv$ is written $\bv^{t}$. The meaning of
algebraic
 symbols can change from section to section; for example the (semi)group $\G$ and quadratic form $Q$
 will
  vary
   %from section to section, 
   depending on the context.
%\\

\subsection*{Acknowledgements}\

I wish to thank Andrew Granville for %inviting and 
encouraging me to pen these notes and for his
insightful
and detailed
 input on 
 %the 
 %an 
% early versions. %s. % of this file
various drafts.
%I am grateful
Thanks  to  Mel Nathanson for %asking 
inviting
me to give a mini-course at CANT 2012, %thanks to 
as a result of
which these notes were finally assembled.
%, and to Steven J. Miller for TeXing my lectures in real time. 
%Thanks to 
I am %infinitely 
grateful to
Peter Sarnak for %his continued inspiration, counsel, and for 
introducing me to Apollonian gaskets and infinite volume spectral methods, %. Thanks also 
to Hee Oh for introducing me to homogeneous dynamics, and to Dorian Goldfeld for his constant support and advice.
Thanks %; and 
to Elena Fuchs, Aryeh Kontorovich, and Sam Payne for comments on an earlier draft. Most of all, 
%words cannot express 
I owe 
a huge debt of gratitude 
to Jean Bourgain for his
 %kind and 
 generous 
 tutelage
 %, encouragement, 
 and collaboration.
\\

\newpage

\section{Zaremba's Conjecture}\label{sec:Z}

 \begin{figure}
        \begin{subfigure}[t]{0.46\textwidth}
                \centering
		\includegraphics[width=\textwidth]{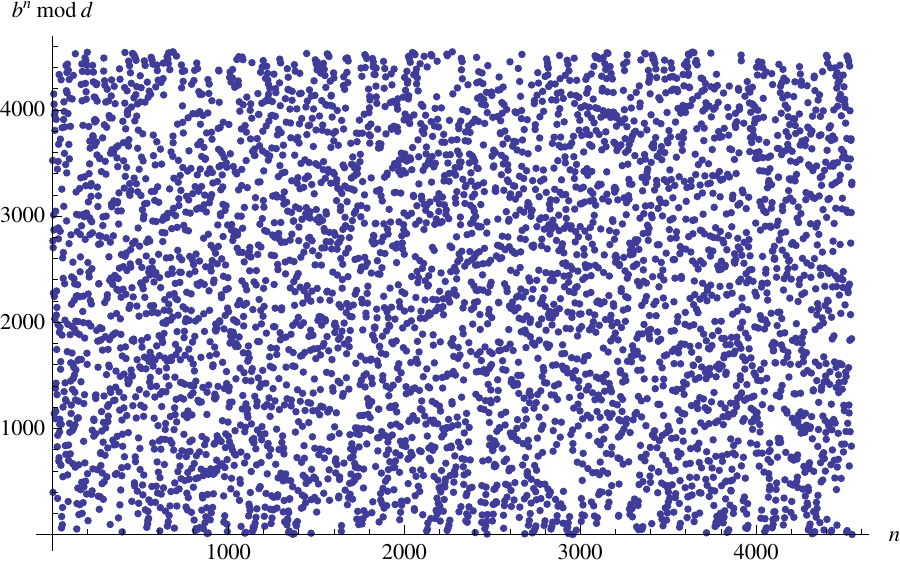}
                \caption{Multiplier $b=3523$.}
                \label{fig:bnModDA}
        \end{subfigure}%
%\qquad
\qquad
        \begin{subfigure}[t]{0.46\textwidth}
                \centering
		\includegraphics[width=\textwidth]{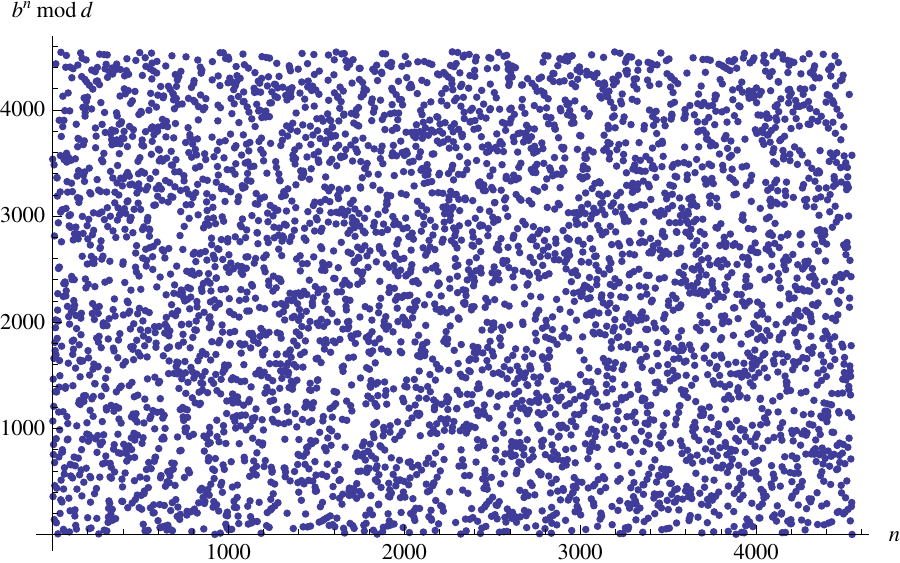}
                \caption{Multiplier $b=3535$.}
                \label{fig:bnModDB}
        \end{subfigure}

\caption{Graphs of the map \eqref{eq:nToBn} with prime modulus $d=4547$, and multiplier $b$ as shown.}
\label{fig:bnModD}
\end{figure}

%The digital age needs computers to 
%produce 
%spit out
%gazillions of random numbers for countless applications, 
%but being Turing machines (at least for now), the best we can ask
%of
% computers
% is
%  to produce 
  %are 
Countless applications require
  {\it pseudo-random numbers}: deterministic algorithms which ``behave randomly.'' 
  %whatever that means. 
  Probably the simplest, oldest, and best known among these is the so-called {\it linear congruential method}: For some starting seed $x_{0}$, iterate the map
\be\label{eq:linCong}
x\mapsto bx+c\quad(\mod d)
.
\ee
Here $b$ is called the multiplier, $c$ the shift, and $d$ the modulus. For simplicity, we consider the homogeneous case $c=0$. % Now number theory interferes, but we can essentially ignore it if we 
To have as long a sequence as possible,
take $d$ to be prime, and $b$ a primitive root mod $d$, that is, a generator of the cyclic group $(\Z/d\Z)^{\times}$. In this case we may as well start with the seed $x_{0}=1$; then the iterates of \eqref{eq:linCong} are nothing more than the map
\be\label{eq:nToBn}
n\mapsto b^{n}\quad(\mod d)
.
\ee
We show graphs of this map in Figure \ref{fig:bnModD} for the prime $d=4547$, with two choices of roots $b=3523$ and $b=3535$.
%Our eyes confirm i
In both cases, % that 
the graphs ``look'' random, in that, given $b$ and $n$,
 it is hard to guess where $b^{n}(\mod d)$ will lie (without just computing% it
 ).
 % (except that both sequences end at $1$, a consequence of Fermat's little theorem).
Similarly, given $b$ and $b^{n}(\mod d)$, it is 
%very 
typically
%hard 
difficult
to determine $n$; this is the classical
 %and difficult 
 problem of computing a discrete logarithm.

\begin{figure}
        \begin{subfigure}[t]{0.35\textwidth}
                \centering
		\includegraphics[width=\textwidth]{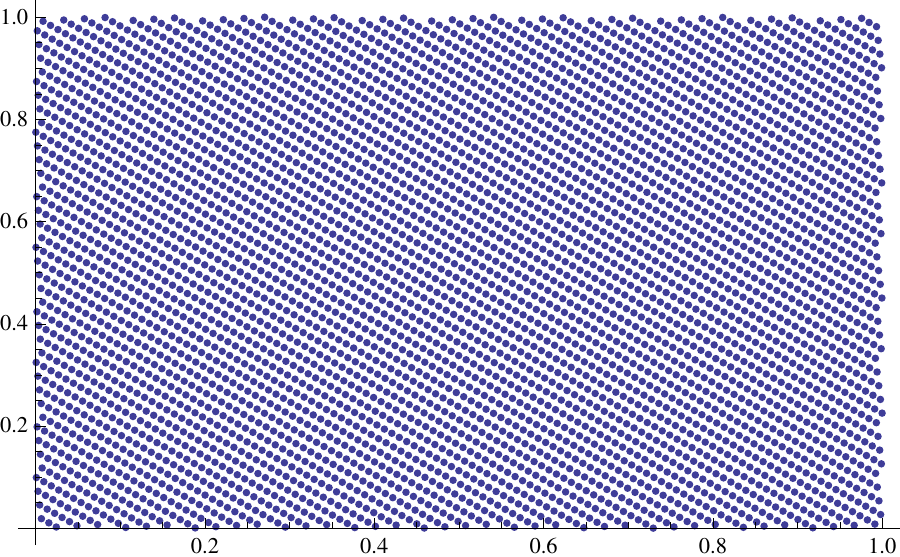}
                \caption{Multiplier $b=3523$.}
                \label{fig:pairsA}
        \end{subfigure}%
%\qquad
\qquad
        \begin{subfigure}[t]{0.35\textwidth}
                \centering
		\includegraphics[width=\textwidth]{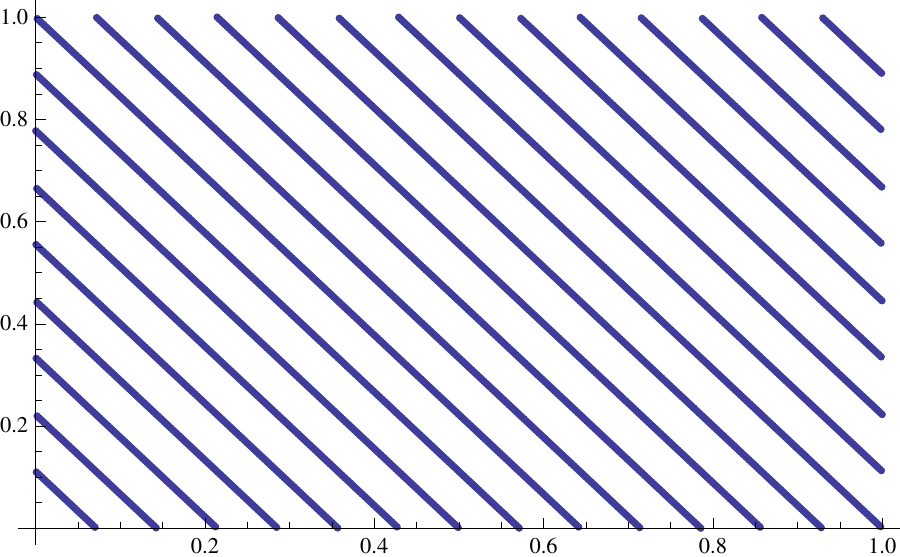}
                \caption{Multiplier $b=3535$.}
                \label{fig:pairsB}
        \end{subfigure}

\caption{Plots of the points \eqref{eq:bnSeq} for the same choices of modulus $d=4547$ and multipliers as in Figure \ref{fig:bnModD}.}
\label{fig:pairs}
\end{figure}

%Second to our eyes, the next simplest (and %slightly 
A slightly more rigorous statistical test for randomness is the serial correlation of pairs: how well can we guess
where
 $b^{n+1}$
 is, knowing $b^{n}$? To this end, we plot in Figure \ref{fig:pairs} these pairs, or what is the same, the pairs
\be\label{eq:bnSeq}
\left\{
\left(
{b^{n}\over d}
,
{b^{n+1}\over d}
\right)
(\mod 1)
\right\}_{n=1}^{d}%-1}
\subset
\R^{2}/\Z^{2}
\ee
in the unit square, with the previous choices of modulus and multiplier. Focus first on Figure \ref{fig:pairsA}: it looks like a fantastically equidistributed grid. Keep in mind that the mesh in each coordinate is of size $1/d\approx1/4000$, so we have $(4000)^{2}$ points from which to choose, yet we are only plotting $4000$ points,  square-root the total number of options.

On the other hand, look at Figure \ref{fig:pairsB}: these parameters make a terrible random number generator! Knowing $b^{n}$, we have about a $1:10$ chance of guessing $b^{n+1}$, not $1:4000$.
\\

%This %type of problem 
A related
phenomenon
%is
also
 %important
appears
  in
  two-dimensional
   numerical integration: %ng a two
%dimensional function: 
Suppose that you wish to integrate a ``nice'' function $f$ on 
$\R^{2}/\Z^{2}\cong
[0,1) \times [0,1)$, say of bounded variation, $V(f)<\infty$. 
The idea is to take a large % (but practically sized)
%finite
sample of points $\cZ$ in $\R^{2}/\Z^{2}$, and approximate the integral by the average
of $f(z)$,   $z\in\cZ$.  For this to be a good approximation one obviously
needs that $f$ does not vary much in a small ball, and that the points
of $\cZ$ are well-distributed throughout $\R^{2}/\Z^{2}$. In fact, Koksma and Hlawka
showed, rather beautifully, that this is all that one needs to take
into account:
% \begin{thm}[Koksma-Hlawka Inequality]
$$
\left|
\int_{0}^{1}
\int_{0}^{1}
f(x,y)
dx\,dy
-
{1\over |\cZ|}
\sum_{z\in\cZ}
f(z)
\right|
\le%l
V(f)
\cdot
\Disc(\cZ)
.
$$
%\end{thm}
Here
%$V(f)$ is the total variation of $f$, and
 $\Disc$ is the {\it discrepancy} of the set $\cZ$, %; %; recall this is computed
%which
%the latter is
 defined
 as follows. Take a rectangle 
 $
 R=[a,b]\times[c,d]\subset \R^{2}/\Z^{2}.
 $ 
 One would like the fraction of points in $R$ to be close to its area, so set
$$
\Disc(\cZ):=\sup_{R\subset \R^{2}/\Z^{2}}
\left|
{\#(\cZ\cap R)\over \#\cZ}
-
\Area(R)
\right|
.
$$

It is elementary that for a growing family $\cZ^{(k)}\subset \R^{2}/\Z^{2}$, $|\cZ^{(k)}|\to\infty$, the discrepancy $\Disc(\cZ^{(k)})$ decays to $0$ if and only if $\cZ^{(k)}$ becomes equidistributed in $\R^{2}/\Z^{2}$. But the discrepancy itself is a finer measure of the rate of this decay.
For example, %it is %immediate 
%easy to see 
%that
%clearly
observe that for any finite sample set $\cZ$, we have the lower bound
 $\Disc(\cZ)\ge1/|\cZ|$. Indeed, take a family of rectangles $R%_{\vep}
$ zooming in on a single point in $\cZ$; the proportion of points in $R$  is always $1/|\cZ|$, while the area of $R%_{\vep}
$ can be made arbitrarily small. 
It turns out there is a sharpest possible lower bound, due to Schmidt \cite{Schmidt1972}:
\be\label{eq:Schmidt}
\text{For any finite $\cZ\subset S$,}
\qquad
\Disc(\cZ)\gg {\log |\cZ|\over |\cZ|}
.
\ee
%The 

For standard Monte Carlo integration, %method 
one
often
samples $z\in\cZ$ according to the uniform measure;
%Note also that for a  %``purely random'' 
%uniformly sampled
%set of points $\cZ$,
 the Central Limit Theorem then predicts that
\be\label{eq:Disc12}
\Disc(\cZ)\approx{1\over |\cZ|^{1/2}},
\ee
ignoring $\log\log$ factors. So 
comparing \eqref{eq:Disc12} to \eqref{eq:Schmidt},
%from this point of view, 
%``purely random'' 
it is clear that
uniformly sampled
sequences are far from optimal in numerical integration. %from 
%
%numerical integration.
%Likewise
Alternatively, one could 
%just
take $\cZ$ to be an evenly spaced $d$-by-$d$ %evenly spaced 
grid, 
$$
\cZ=\{(i/d,j/d):0\le i,j<d\}.
$$
But then the rectangle $[\vep,1/d-\vep]\times[0,1]$ contains no grid points while its area is almost $1/d=1/|\cZ|^{1/2}$,
again giving \eqref{eq:Disc12}. % This is again far from %optimal
%\eqref{eq:Schmidt}.

In the {\it qausi} Monte Carlo method,
rather than sampling uniformly, 
%one chooses
one tries to find
 a special sample set $\cZ$ to come as close as possible to the %ideal
% extremal 
optimal discrepancy
  \eqref{eq:Schmidt}.
  Ideally, such a set $\cZ$ would %itself 
  also be quickly and easily
constructible by a computer algorithm.
Not surprisingly, the set $\cZ$ illustrated in Figure \ref{fig:pairsA} makes an excellent sample set.
%
%, instead of sampling uniformly, tries to find a sequence $\cZ$ (or increasing family of such)
%
%
%\newpage
%
%, and i
It was this problem which led Zaremba to
 his theorem and 
 %the 
 conjecture, %we
described below.
\\

Returning % from this interlude 
to our initial discussion, observe that the sequence \eqref{eq:bnSeq} is essentially (since $b$ is a generator) the same as
\be\label{eq:Zd}
%\cZ_{d}=
\cZ_{b,d}:=
\left\{
\left(
{{n}\over d}
,
{b{n}\over d}
\right)
\right\}_{n=1}^{d}
(\mod 1)
.
\ee
And this is nothing more than a graph of our first map %s 
\eqref{eq:linCong}.
%Then
Now it is 
%obvious 
clear
that both Figures \ref{fig:pairsA} and \ref{fig:pairsB} are %projective 
``lines'', but the first must be ``close to a line with irrational slope,'' causing the equidistribution. 
This Diophantine property is best described in terms of continued fractions, %which we now recall.
as follows.

For $x\in(0,1)$, we use the notation 
$$
x=[a_{1},a_{2},\dots]
$$ 
for the continued fraction expansion
$$
x=\cfrac{1}{a_{1}+\cfrac{1}{a_{2}+\ddots}}
$$
The %numbers 
integers
$a_{j}\ge1$ are called {\it partial quotients} of $x$.
%, or digits. 
%As is
 %very 
 %$well-known, r
 Rational numbers have finite continued fraction expansions.

%Thus we are immediately
One is then immediately 
 prompted to 
 %check
 study
 the continued fraction 
%expansion 
expansions of the ``slopes'' $b/d$ in Figure \ref{fig:pairs}:
\beann
 {3523}/{4547}&=&[1,3,2,3,1,2,3,2,1,3],\\
 {3535}/{4547}&=&[1,3,2,35,1,1,1,4].
\eeann
Note the very large partial quotient $35$ in the middle of the second expression, while the partial quotients in the first
are all at most $3$.
 Observations of this kind naturally led Zaremba to the following
\begin{thm}[{Zaremba 1966 \cite[Corollary 5.2]{Zaremba1966}}]
%Let
Fix
 $(b,d)=1$ with $b/d=[a_{1},a_{2},\dots,a_{k}]$ and 
%all $a_{j}\le A$. 
let $A:=\max a_{j}$.
Then for $\cZ_{b,d}$ given in \eqref{eq:Zd},
\be\label{eq:Zarem}
\Disc(\cZ_{b,d})\le
\left(
{4A\over \log (A+1)}
+
{4A+1\over \log d}
\right)
{\log d\over d}
.
\ee
\end{thm}

Since $|\cZ_{b,d}|=d$, comparing \eqref{eq:Zarem} to \eqref{eq:Schmidt} shows that the sequences \eqref{eq:Zd} are essentially best possible, up to the ``constant'' $A$, cf. Figure \ref{fig:pairsA}. 
%Of course 
But
the previous sentence is complete nonsense: $A$ is not constant at all; it depends on $d$,%
\footnote{%Of course 
%Note that
The value $A$ also depends on $b$, but %it is clear that 
the important variable for applications is 
$|\cZ_{b,d}|=d$.} 
cf. Figure \ref{fig:pairsB}.
\\

%$$
%\text{Or must it?}
%$$

With this motivation, Zaremba predicted that in fact $A$ {\it can} be taken constant:
%\\

\phantomsection
\begin{conjz}[{Zaremba 1972 \cite[p. 76]{Zaremba1972}}]\label{conj:Z}
Every natural number is the denominator of a reduced fraction whose partial quotients are absolutely bounded.

That is, there exists some absolute $A>1$ so that for each $d\ge1$, there is some $(b,d)=1$, so that $b/d=[a_{1},\dots,a_{k}]$ with  $\max a_{j}\le A$.
\end{conjz}

Zaremba even suggested a sufficient value for $A$, namely $A=5$. So this is really a problem 
%the
 %ancient Greeks could have understood
that could have been posed 
  in Book VII of the {\it Elements} (after Euclid's algorithm):
 % : 
 using the 
 %digits
 partial quotients
  $a_{j}\in\{1,\dots,5\}$, does the set of (reduced) fractions with expansion $[a_{1},\dots,a_{k}]$ contain every integer as a denominator? The reason for Zaremba's guess $5$ is simply that it is false for $A=4$, as we now explain. First some more notation.

Let $\sR_{A}$ be the set of rationals with the desired property that all partial quotients are at most $A$:
$$
\sR_{A}:=
\left\{
\frac bd
=
[a_{1},\dots,a_{k}]
:
(b,d)=1,\text{ and } a_{j}\le A,\ \forall j%=1,\dots,k
\right\}
,
$$
and let $\sD_{A}$ be the set of denominators which arise:
$$
\sD_{A}
:=
\left\{
d:\exists (b,d)=1\text{ with }\frac bd\in\sR_{A}
\right\}
.
$$
Then Zaremba's conjecture is that $\sD_{5}=\N$, and we claim that this is false for $\sD_{4}$. Indeed, $6\notin\sD_{4}$: the only numerators to try are $1$ and $5$, but 
the continued fraction expansion of
$1/6$
is just
$[6]$, and $5/6=[1,5]$, so the largest partial quotient in both is too big.

That said, there are only two other numbers, $54$ and $150$, known to be missing from $\sD_{4}$ (see \cite{OEISD4}), leading one to ask what happens if 
%one allows 
a finite number of exceptions
is permitted. Indeed, Niederreiter \cite[p. 990]{Niederreiter1978} conjectured in $1978$ that
for $A=3$,
 $\sD_{3}$ already contains every sufficiently large number; we write this as
$$
\sD_{3}\supset \N_{\gg1}.
$$
With lots more computational capacity and evidence, Hensley almost $20$ years later \cite{Hensley1996} conjectured even more boldly that the same holds already for $A=2$:
\be\label{eq:D2N}
\sD_{2}\supset\N_{\gg1}
.
\ee
Lest the reader be tempted to one-up them all, let us consider the case $A=1$. %Obviously 
%It is elementary that
Here
$\sR_{1}$ contains only
continued fractions of the form $[1,\dots,1]$, and these are
 quotients of consecutive Fibonacci numbers $F_{n}$,
$$
\sR_{1}=\{F_{n}/F_{n+1}\}.
$$
%and
So $\sD_{1}=\{F_{n}\}$ is just the Fibonacci numbers, %so 
and
this is an exponentially thin sequence.
\\

In fact, Hensley conjectured something much stronger than \eqref{eq:D2N}.
%, as we 
%will  
%explain. 
First
% we need just a little 
 some more notation. Let $\sC_{A}$ be the set of limit points of $\sR_{A}$,
$$
\sC_{A}:=\{[a_{1},a_{2},\dots
%,a_{k},\dots
]:
%\text{ all }
a_{j}\le A,\
\forall j%\ge1
\}.
$$
This is a Cantor-like set with some Hausdorff dimension
\be\label{eq:gdZ}
\gd_{A}:=\dim(\sC_{A})
.
\ee
To get our bearings, consider again the case $A=1$. Then $\sC_{1}=\{1/\varphi\}$ is just the singleton consisting of the reciprocal of the golden mean, and %obviously 
hence
$\gd_{1}=0$.

%%%%%%%%%%%%%%%%%%%%%%%%%%%%%%%%%%%%%%%%%%%%%%%%%%%%
%%%%%%%%%%%%%%%%%%%%%%%%%%%%%%%%%%%%%%%%%%%%%%%%%%%%
\begin{figure}
%\begin{center}
\setlength{\unitlength}{0.12cm}
\begin{picture}(-50,10) %%%%%%%%%%%%%% C^{1}
%\put(-60,-2){$:$}
\put(-80,-2){$\sC_{2}^{(1)}=\{[a_{1},\dots]:a_{1}\le2\}$:}

\thinlines
%\dottedline
\put(-30,0){\line(1,0){60}}

\put(-30,1){\line(0,-1){2}}
\put(-31,-5){$0$}

\put(30,1){\line(0,-1){2}}
\put(29,-5){$1$}

%\put(-30,0){\dots\dots\dots}
%\thicklines

\put(-10,1){\line(0,-1){2}}
\put(-11,-5){$\frac13$}

\put(-9.5,3){${a_{1}=2\atop{\overbrace{\hskip1cm}^{}}}$}
\put(.5,3){${a_{1}=1\atop{\overbrace{\hskip1.35in}^{}}}$}

\put(0,1){\line(0,-1){2}}
\put(0,-5){$\frac12$}

\linethickness{.1cm}

\put(-10,0){\line(1,0){10}}

\put(0,0){\line(1,0){30}}

\end{picture}

\

\begin{picture}(-50,10) %%%%%%%%%%%%%%%%%%%%% C^{2}
%\put(-60,-2){$:$}
\put(-80,-2){$\sC_{2}^{(2)}=\{[a_{1},\dots]:a_{1},a_{2}\le2\}$:}

%\dottedline
%\put(-30,0){\line(1,0){20}}
%\put(-30,0){\dots\dots\dots}
%\thicklines

\thinlines
%\dottedline
\put(-30,0){\line(1,0){60}}

\put(-30,1){\line(0,-1){2}}
\put(-31,-5){$0$}

\put(30,1){\line(0,-1){2}}
\put(29,-5){$1$}

%\put(-30,0){\dots\dots\dots}
%\thicklines

\put(-10,1){\line(0,-1){2}}
\put(-12,-5){$\frac13$}

\put(-10,3){${a_{2}\le2\atop{\overbrace{\hskip0cm}^{}}}$}

\put(-5,1){\line(0,-1){2}}
\put(-7.5,-5){$\frac5{12}$}

\put(-3.3,1){\line(0,-1){2}}
\put(-3.7,-5){$\frac49$}

\put(20,1){\line(0,-1){2}}
\put(19,-5){$\frac56$}

\put(0,1){\line(0,-1){2}}
\put(0,-5){$\frac12$}

\put(0,3){${a_{2}\le2\atop{\overbrace{\hskip.9in}^{}}}$}

\put(15,1){\line(0,-1){2}}
\put(14,-5){$\frac34$}

\linethickness{.1cm}

\put(-10,0){\line(1,0){6.7}}

\put(0,0){\line(1,0){20}}

\end{picture}

\vskip.2in
\vdots
\vskip-.3in

\begin{picture}(-50,10)%%%%%%%%%%%%%%%%%%%% C^{infty}
%\put(-60,-2){$:$}
\put(-80,-2){$\sC_{2}^{(\infty)}=\sC_{2}$:}

\put(-30,-2.5){\includegraphics[width=7.5cm]{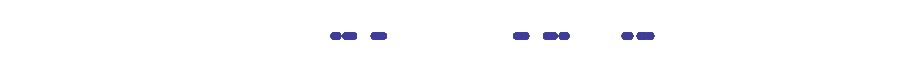}}

\thinlines
%\dottedline
\put(-30,0){\line(1,0){60}}

\put(-30,1){\line(0,-1){2}}
\put(-31,-5){$0$}

\put(30,1){\line(0,-1){2}}
\put(29,-5){$1$}

\end{picture}

\

%\end{center}
\caption{The Cantor set 
$
\sC_{2}=\bigcap\limits_{k=1}^{\infty}\sC_{2}^{(k)}
,
$ 
where 
$
\sC_{2}^{(k)}
=
\{
[a_{1},\dots,a_{j},\dots,a_{k},\dots]
:
a_{j}\le A
\text{ for all $1\le j\le k$}
\}
$
restricts only the first $k$ partial quotients.
}
\label{fig:C2}
\end{figure}
%%%%%%%%%%%%%%%%%%%%%%%%%%%%%%%%%%%%%%%%%
%%%%%%%%%%%%%%%%%%%%%%%%%%%%%%%%%%%%%%%%%%

Now %try 
take
$A=2$. Consider the unit interval $[0,1]$. The numbers in the range $(1/2,1]$ have first 
partial quotient 
%digit
$a_{1}=1$, and those in $(1/3,1/2]$ have first 
%digit 
partial quotient
$a_{1}=2$. The remaining interval $[0,1/3]$ has numbers whose
 %already has
  first 
%digit
partial quotient
 %that
  is 
  already
  too big, and thus is cut out. We repeat in this way, cutting out intervals for each 
 %digit
 partial quotient, and arriving at $\sC_{2}$; see Figure \ref{fig:C2}.
There is a substantial literature estimating the dimension $\gd_{2}$ which we will not survey, but the current record is due to Jenkinson-Pollicott \cite{JenkinsonPollicott2001}, whose superexponential algorithm estimates
\be\label{eq:gd2}
\gd_{2}
%\approx
=
0.5312805062772051416244686\dots
\ee
If we relax the bound $A$, the Cantor sets increase, as do their dimensions. In fact, Hensley \cite{Hensley1992} determined the asymptotic expansion, which to first order is
\be\label{eq:gdAbig}
\gd_{A}=1-{6\over \pi^{2}A}+o\left(\frac1A\right),
\ee
as $A\to\infty$.
% (His result is %actually 
%sharper, giving more lower order terms which we suppress.) 
In particular, the dimension can be made arbitrarily close to $1$ by taking $A$ large.
\\

We can now explain Hensley's stronger conjecture. His observation is that one need not only consider restricting the partial quotients $a_{j}$ to the full interval $[1,A]$; one can allow more flexibility by 
%picking 
fixing
any finite ``alphabet'' $\cA\subset\N$, and restricting the 
%digits
partial quotients to
the ``letters''
 %lie 
 in this alphabet. To this end, let $\sC_{\cA}$ be the Cantor set
$$
\sC_{\cA}:=\{[a_{1},a_{2},%\dots,a_{k},
\dots]: 
%\text{ all }
a_{j}\in\cA,\
\forall j\ge1
\}
,
$$
and similarly let $\sR_{\cA}$ be the partial convergents to $\sC_{\cA}$, $\sD_{\cA}$ the denominators of $\sR_{\cA}$, and $\gd_{\cA}$ the
Hausdorff
 dimension of $\sC_{\cA}$. 
%We can combine the two notations, using $\cA=A$ for the alphabet $\cA=\{1,\dots,A\}$. 
 Then Hensley's
 %elegant 
 %conjecture asserts 
 claim is
 the following % that
\begin{conj}[Hensley 1996 {\cite[Conjecture 3, p. 16]{Hensley1996}}]
\be\label{eq:gdcA}
\sD_{\cA}\supset\N_{\gg1}
\qquad\qquad
\Longleftrightarrow
\qquad\qquad
\gd_{\cA}>1/2
.
\ee
\end{conj}

Observe in particular that $\gd_{2}$ in \eqref{eq:gd2} exceeds $1/2$, and hence \eqref{eq:gdcA} implies that \eqref{eq:D2N} holds.
\\

%We have explained \eqref{eq:D2N} via, by why should the latter hold? 

\begin{figure}
        \begin{subfigure}[t]{0.46\textwidth}
                \centering
		\includegraphics[width=\textwidth]{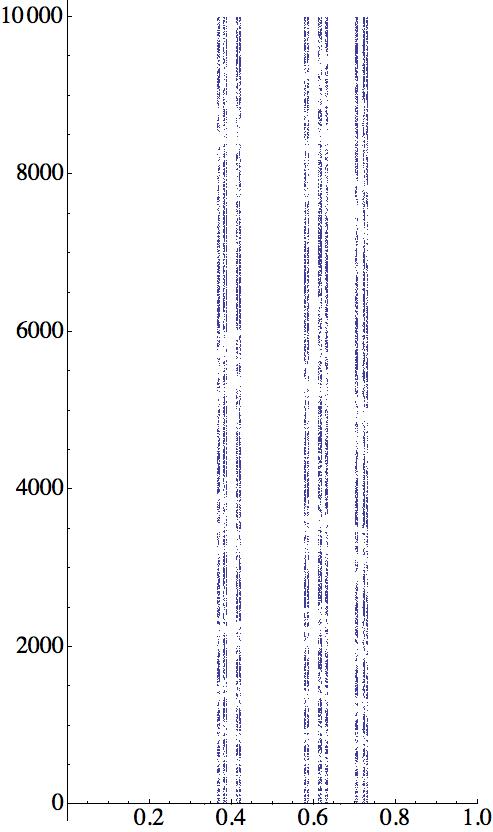}
                \caption{$A=2$, $N=10 000$.}
                \label{fig:RAN2}
        \end{subfigure}%
%\qquad
\qquad
        \begin{subfigure}[t]{0.46\textwidth}
                \centering
		\includegraphics[width=\textwidth]{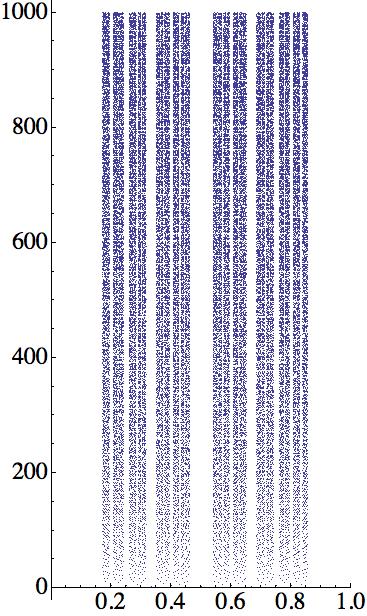}
                \caption{$A=5$, $N=1000$.}
                \label{fig:RAN5}
        \end{subfigure}

\caption{For each $b/d\in\sR_{A}(N)$, plot $b/d$ versus $d$, with $A$ and truncation parameter $N$ as shown.}
\label{fig:RAN}
\end{figure}

Here is some %more
heuristic
 evidence in favor of  \eqref{eq:gdcA}.
%To see things more clearly, l
Let us
  visualize the set $\sR_{\cA}$ of rationals, by grading each fraction according to the denominator. That is,
  %we show 
  plot each fraction $b/d$ at height $d$, showing the set
  % by plotting the 
  %points
%pairs
\be\label{eq:bdToD}
\left\{  \left(\ \frac bd\ ,\ d\right):
\ \frac bd\in\sR_{\cA}
,\quad
(b,d)=1
\right\}.
  \ee
We show this plot  in Figure \ref{fig:RAN2} for $\cA
%=A
=\{1,2\}$
%(that is, $\cA=\{1,2\}$)
 truncated at height $N=10 000$, and  in Figure \ref{fig:RAN5} for $\cA=%A=
 \{1,2,3,4,
 5\}$ truncated at height $N=1000$.
%
% truncate $\sR_{\cA}$ at height $N$, by which we mean defining
%Let us 
We 
give a name to this truncation, defining
$$
\sR_{\cA}(N):=
\left\{
\frac bd\in\sR_{\cA}
:
(b,d)=1,\ 1\le b<d<N
\right\}
.
$$
%This is still just a set of rationals, but we can

Observe first that the limiting vertical
 %tentacles
 directions in which the figures grow are precisely
  %are reaching up towards 
  the Cantor sets $\sC_{\cA}$; compare Figures \ref{fig:RAN2} and \ref{fig:C2}. Moreover, note that %, 
  if 
  %we have plotted 
  at least one point has been placed at height $d$, then $d\in\sD_{\cA}$. That is, the 
  %key issue boils down to the question: 
  ``beef'' 
  %here is:
  of this problem boils down to:
  what
  % happens 
  %if we 
  is the
  projection of the plots in Figure \ref{fig:RAN} %on
  to the $y$-axis? In particular, does every (sufficiently large) integer appear? % in such a projection?

%Before discussing the projection, we should first ascertain 
The first question to address is:
%the amount of ``mass'' at our disposal; that is, 
how big is $|\sR_{\cA}(N)|$, %or 
that is,
how many points are being plotted in Figures \ref{fig:RAN2} and \ref{fig:RAN5}?
Hensley
%
%\begin{thm}[Hensely 1989 
\cite{Hensley1989} %]
showed that,
as $N\to\infty$,
\be\label{eq:N2del}
\#\sR_{\cA}(N)
\quad
\asymp%_{\cA}
\quad
N^{2\gd_{\cA}}
,
\ee
where the implied constant can depend on $\cA$.
%\end{thm}
%
%This is proved via the thermodynamical formalism, %about which we will say more later, see \S\ref{sec:therm}.
%which we do not have sufficient space to exposit.
%
(Hensley proved this for the alphabet $\cA=\{1,2,\dots,A\}$, but the same proof works for
an
 arbitrary finite  $\cA$.)

Now, the $\Longrightarrow$ direction of \eqref{eq:gdcA} is trivial. Indeed, let
$$
\sD_{\cA}(N):=\sD_{\cA}\cap[1,N],
$$
so that the left hand side of \eqref{eq:gdcA} 
%can be rewritten as
is equivalent to
\be\label{eq:Hens2}
\#\sD_{\cA}(N)=N+O(1), \qquad\qquad\text{ as $N\to\infty$}
%\qquad%\qquad
%\Longleftrightarrow
%\qquad%\qquad
%\gd_{\cA}>1/2
.
\ee
Then it is clear that $\#\sR_{\cA}(N)$ counts $d$'s with multiplicity, whereas $\#\sD_{\cA}(N)$ counts each appearing $d$ only once; hence
\be\label{eq:Mult1}
\#\sD_{\cA}(N)
\le
\#\sR_{\cA}(N)
\overset{\eqref{eq:N2del}}{\ll}
N^{2\gd_{\cA}}
.
\ee
So if 
%the %left 
%assertion of
 \eqref{eq:Hens2} holds, then \eqref{eq:Mult1} implies that $2\gd_{\cA}$ must be at least $1$. 

%But %of course 
%the previous paragraph is essentially vacuous; 
A caveat: 
we do not know how to verify 
 \eqref{eq:Hens2} 
 %the left assertion of \eqref{eq:gdcA} 
for a single alphabet! Nevertheless the content of Hensley's Conjecture is clearly the opposite $\Longleftarrow$ direction. %Let us now present 
Here is some evidence in favor of this claim.
\\

%To this end, r
%Recall the following
%Here is a
An old % 1954 
theorem of 
Marstrand's \cite%[p. 267]
{Marstrand1954} %: Let 
states the following. %that i
%f
Let
$E\subset[0,1]\times[0,1]$ 
be 
%is
a Hausdorff measurable set having Hausdorff dimension $\ga>1$.
%Martrand's theorem says that 
%this projection
T%t
hen
%$$
%\text{
 the projection of $E$ 
into a line of slope $\tan\gt$
%$E$ 
%projects into a
is
 ``large,''
 for Lebesgue-almost every $\gt\in\R/2\pi\Z$.
 %, 
 %set %in almost 
%must be ``large''
%when collapsed 
%in almost every direction.
%}
%$$
%wh
Here
``large'' means
 %not only of full Hausdorff dimension, but 
 of positive Lebesgue measure.
 %,
%for almost every slope.
%That is,
%and ``almost every direction'' refers to
%for almost every angle $\gt$, the projection of $E$ into a line of slope $\gt$.
% given a 
%projective (i.e., through the origin) 
%line $L$,
 %(or what is the same, a slope in the unit circle)
% in the same plane, 
% let $\Proj_{L}(E)$
%denote the orthogonal projection of $E$ on to $L$. 
%Then 
%
%So by analogy, %and {\it only} by analogy, 
%we can consider
%Then
%to be a 
One 
may
thus
 heuristically
 %be 
% identify
 %i
 %ed 
 think of
 \eqref{eq:bdToD} 
%with 
as
$E$ above, with \eqref{eq:N2del} suggesting the ``dimension'' $\ga=2\gd_{\cA}$. Then $\sD_{\cA}$ is the projection of this $E$ to the $y$-axis, and
%, if this a generic direction, 
it should be ``large'' according to the analogy. 
%Of course 
Marstrand's theorem says nothing about an individual line, and %doesn't 
does not
apply to the countable set \eqref{eq:bdToD}, so the analogy cannot be furthered in any meaningful way. 
Nevertheless, we see the condition $\ga>1$ is converted into $2\gd_{\cA}>1$, giving evidence for the 
$\Longleftarrow$ direction of \eqref{eq:gdcA}.

For another heuristic, if one uniformly samples $N^{2\gd}$ pairs $(b,d)$ out of the integers up to $N$, a given $d$ is expected to appear with multiplicity roughly $N^{2\gd-1}$. For $\gd>1/2$ and $N$ growing, this multiplicity will be positive with probability tending to $1$.

This heuristic does not rule out the possible conspiracy that only very few (about $N^{2\gd-1}$) $d$'s actually appear, each with very high (about $N$) multiplicity. But %this %actually 
such an argument %gives
leads to
%gives 
another bit of evidence towards \eqref{eq:gdcA}: 
%; namely, 
since the 
multiplicity of any $d<N$ is at most $N$, we have the
%We give even some more evidence by proving an 
elementary lower bound % on $\#\sD_{\cA}(N)$. 
%Recall again that $\#\sR_{\cA}(N)$ counts each  $d$ with multiplicity, and for $d<N$ %the 
%this
%multiplicity is at most $N$.
%Hence we immediately have that
$$
\#\sD_{\cA}(N)\ge\frac1N\#\sR_{\cA}(N)
\overset{\eqref{eq:N2del}}{\gg}
\frac1N
N^{2\gd_{\cA}}
=
N^{2\gd_{\cA}-1}
.
$$
So if $\gd_{\cA}>1/2$, 
then
the %re are 
set $\sD_{\cA}$ %is 
already %infinite, and moreover 
grows at least at a %polynomial 
power
rate.
Furthermore, for any fixed $\vep>0$, one can take some $\cA=\cA(\vep)$ sufficiently large so that
$2\gd_{\cA}-1>1-\vep$.
 For example, 
  using 
 \eqref{eq:gdAbig}, we can
take $\cA=\{1,2,\dots,A\}$ where
 $$
 %\cA=
 A>{12\over\pi^{2}\vep}(1+o(1)).
 $$
Here $o(1)\to0$ as %with 
 $\vep\to0$.
 Hence one can produce $N^{1-\vep}$ points in $\sD_{\cA}(N)$, which is already 
 substantial progress towards \eqref{eq:Hens2}.
 %, one can make 
 %for any $\vep$, by choosing $\cA$ sufficiently large. 

But unfortunately, Hensley's conjecture \eqref{eq:gdcA}, as stated, is false.
%Indeed, in
\begin{lem}[Bourgain-K. 2011 
{\cite[Lemma 1.19]{BourgainKontorovich2011a}}%
]%\label{lem:BK}
%it was observed that
%For
%
 %t
The alphabet 
$\cA=\{2,4,6,8,10\}$ %, 
has
%the 
dimension $\gd_{\cA}=0.517\dots$, which exceeds $1/2$, but does not contain every sufficiently large number.
\end{lem}
\pf
%
%But % for this alphabet, 
%then 
%i
The dimension can be computed by the Jenkinson-Pollicott algorithm used to establish \eqref{eq:gd2}.
It is an elementary calculation %with
from the definitions to show for this alphabet that %if 
every fraction in $\sR_{\cA}$ is of the form $2m/(4n+1)$ or $(4n+1)/(2m)$,
%that is,
and so 
%$\sD_{\cA}$ is a subset of $2\Z\cup\{1\mod 4\}$.
%Hence
%.
%
%, while
 $\sD_{\cA}\equiv\{0,1,2\}(\mod4)$.
 %, and 
 Hence $\sD_{\cA}$ does not contain every sufficiently large number.
 \epf
%\end{lem}

That is, there can be congruence obstructions, in addition to the condition on dimension. This suggests instead a closer analogy with Hilbert's 11th problem on numbers represented by %integral (or rational) 
quadratic forms.
%, the solution being
%Hasse's local-global principle. % states that 
%Nevertheless, a version of Hensley's conjecture, suitably modified to account for local obstructions, may be plausible. 
%
%Following 
According to
this analogy, we make the following
\begin{Def} 
Call $d$ {\it represented} by the given alphabet $\cA$ if $d\in\sD_{\cA}$. 
Also, call $d$ {\it admissible} for the alphabet $\cA$ if it is everywhere locally represented, %in the sense 
meaning
that $d\in\sD_{\cA}(\mod q)$ for all $q\ge1$.
\end{Def}
One can then modify Hensley's conjecture to state  that, if $\gd_{\cA}$ exceeds $1/2$, then every sufficiently large admissible number is represented, akin to Hasse's local-to-global principle. 
\begin{rmk}\label{rmk:A12}
We will
 %show 
explain
 in \S\ref{sec:locZ} 
 %that for
 %why 
 that
 the alphabet $\cA=%A=2$ %
 \{1,2\}$
 %, there are 
 has no local obstructions, so \eqref{eq:D2N} 
is still plausible as it stands.
\end{rmk}

%Now we can present 
Here is
some %partial 
progress towards the %is 
conjecture.

\phantomsection
\begin{thmz}[Bourgain-K. 2011 \cite{BourgainKontorovich2011a}]\label{thm:Z}
Almost every natural number is the denominator of a reduced fraction whose partial quotients are
 %absolutely
  bounded by $50$.
\end{thmz}

Here ``almost every'' is in the sense of density: 
for $\cA=\{1,2,\dots,
%A=
50\}$,
%Let $\sA_{\cA}$ denote the set of admissible numbers for $\cA$; then the theorem states that 
$$
{
1
\over
N
%\#(\sA_{\cA}\cap[1,N])
}
\#(\sD_{\cA}\cap[1,N])
\to
1
,
$$
as $N\to\infty$.
%
%
%\begin{thm}[Bourgain-K. 2011 ]\label{thm:BK11}
The %actual statement is: % that 
proof in fact shows that
%there is a $\gd_{0}<1$ so that,
%F
for any alphabet $\cA$ having sufficiently large dimension 
\be\label{eq:gdMin}
\gd_{\cA}>
\gd_{0}
%1-5/312\approx0.98
,
\ee
almost every 
admissible number is represented, where
%\end{thm}
%
%The proof of
%Theorem \ref{thm:BK11} %says moreover that  one can take 
%gives 
%shows that
the value
\be\label{eq:gd0IsZ}
\gd_{0}=
1-5/312\approx0.98
\ee
is sufficient. % (of course this is far from the desired value $\gd_{0}=1/2$).
%In light of Hensley's sharper version of \eqref{eq:gdAbig}, it seems
 Using Hensley's  \eqref{eq:gdAbig}, the value $%\cA=
 A=50$ 
 seems to satisfy
  %suffice for 
 % have 
 % the condition 
\eqref{eq:gdMin}. 
%Moreover, %if 
%it is a 
The reason Theorem \hyperref[thm:Z]{Z} needs no mention of admissibility is that
%In fact, %that
any 
alphabet $\cA$
%for 
%satisfying
%\eqref{eq:gdMin} 
 %does not 
 with such a large dimension
 \eqref{eq:gd0IsZ}
 must
 contain
 both  $1$ and $2$;
 % the partial quotient 
%to be possible; an alphabet  
missing even one of these %two 
%partial quotients 
letters
will
 drop the 
%have 
%give
%result in
%too small a
 dimension by too much%
 . Hence there are actually no local obstructions in the theorem, cf. Remark \ref{rmk:A12}. 
%To compare more immediately to Zaremba's original Conjecture \hyperref[conj:Z]{Z}, we restate
%Theorem \ref{thm:BK11}
%which can be restated 
%in the following way.

To explain the source of this progress, we %now 
reformulate Zaremba's problem 
in a way that highlights the role 
%into one involving 
of 
the hitherto unmentioned
 ``thin orbit'' lurking underneath.
\\

\subsection{Reformulation}\

The key to the above progress is the old and elementary observation that
$$
\frac bd=[a_{1},\cdots,a_{k}]
$$
is equivalent to
\be\label{eq:mats}
\mattwo*b*d =
 \mattwo011{a_{1}}
\cdots 
 \mattwo011{a_{k}}
.
\ee
With this observation, it is natural to introduce the semigroup generated by matrices of the above form with partial quotients restricted to the %allowed 
given
alphabet. Let
\be\label{eq:Ggen}
\G=
\G_{\cA}:=
\<
\mattwo011a:a\in\cA
\>^{+}
,
\ee
where the ``$+$'' denotes generation as a semigroup (no inverse %s
matrices).
Then 
%clearly 
the orbit 
\be\label{eq:cQZ}
\cO=\cO_{\cA}:=\G%_{\cA}
\cdot\bv_{0}
\ee
with 
\be\label{eq:bv0Z}
\bv_{0}=(0,1)^{t}
\ee
%consists of 
isolates
the set of second columns in $\G%_{\cA}
$, and
from \eqref{eq:mats}
 is hence in bijection with the set $\sR_{\cA}$.
The ``thinness''  of the orbit is explained by \eqref{eq:N2del}, which implies that
$$
\#\{\bv\in\cO:\|\bv\|<N\}\asymp N^{2\gd_{\cA}}
,
$$
as $N\to\infty$. %Of course i
If $\cO$ consisted of all %coprime 
integer
pairs $(b,d)^{t}$, the above count would be replaced by $N^{2}$, ignoring constants.
So this is the reason we call $\cO$ thin: it contains many fewer points than the ambient set in which it naturally sits.
 
From \eqref{eq:mats} again, the set $\sD_{\cA}$ is %just 
nothing more than
the set of bottom right entries of matrices in $\G_{\cA}$. This can be isolated via: 
\be\label{eq:DAvGv}
\<\bv_{0},\cO\>
=
\<\bv_{0},\G\cdot\bv_{0}\>
=
\sD_{\cA}
,
\ee
where the inner product is the standard one on $\R^{2}$.
Thus $d$ is represented if and only if there is a $\g\in\G$ so that 
\be\label{eq:innProd}
d=\<\bv_{0},\g\cdot\bv_{0}\>
,
\ee
with $\bv_{0}
$ 
given in \eqref{eq:bv0Z}.

\subsection{Local Obstructions}\label{sec:locZ}\

One can now easily understand
 Remark \ref{rmk:A12}, and
 the source of %the 
 any potential
 local obstructions. % in Lemma \ref{lem:BK}.
 %, and
 The key observation, via \eqref{eq:DAvGv}, is that to understand $\sD_{\cA}(\mod q)$, one needs only to understand the reduction of $\G(\mod q)$. And the latter can be analyzed by some algebra, namely the so-called {\it strong approximation} property; see e.g. \cite{Rapinchuk2012} for a comprehensive survery.
 This is 
 %a type of Chinese Remainder Theorem 
a property which determines when the reduction mod $q$ map is onto. For  general algebraic groups this is a deep theory, the first %complete 
proof \cite{MatthewsVasersteinWeisfeiler1984} using the classification of finite simple groups. But for $\SL_{2}$, the proofs are elementary, see e.g. \cite{DavidoffSarnakValette2003}. 

First observe that $\G$ sits inside the integer points of the {\it algebraic} group $\GL_{2}$, %(\Z)$, 
meaning that any %integer 
solution in $\Z$ to %a certain %set of 
the
polynomial 
equation %s
%, namely
 $(ad-bc)m=1$ gives an element %in 
$
  \bigl( \begin{smallmatrix}
a&b\\ c&d
\end{smallmatrix} \bigr)
\in
\GL_{2}(\Z)$, and vice-versa. Actually $\GL_{2}$
does not have strong  approximation, 
(e.g. the determinant in $\GL_{2}(\Z)$ can only be $\pm1$, while in $\GL_{2}(\Z/5\Z)$ it is $1,2,3$ or $4$, so the reduction map is not onto). 
S%
%s
o we first pass to $\SL_{2}$, as follows. % which does (see below). %
 %have strong approximation
%That is,
%t
The generators in \eqref{eq:Ggen} all have determinant $-1$, so the product of any two has determinant $+1$.
We make these products the generators for a subsemigroup $\tilde\G$ of $\G$, that is,
set
% replace
%%ing 
%$\G$ by 
$\tilde\G:=\G\cap\SL_{2}$.
%,
%%. That is, 
%or in other words,
%take only even words in the generators (which all have determinant $-1$) in \eqref{eq:Ggen}.  
We recover the original $\G$-orbit $\cO$ in \eqref{eq:cQZ} by a finite union of $\tilde\G$-orbits.
The limiting Cantor set and its Hausdorff dimension are unaffected.
% in an obvious way. % such. 
%Since every rational number has two continued fraction expansions, 

Then strong approximation says essentially that for $p$ a sufficiently large prime, and $q=p^{e}$ any %prime 
$p$
power, the reduction of $\tilde\G$ mod $q$ is all of $\SL_{2}(\Z/q\Z)$.
(%Note that i
It does not matter that $\tilde\G$ is %not itself a 
only a semigroup; 
%its
upon
reduction mod $q$, %makes it into 
%results in
%one
%is
it
becomes
%one
a group.)
 Moreover for {\it ramified} primes $p$ (those for which the reduction mod $p$ is not onto), the reduction mod sufficiently large powers of $p$ {\it stabilizes} after some finite power. This means that
 there is some $e_{0}=e_{0}(p,\tilde\G)$ so that the following holds.
For any $e>e_{0}$, 
   if % the reduction of 
   $M\in\SL_{2}(\Z/p^{e}\Z)$ %such that 
   is such that its reduction is in 
   %mod $p^{e_{0}}$ is in 
   $\tilde\G(\mod p^{e_{0}})$, then $M$ is also in 
%the reduction 
$\tilde\G(\mod p^{e})$.
% is the full preimage in $\SL_{2}(\Z/p^{e}\Z)$ of $\tilde\G(\mod p^{e_{0}})$ under the natural projection 
% $$
% \Z/p^{e}\Z\to\Z/p^{e_{0}}\Z.
% $$ 
(These statements are best made in the language of $p$-adic numbers, which we avoid here.)
A key ingredient %here 
is that, while $\tilde\G$ is some strange subset of $\SL_{2}(\Z)$, it is nevertheless {\it Zariski dense} in $\SL_{2}$. 
This  means that 
if $P(a,b,c,d)$ is a polynomial
 %in the entries of a $2\times2$ matrix $\g$
  which vanishes for every $
  \bigl( \begin{smallmatrix}
a&b\\ c&d
\end{smallmatrix} \bigr)
  %\g
  \in\tilde\G$, then $P$ also vanishes on all matrices in $\SL_{2}$ with entries in $\C$.

In the above, ``sufficiently large,'' both for primes $p$ to be unramified, and 
the
stabilizing powers $e_{0}$ of ramified primes, 
can be effectively computed in terms of the generators.
Then for an arbitrary modulus $q=p_{1}^{e_{1}}\cdots p_{k}^{e_{k}}$, the reduction mod $q$ can be pieced  together from those mod $p_{j}^{e_{j}}$ using 
a type of Chinese Remainder Theorem 
%that can be 
%when
%applied to 
%algebraic 
groups called Goursat's Lemma. %. The lemma states that 
%(%A key ingredient 
%We use
%here is that $\PSL_{2}(\F)=\SL_{2}(\F)/\{\pm I\}$ is simple, for any finite field $\F$ 
%a finite field 
%of cardinality at least $5$.)
This leaves some %elementary 
finite group theory to determine completely the reduction of $\tilde\G$ mod any $q$, and hence explains all local obstructions via \eqref{eq:DAvGv}.
\\

%\newpage

%With this reformulation, w
We now leave Zaremba's problem, and return to sketch a proof of Theorem \hyperref[thm:Z]{Z} in \S\ref{sec:C}.

\newpage

\section{Integral Apollonian Gaskets}\label{sec:A}

\begin{figure}
        \begin{subfigure}[t]{0.3\textwidth}
\includegraphics[width=\textwidth]{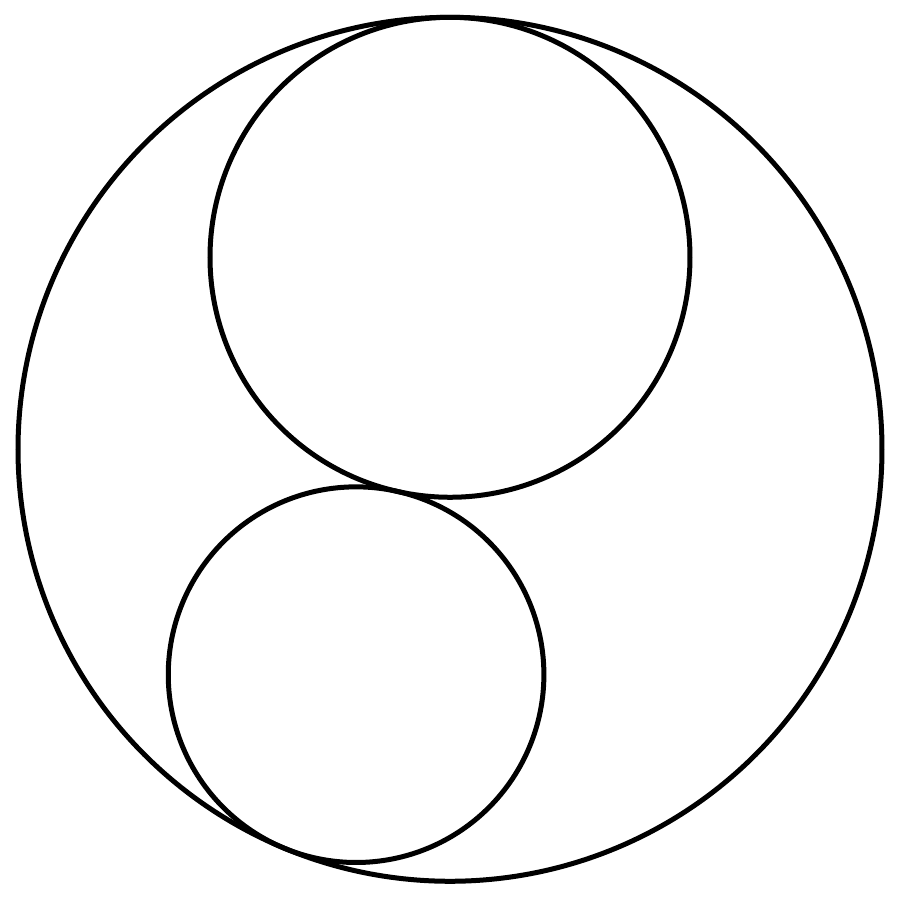}
                \caption{Three mutually 
                tangent circles}
                \label{fig:ApGen0}
        \end{subfigure}%
\quad
        \begin{subfigure}[t]{0.3\textwidth}
\includegraphics[width=\textwidth]{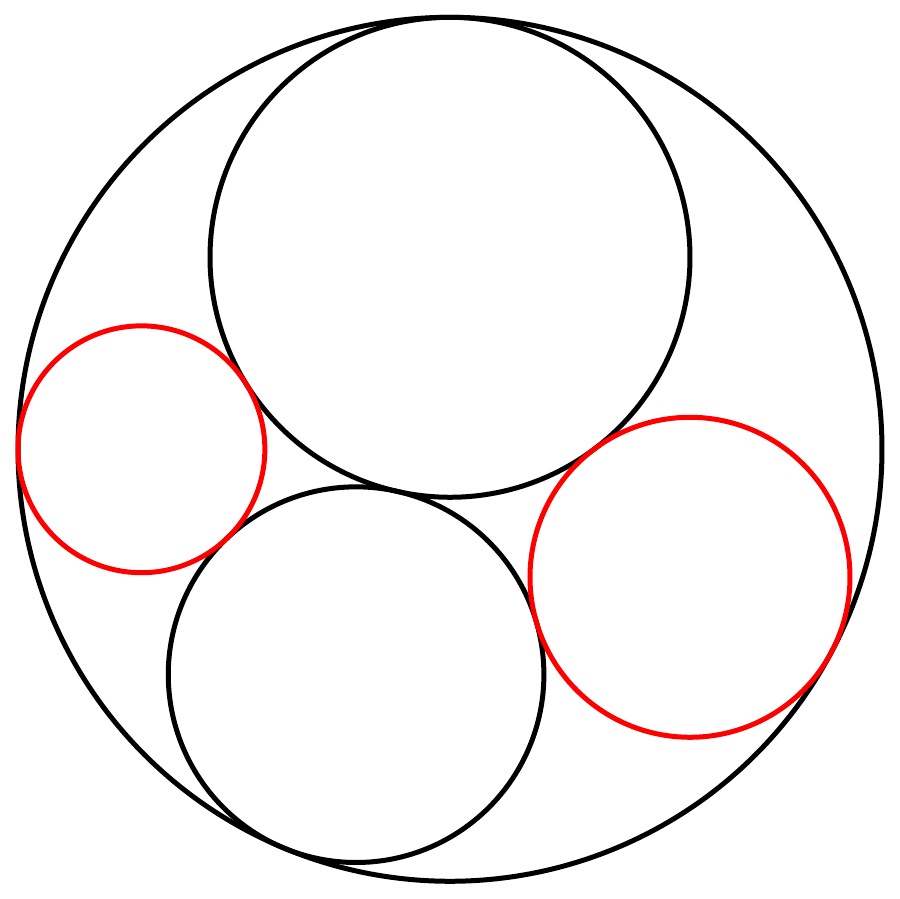}
                \caption{Two more tangent circles}
                \label{fig:ApGen1}
        \end{subfigure}
\quad
        \begin{subfigure}[t]{0.3\textwidth}
\includegraphics[width=\textwidth]{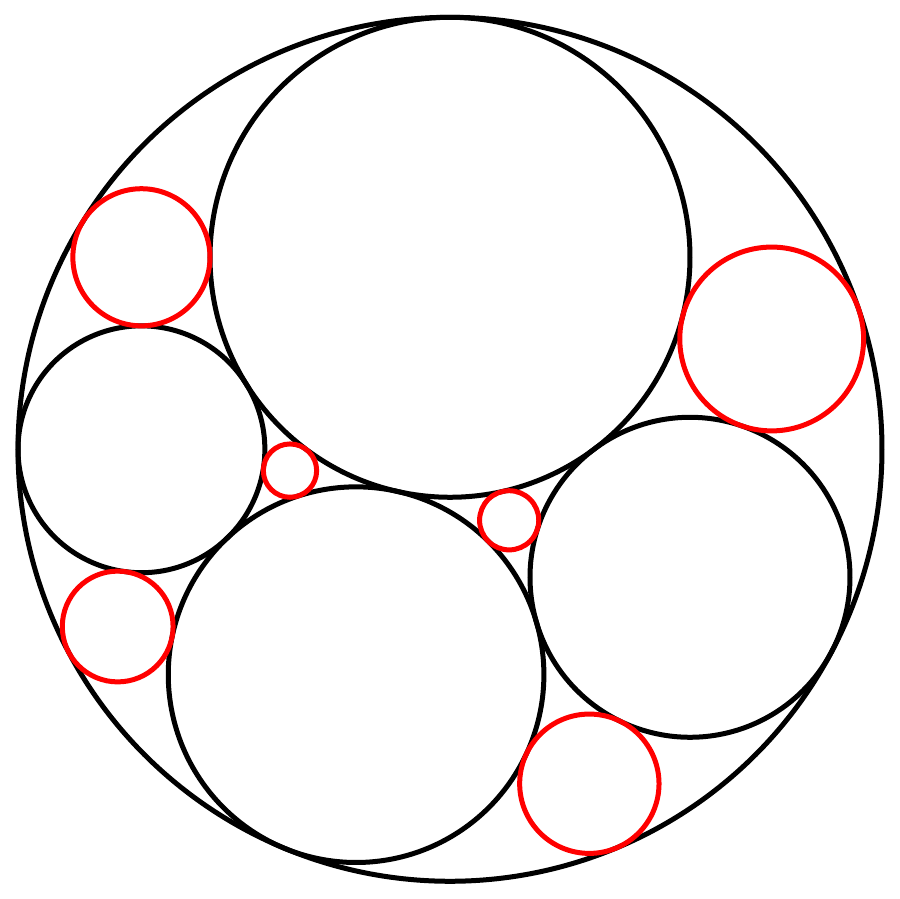}
                \caption{Six more tangent circles}
                \label{fig:ApGen2}
        \end{subfigure}

\caption{Tangent circles}
\label{fig:ApGen}
\end{figure}

Apollonius of Perga (ca 262-190 BC) wrote a two-volume book on Tangencies, solving in every conceivable configuration
the following general problem: Given three circles in the plane, any of which may have radius zero (a point) or infinity (a line), construct a circle tangent to the given ones. The volumes were lost but the statements survived via a survey of the work by Pappus. %; as late as the 15th century, some mathematicians doubted whether the problem was solvable 
In the special case when the given three circles are themselves mutually tangent with disjoint points of tangency  (Figure \ref{fig:ApGen0}), Apollonius proved that 
\be\label{eq:Apoll}
\text{
there are exactly two solutions}
\ee 
to his problem (Figure \ref{fig:ApGen1}). Adding these new circles to the configuration, one has many other triples of tangent circles, and Apollonius's construction can be applied to them (Figure \ref{fig:ApGen2}). Iterating in this way {\it ad infinitum}, as apparently was first done in Leibniz's  notebook, gives rise to a circle packing,
the closure of which
 has become known 
 in the last century
 as 
% is now called
 an {\it Apollonian gasket}.
% (As far as the author can tell, this term was coined in \cite{KasnerSupnick1943}.)
We restrict our discussion henceforth to bounded gaskets, such as that 
  %One such is
   illustrated in Figure \ref{fig:Apoll}; 
%where 
there
the number shown inside a circle is its curvature, that is, one over its radius.
Such pictures have received considerable attention recently,
see e.g. \cite{LagariasMallowsWilks2002, GrahamLagarias2003, GrahamLagariasMallowsWilksYanI, GLMWYII, GLMWYIII, ELag2007,
SarnakToLagarias, Sarnak2008a, BourgainGamburdSarnak2010, KontorovichOh2011, Oh2010,  BourgainFuchs2011, Sarnak2011, Fuchs2011, FuchsSanden2011, OhShah2012, Vinogradov2012, LeeOh2012, BourgainKontorovich2012}.
We
 %briefly 
 %now survey some of these results, 
 will
 focus %ing 
 our discussion
 on the following two problems:
 
 \begin{enumerate}
\item
The Counting Problem: 
%The first question one might ask is:
For  a fixed %Apollonian 
gasket $\sG$,
 how quickly do the circles shrink, or alternatively, how many circles are there in $\sG$ with curvature bounded by a growing parameter $T$? 
 
 \item
 The Local-Global Problem: Suppose $\sG$ is furthermore {\it integral}, meaning that the curvatures of all circles in it are integers, such as the gasket in Figure \ref{fig:Apoll}. How many distinct integers appear up to a growing parameter $N$? That is, count curvatures up  to $N$, but without multiplicity.
  \end{enumerate}
  
  Problem (2) does not yet look like a local-global question, but will soon turn into one. We first address Problem (1) in more detail.
\\

\subsection{The Counting Problem}

\subsubsection{Preliminaries}\

 Some notation: for a typical circle $C$ in a fixed bounded gasket $\sG$, let $r(C)$ be its radius and 
$$
b(C)=1/r(C)
$$ 
its curvature (or ``bend''). Let
\be\label{eq:NGT}
\cN_{\sG}(T):=\#\{C\in\sG:b(C)<T\}
\ee
be the desired counting function.
To study this quantity, one might introduce an ``$L$-function'':
\be\label{eq:Lfun}
\cL_{\sG}(s):=
\sum_{C\in\sG}
{1\over b(C)^{s}}
=
\sum_{C\in\sG}
{r(C)^{s}}
.
\ee
Since the sum of the areas of inside circles in $\sG$ yields the
% whole of $\sG$
area of the bounding circle,
% in $\sG$ is bounded. 
%T
%his
the  series 
$\cL_{\sG}$
converges for $\Re(s)\ge2$.
%, s
%
It has some
{\it abscissa of convergence} $\gd$,
meaning $\cL_{\sG}$ converges for $\Re(s)>\gd$ and diverges for $\Re(s)<\gd$.
% and i
 It % is not hard to convince oneself
 follows
  from \eqref{eq:Lfun}
that $\gd$ is %in fact 
%actually
the Hausdorff dimension of the gasket $\sG$ \cite{Boyd1973a}.
In fact, Apollonian gaskets are 
 rigid, in the sense that
one can be %sent 
mapped
to any other by 
M\"obius transformations.
The latter are conformal (angle preserving) motions of the complex plane, %transforming under
sending
 $z\mapsto(az+b)/(cz+d)$, 
$
ad-bc=1
$%
.
%
%\newpage
%
% move one to any other, and h
%
%\newpage
%
Hence $\gd$ is a universal constant; %e.g. 
McMullen \cite{McMullen1998} estimates  that 
\be\label{eq:gdIs}
\gd=1.30568\dots
\ee
From such considerations,
% and standard so-called ``Tauberian'' arguments, 
%it is not hard to conclude, as did
 Boyd \cite{Boyd1982} was able to conclude
that
$$
{\log\cN_{\sG}(T)\over \log T}
\to
\gd,
$$
as $T\to\infty$. 

To refine this crude estimate to an asymptotic formula for $\cN_{\sG}(T)$, 
%one might try to
%
% show that $\cL_{\sG}(s)$ has analytic continuation beyond the region of absolute convergence, and take a residue at the %pole $s=\gd$. One approach towards this might be to 
% 
the author and Oh
 \cite{KontorovichOh2011}
 %find 
% found
established
 a ``spectral interpretation'' for $\cL_{\sG}$, %and exactly this is achieved in, 
 proving:
\be\label{eq:KO}
\cN_{\sG}(T)\sim c%_{\sG}
\cdot T^{\gd}, 
\ee
for some $c=c({\sG})>0$, as $T\to\infty$. 
%In fact t
(This asymptotic was recently   refined further
% independently and simultaneously %by 
in
Vinogradov's thesis \cite{Vinogradov2012} and 
independently
by Lee-Oh \cite{LeeOh2012}, giving %an effective power savings 
%sharp 
lower order
error terms.)
% We 
%now 
%briefly 
%will spend the rest 
The remainder
of this subsection
is devoted to
explaining 
this 
spectral interpretation and 
%giving 
highlighting
some of the ideas
 %towards 
going into the  proof of \eqref{eq:KO}.
\\

\subsubsection{Root quadruples and generation by reflection}\

It is easy to see \cite
[p. 14]
{GrahamLagarias2003}
that
%E
each such gasket $\sG$ contains a {\it root configuration} $\cC=\cC(\sG):=(C_{1},C_{2},C_{3},C_{4})$ of four largest mutually tangent circles in $\sG$.
 Let 
\be\label{eq:rootQ}
\bv_{0}=\bv_{0}(\sG)=(b_{1},b_{2},b_{3},b_{4})^{t}
\ee 
with $b_{j}=b(C_{j})$ be the {\it root quadruple} of corresponding curvatures.
The bounding circle, being internally tangent to the others, is given opposite orientation
%; orienting 
to make all interiors disjoint; %, the bounding circle is given 
%and
this is accounted for by giving it
negative curvature. 
For example in Figure \ref{fig:Apoll}, the root quadruple is 
\be\label{eq:bv0G}
\bv_{0}=(-10,18,23,27)^{t},
\ee 
where the bounding circle has radius $1/10$.
%The root quadruple 
%for such a bounded packing 
%satisfies
%$b_{1}+b_{2}+b_{3}+b_{4}>0$, $b_{1}<0<b_{2}\le b_{3}\le b_{4}$, and $b_{1}+b_{2}+b_{3}\ge b_{4}$

\begin{figure}
		\includegraphics[height=2in]{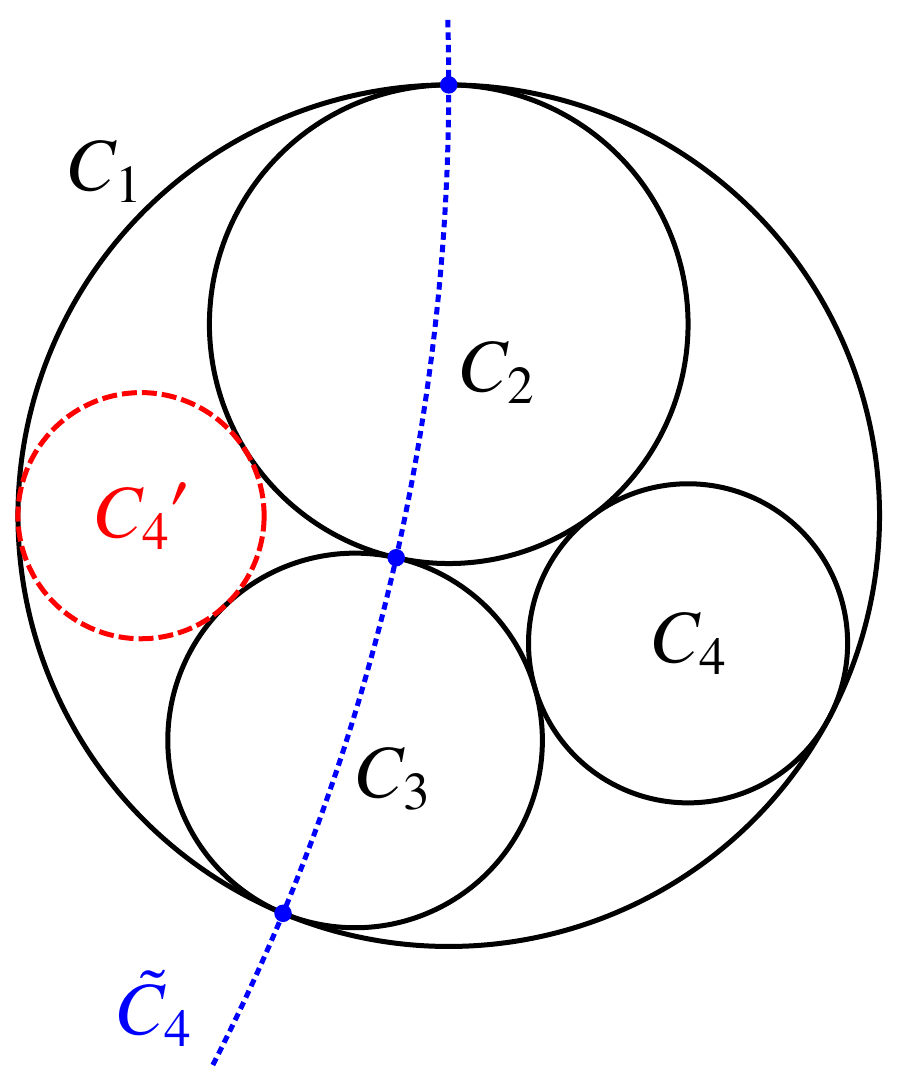}
                \caption{Generation from a root quadruple}
                \label{fig:Dual}
\end{figure}

Three tangent circles, say $C_{1},C_{2},C_{3}$ have three points of tangency, and determine a {\it dual circle} $\tilde C_{4}$ passing through these points, see Figure \ref{fig:Dual}. Thus the root configuration $\cC$ determines a {\it dual configuration} $\tilde\cC=(\tilde C_{1},\tilde C_{2},\tilde C_{3},\tilde C_{4})$ of four mutually tangent circles, orthogonal to those in $\cC$, see Figure \ref{fig:Duals}.
Reflection 
%in the plane 
through $\tilde C_{4}$ fixes $C_{1},C_{2},$ and  $C_{3}$, and sends $C_{4}$ to $C_{4}'$, the other solution to Apollonius's problem \eqref{eq:Apoll}, see Figure \ref{fig:Dual}. %In this way, r
Starting with the root configuration, repeated reflections through the dual circles give the whole circle packing.
\\

\begin{figure}
%        \begin{subfigure}[t]{0.32\textwidth}
\centering
\includegraphics[width=.5\textwidth]{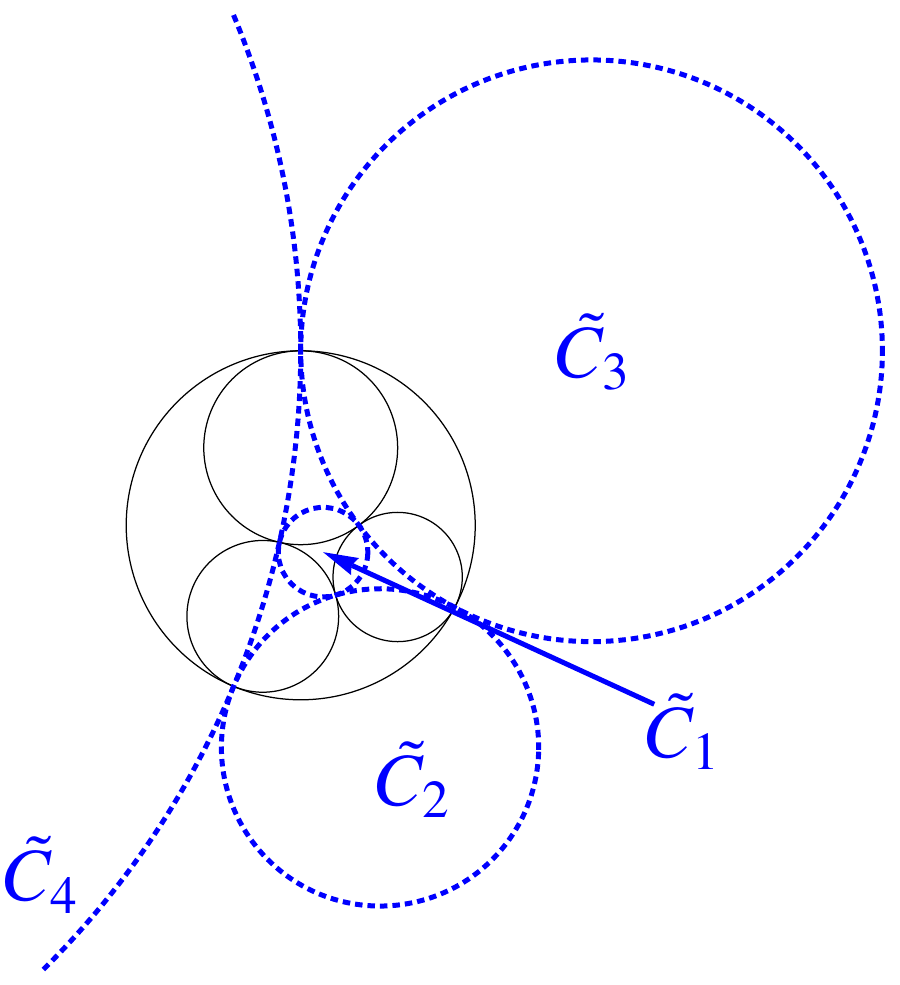}

                \caption{Root and dual configurations}
                \label{fig:Duals}
        \end{%sub
        figure}%

%Descartes theorem, integrality

\subsubsection{%Poincar\'e extension 
Hyperbolic space
and the group $\cA$}\

%We take the 
Following Poincar\'e, 
we
% take the 
extend %sion 
%move 
these
circle 
reflections
to the hyperbolic upper half space, 
\be\label{eq:bH3}
\bH^{3}:=\{(x_{1},x_{2},y):x_{1},x_{2}\in\R,y>0\}
,
\ee
replacing the action of the dual circle $\tilde C_{j}$ by a reflection
 through a (hemi)sphere
 $\fs_{j}$ 
  whose equator is $\tilde C_{j}$ (with $j=1,\dots,4$). We abuse notation, writing $\fs_{j}$ for both the hemisphere and the 
conformal map reflecting through $\fs_{j}$.  
The %so-called {\it Apollonian group}
group
\be\label{eq:Gapp}
\cA:=\<
\fs_{1},
\fs_{2},
\fs_{3},
\fs_{4}
\>
<
\Isom(\bH^{3})
,
\ee
generated by these reflections % through these four hemispheres
acts %ing 
discretely on $\bH^{3}$;
it is
%as 
a
 so-called {\it Schottky} group,
in that 
the four
generating
 spheres have disjoint interiors.
 %, and hence reflections through them 
%generate a
%
%

 The $\cA$-orbit of any fixed base point $p_{0}\in\bH^{3}$ has a limit set in the boundary $\dd\bH^{3}$, which is easily seen to be the original gasket, see Figure \ref{fig:ApOrb}.
 A {\it fundamental domain}
 for an %discrete 
 action 
 is a region 
 \be\label{eq:fundDom}
 \gW\subset\bH^{3}
 \ee 
 so that any point in $\bH^{3}$ can be sent to $\gW$ in 
 an essentially unique way;
% 
% \newpage
% 
% 
  for %  this action 
  the action of $\cA$, one can take $\gW$ to be
%being 
%is
%then
the exterior of the four %dual 
%(geodesic) 
hemispheres.
%, these hemispheres are totally geodesic.
%
 %(and hence $\G$ is geometrically finite). 
 %Two facts are evident from 
 %
 To see this, observe that if a point $p=(x_{1},x_{2},y)\in\bH^{3}$ is inside one of the spheres $\fs_{j}$, then its reflection $\fs_{j}(p)$ is outside of $\fs_{j}$ and has a strictly larger $y$-value. This does not guarantee that $\fs_{j}(p)$ is 
 outside all of the other spheres, but if it is inside some $\fs_{k}$, then
 %not inside one of the other spheres $\fs_{k}$, but its 
 reflection through $\fs_{k}$ will again have even higher $y$-value.
 This procedure must halt after finitely many iterations, since
 the only limit points of $\cA$ are in the boundary $\dd\bH^{3}$ where $y=0$. And it halts only when the image is outside of the four geodesic hemispheres. %To see uniqueness, observe that any point in $\gW$ 
 Uniqueness follows since any reflection $\fs_{j}$ takes a point in $\gW$ to a point inside $\fs_{j}$, that is, not in $\gW$.

 Two facts are evident from %the %domain
the above: first of all, $\cA$ is {\it geometrically finite}, meaning it has a fundamental domain bounded by a finite number
 (here it is four) of geodesic\footnote{%
%Recall a 
A geodesic in hyperbolic space is a straight vertical line or a  semicircle orthogonal to the boundary $\dd\bH^{3}$.}
 hemispheres; on the other hand, $\cA$ has {\it infinite} co-volume, that is, 
 any fundamental domain
 %$\cA\bk\bH^{3}$ 
 has infinite  volume with respect to the hyperbolic % volume form 
 measure
 $$
 y^{-3}dx_{1}dx_{2}dy
$$ 
 in the coordinates \eqref{eq:bH3}.
Note %also 
moreover
that $\cA$ has the
 structure %is 
%easily seen to be 
%that 
of a % so-called 
{\it Coxeter group}, % meaning it is 
being
free save the relations $\fs_{j}^{2}=I$ for the generators.
%Hence it is the ``alphabet'' which  generates our Apollonian gasket; words in the generating ``letters'' $\fs_{j}$ . 
It is also %clearly 
the symmetry group of all M\"obius transformations fixing $\sG$.
\\

%\qquad
        \begin{%sub
        figure}%[t]{0.65\textwidth}
                \centering
\includegraphics[width=\textwidth]{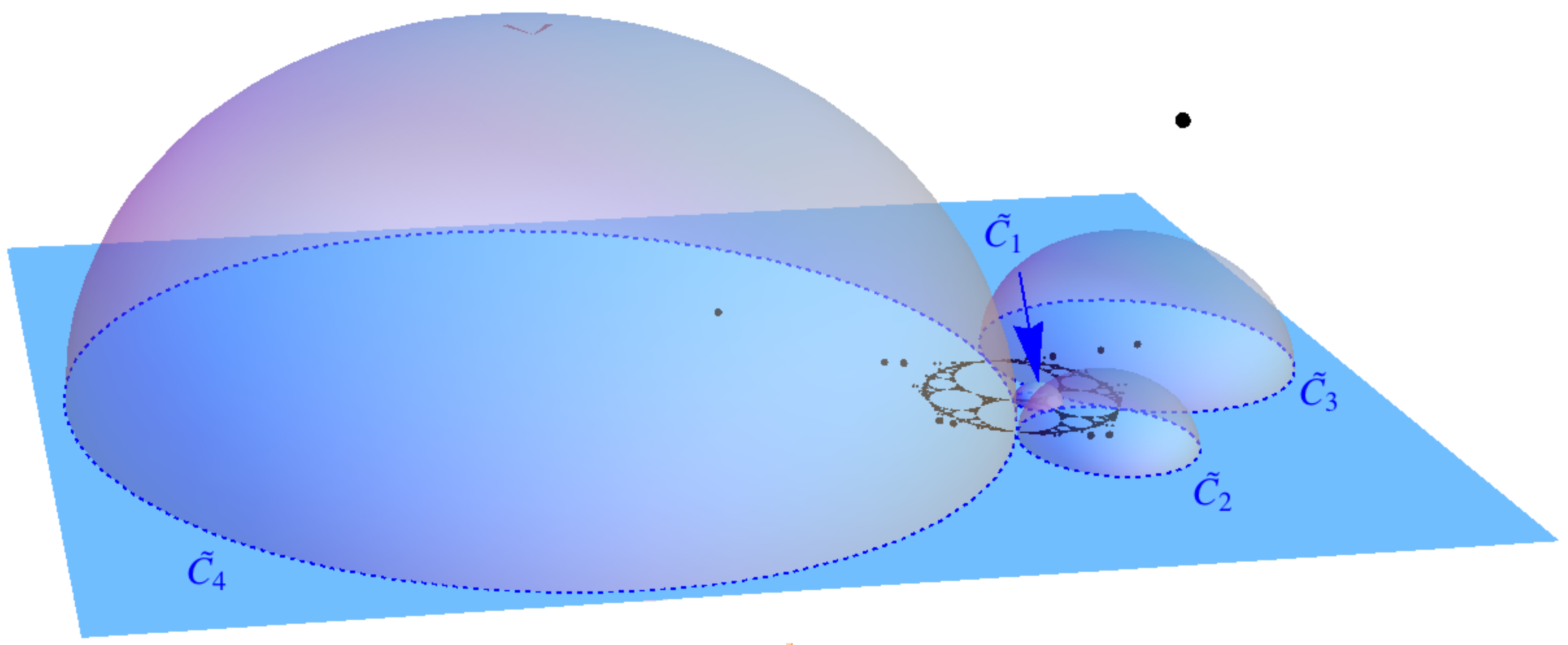}

                \caption{Poincar\'e extension: an $\cA$-orbit in $\bH^{3}$}
                \label{fig:ApOrb}
        \end{%sub
        figure}

%\caption{Poincar\'e extension}
%\label{fig:Poin}
%\end{figure}

\subsubsection{Descartes' Circle Theorem
%,
 %the Apollonian group,
  and integral gaskets}\

Next we need an
%other 
observation due to Descartes in the year 1643 \cite[pp. 37-50]{Descartes1901} (though his proof had a gap \cite{Coxeter1968}), that a quadruple $\bv=(b_{1},b_{2},b_{3},b_{4})^{t}$ of signed curvatures of four mutually tangent circles lies on the cone
\be\label{eq:cone}
Q(\bv)=0
,
\ee
where $Q$ is the so-called ``Descartes quadratic form''
\be\label{eq:QDes}
Q(\bv):=
2\left(
b_{1}^{2}+
b_{2}^{2}+
b_{3}^{2}+
b_{4}^{2}
\right)
-
\left(
b_{1}+
b_{2}+
b_{3}+
b_{4}
\right)^{2}
.
\ee
By a real linear change of variables, $Q$ can be diagonalized to the form 
$$
x^{2}+y^{2}+z^{2}-w^{2},
$$ 
that is,
it has
%of 
signature $(3,1)$. % over 
%the real numbers. %$\R$. 
Arguably the most beautiful formulation of 
%this
 %result 
%fact 
Descartes' Theorem
(rediscovered on many separate occasions% since Descartes
) is the following excerpt  from Soddy's
1936 {\it Nature} poem \cite{Soddy1936}:
\\

{\scriptsize
\quad Four circles to the kissing come. / The smaller are the bender. /

\quad The bend is just the inverse of /  The distance from the center. /

\quad Though their intrigue left Euclid dumb /  There's now no need for rule of thumb. /

\quad Since zero bend's a dead straight line / And concave bends have minus sign, /

\quad The sum of the squares of all four bends / Is half the square of their sum.\\
}

If $b_{1},b_{2}
$ and $b_{3}$
are given, then \eqref{eq:cone} is a quadratic equation in $b_{4}$ with two solutions, $b_{4}$ and $b_{4}'$, say; this is an algebraic proof of Apollonius's theorem \eqref{eq:Apoll}. It is then an elementary exercise to see that
$$
b_{4}+b_{4}'
=
2(b_{1}+b_{2}+b_{3})
.
$$
In other words, if the quadruple $(b_{1},b_{2},b_{3},b_{4})^{t}$ is given, then one obtains the quadruple 
%$(b_{1},b_{2},b_{3},b_{4}')$
with $b_{4}$ replaced by $b_{4}'$ via a linear action:
$$
\bp
1&&&\\
&1&&\\
&&1&\\
2&2&2&-1
\ep
\cdot
\bp
b_{1}\\
b_{2}\\
b_{3}\\
b_{4}
\ep
=
\bp
b_{1}\\
b_{2}\\
b_{3}\\
b_{4}'
\ep
.
$$
Hence we have given an algebraic realization to the geometric action of $\tilde C_{4}$ (or $\fs_{4}$) on the root quadruple, see again Figure \ref{fig:Dual}. Call the above $4\times4$ matrix $S_{4}$. 
Of course one could also send other $b_{j}$ to $b_{j}'$ keeping the three complementary curvatures fixed, via the matrices
\be\label{eq:SjAp}
S_{1}
=
\bp
-1&2&2&2\\
&1&&\\
&&1&\\
&&&1
\ep
,
S_{2}
=
\bp
1&&&\\
2&-1&2&2\\
&&1&\\
&&&1
\ep
,
S_{3}
=
\bp
1&&&\\
&1&&\\
2&2&-1&2\\
&&&1
\ep
.
\ee
Moreover one can iterate these actions, so we introduce the so-called {\it Apollonian group}
$\G$, isomorphic to $\cA$, generated by the $S_{j}$:
\be\label{eq:GamAp}
\G:=
\<
S_{1},
S_{2},
S_{3},
S_{4}
\>
.
\ee
Then the orbit 
\be\label{eq:cOAp}
\cO:=\G\cdot\bv_{0}
\ee
of the root quadruple $\bv_{0}$ under the Apollonian group $\G$ 
consists of all quadruples corresponding to curvatures of four mutually tangent circles in the gasket $\sG$.
%This 
We can now
explain %s 
the integrality of all curvatures in Figure \ref{fig:Apoll}: the group $\G$ has only integer matrices, so if the root quadruple $\bv_{0}$ (or for that matter any % quadruple of 
curvatures of 
four mutually tangent circles in $\sG$) is integral, then all curvatures in $\sG$ are integers! 
This fact seems to have been first observed by Soddy \cite{Soddy1937}.
\\

\subsubsection{Reformulating the counting statement, and the thin orbit}\

Moreover, note that starting with $\bv_{0}$, any new circle generated by a reflection is the smallest in its configuration, and hence has largest curvature. That is, for $\bv=\g\cdot\bv_{0}\in\cO$, writing $\g\in\G$ as a reduced word in the generators $\g=S_{i_{k}}\cdots S_{i_{1}}$, the last multiplication by $S_{i_{k}}$ changes one entry, which is the largest entry in $\bv$.
Hence, setting $\|\bv\|_{\infty}$ to be the max-norm, and for $T$ large, we can rewrite $\cN_{\sG}(T)$ in \eqref{eq:NGT} 
%for large $T$ 
as 
\be\label{eq:NGTorb}
\cN_{\sG}(T)
=
4
+
\#
\left\{
\bv\in
\cO %\G\cdot\bv_{0}
:
\bv\neq\bv_{0},\
\|\bv\|_{\infty}<T
\right\}
.
\ee
Here the first ``4'' accounts for the root quadruple $\bv_{0}$.
% (and of course one should take $T>\|\bv_{0}\|_{\infty}$ in \eqref{eq:NGTorb}).

We have thus
converted the circle counting problem into something
seemingly
 more tractable: the counting problem for a $\G$-orbit.
That said, we clearly 
 need % to have 
 a better understanding of the group $\G$.
%
%This group is 
%clearly 
%isomorphic to $\cA$, and
%
%
%By
Returning to the
Descartes
 form $Q$ in \eqref{eq:QDes}, we 
 %see 
 have
 by
 construction (and one can check directly) that for each $j=1,\dots,4$,
$$
Q(S_{j}\cdot\bv)=Q(\bv),
$$
for any $\bv\in\R^{4}$. That is, each generator $S_{j}$  lies in the so-called {\it orthogonal group} preserving the quadratic form $Q$,
$$
O_{Q}%(\R)
:=
\left\{
g\in\GL_{4}%(4,\R)
:
Q(g\cdot \bv)
=
Q(\bv),\
\forall\bv%\in\R^{4}
\right\}
.
$$
Hence $\G$ also sits inside $O_{Q}$, and moreover inside $O_{Q}(\Z)$, the group of matrices in $O_{Q}$ with integer entries. 
The latter is a %very 
well understood {\it algebraic} group, % ************ 
again meaning that any solution to a certain set of polynomial equations gives an element in $O_{Q}%(\Z)
$, and vice-versa. But $\G$ is quite a mysterious group, in particular having infinite index in $O_{Q}(\Z)$ (this fact is equivalent to $\cA$ %acting on $\bH^{3}$ 
having infinite co-volume). 
%Also %recall 
It is also worth noting
here that the general membership problem in a group is known to be undecidable \cite{Novikov1955}, so 
presenting a matrix group via its
%just giving 
generators leaves much to be desired.\footnote{%
That said, for our particular group $\G$, one can use a reduction algorithm to root quadruples to determine membership.}
%so the presentation of $\G$ in terms of its generators does not hel

Just as in Zaremba's problem, we can now again call this orbit $\cO$ {\it thin}; indeed,
for
 the counting problem with $\G$ replaced by the full $O_{Q}(\Z)$, % is 
%now classical
%well-known, 
%and 
standard arguments 
in automorphic forms 
or ergodic theory 
\cite{DukeRudnickSarnak1993, %} 
%\cite{
EskinMcMullen1993}
show 
that
\be\label{eq:FullOrb}
\#\{\bv\in O_{Q}(\Z)\cdot\bv_{0}:\|\bv\|_{\infty}<T\}\sim c\ T^{2},\qquad \text{ as $T\to\infty$},
\ee
for some $c>0$. So comparing \eqref{eq:FullOrb} to \eqref{eq:NGTorb}, \eqref{eq:KO}
and \eqref{eq:gdIs}, where the power drops from $T^{2}$ to $T^{\gd}$ with $\gd<2$, we see that the $\G$ orbit is quite degenerate, having many fewer points.
\\
%\newpage

%\newpage

\subsubsection{Sketch  of the counting statement}\label{sec:ApCount}\

Finally, we explain the aforementioned spectral interpretation, by first giving an analogous
elementary example  of a counting statement in another discrete group: the integers. Let us spectrally count the number of integers of size at most $T$:
$$
\cN_{\Z}(T):=\#\{n\in\Z:|n|<T\}.
$$
Of course this is a trivial problem, 
\be\label{eq:cNZ}
\cN_{\Z}(T)=2T\pm1,
\ee 
but it will be instructive to analyze it by harmonic analysis. To this end, let
$$
f(x):=\bo_{\{|x|<1\}},
$$
where $\bo_{\{\cdot\}}$ is the indicator function. Scale $f$ to
$$
f_{T}(x):=f(x/T)=\bo_{\{|x|<T\}},
$$
and periodize it with respect to the discrete group $\Z$:
\be\label{eq:FTis}
F_{T}(x):=\sum_{n\in\Z}f_{T}(n+x)
.
\ee
Then we have %clearly
\be\label{eq:FTN0A}
F_{T}(0)=
\sum_{n\in\Z}
\bo_{\{|n|<T\}}
=\cN_{\Z}(T).
\ee
By construction, $F_{T}(x)=F_{T}(x+1),$ that is, it takes values on the circle 
$
X:=\Z\bk\R,
$ 
and is square-integrable, $F_{T}\in L^{2}(X)$. 
The Laplace operator 
$$
\gD:=-\Div\circ\grad=-{\dd^{2}\over \dd x^{2}}
$$ 
on smooth functions 
%extends acts on 
can be extended to act on 
the whole Hilbert space
$L^{2}(X)$ and is self-adjoint and positive definite (by our choice of sign)
%
%The latter is a Hilbert space 
with respect to the standard inner product 
$$
\<F,G\>=\int_{X}F(x)\bar G(x)dx
.
$$
(Proof: partial integration.) 
Its {\it spectrum} $\Spec(\gD)$ is just the set of its eigenvalues, with multiplicity.
%Undergraduate 
Elementary
Fourier analysis shows that eigenfunctions of $\gD$ invariant under $\Z$-translations
are scalar multiples of 
$$
\vf_{m}:x\mapsto e^{2\pi i m x}
$$ 
for $m\in\Z$. This function has Laplace eigenvalue 
$$
\gl_{m}=+4\pi^{2}m^{2},
$$ %which 
and these
fully span the spectrum (they have multiplicity two, except when $m=0$).
Expanding spectrally gives
\be\label{eq:FTspec}
F_{T}(x)
=
\sum_{\gl_{m}\in\Spec(\gD)}
\<F_{T},\vf_{m}\>\vf_{m}(x)
,
\ee
where equality is in the $L^{2}$-sense.
(Note that  the $\vf_{m}$ are already scaled to have unit $L^{2}$-norm.) The bottom of the spectrum $\gl_{0}=0$ corresponds to the constant function $\vf_{0}(x)=1$, and 
 contributes
 the entire ``main term'' 
 in
\eqref{eq:cNZ}
 to \eqref{eq:FTspec}:
% evaluated at $x=0$:
$$
\<F_{T},\vf_{0}\>\cdot\vf_{0}%(0)
=
\left(
\int_{\Z\bk\R}
\sum_{n\in\Z}
f_{T}(n+x)
\cdot 1\,dx
\right)
\cdot 1
=
T
\int_{\R}
f(x)
dx
=2T
,
$$
after
inserting \eqref{eq:FTis}, 
 a change of variables, and ``unfolding'' $\int_{\Z\bk\R}\sum_{\Z}$ to  just $\int_{\R}$.
%Of course,
That said, %$L^{2}$-functions do not have values 
the equality \eqref{eq:FTspec} is in the $L^{2}$ sense, not pointwise
(we cannot evaluate \eqref{eq:FTspec} at the point $x=0$, as needed in \eqref{eq:FTN0A}), and moreover the rest of the spectrum in \eqref{eq:FTspec}, if bounded in absolute value, does not converge. But there are standard methods (smoothing and later unsmoothing) which overcome these technical irritants. 

A version of the above works with the Apollonian group $\G$ in place of $\Z$, once one overcomes a number of further technical obstructions.
The reader may wish to omit the following paragraph on the first pass; it is not essential to the sequel.

We now need non-abelian harmonic analysis on the space $L^{2}(X)$ with 
$$
X:=\cA\bk\bH^{3},
$$ 
the hyperbolic 3-fold in Figure \ref{fig:ApOrb}. The (positive definite) hyperbolic Laplacian is
$$
\gD=-y^{2}\left(
{\dd^{2}\over \dd x_{1}^{2}}+
{\dd^{2}\over \dd x_{2}^{2}}+
{\dd^{2}\over \dd y^{2}}
\right)
+
y{\dd\over \dd y}
$$
in the coordinates \eqref{eq:bH3}. 
The spectrum 
in this setting, %is governed by 
as studied by Lax-Phillips \cite{LaxPhillips1982},
has both continuous and discrete components 
(though only a finite number of the latter%
%, consistent with the Phillips-Sarnak philosophy \cite{PhillipsSarnak1985} that non-arithmetic groups have a few Maass forms, if any
).
As $X$ has {\it infinite} volume, the constant function is no longer square-integrable, 
and the bottom eigenvalue $\gl_{0}$ is strictly positive.
A beautiful result in Patterson-Sullivan  theory \cite{Patterson1976, Sullivan1984} relates this eigenvalue to the 
Hausdorff dimension of the limiting gasket $\sG$, namely
% via 
$$
\gl_{0}=\gd(2-\gd)
.
$$ %,
%and t
The corresponding base
eigenfunction $\vf_{0}$ 
%is the so-called Patterson-Sullivan eigenfunction, following the 
replaces the role of the constant function. Here we have used crucially that $\cA$ is geometrically finite, and that $\gd>1$, see \eqref{eq:gdIs}. Even this is %not enough
insufficient: because of the non-Euclidean norm $\|\cdot\|_{\infty}$ in \eqref{eq:NGTorb}, one must work not on $X$ but its unit tangent bundle $Y:=T^{1}(X)$. And moreover we do not know how to handle the continuous spectrum directly, applying instead general results in
%the theory of unitary representations of semisimple groups
the
 representation theory
 of semisimple groups
  about ergodic properties of flows %in 
  on
  $Y$. At this point, we will not say more  about the proof, 
 %letting 
 inviting
 the 
 interested
 reader 
 to
 consult the original references \cite{KontorovichOh2011, Vinogradov2012, LeeOh2012}.
\\

%\newpage

\subsection{The Local-%to-
Global Problem}\

%We turn our focus now to a problem in some sense orthogonal to the one studied in the previous subsection. 
Assume now that $\sG$ is
not only
 bounded
  %and 
 but also integral (recall that
 this means it has
 %again, having 
 only integers for curvatures).
Shrinking the gasket by a factor of two will double all of %the 
its
curvatures, making them all even.
So we 
should
%may rule this out,  %after rescaling that this is not the case
rescale %ing 
an integral gasket
 %to be 
 to make it
 {\it primitive}, 
 %that is, 
 %we assume 
meaning
 there is no number other than $\pm1$ dividing all of the curvatures.
%Before we were counting circles 
%Rather than count circles,
%let us now count {\it integers} in $\sG$, that is, count curvatures without multiplicity.
%That is, l
%
 In fact, all of the salient features of the problem persist if we fix $\sG$ to be
  %such as that 
the packing shown  in Figure \ref{fig:Apoll}, and we do so
% from now on. 
henceforth.
Recall that the problem we wish to now address is: How many curvatures are there up to some parameter $N$, counting without multiplicity, that is, counting only distinct curvatures?
% (The parameter $T$ counted circles with multiplicity)?

First some more notation: let $\sB=\sB(\sG)$ be the set of all curvatures of circles $C$ in the gasket $\sG$,
$$
\sB:=\{n\in\Z: \exists C\in\sG\text{ with }b(C)=n \}
,
$$
and call $n$ {\it represented} if $n\in\sB$.
%
%After a moment's reflection
Staring at Figure \ref{fig:Apoll} for a moment or two, one might observe that every curvature
 %there 
in our $\sG$ is
\be\label{eq:cong}
\equiv
2, 3, 6, 11, 14, 15, 18,\text{ or }23\ (\mod 24)
.
\ee
These are the local obstructions for $\sG$, and we call $\sA=\sA(\sG)$ the set of all {\it admissible} numbers $n$ satisfying \eqref{eq:cong}. 
In general, one calls $n$ admissible if, as before, it is everywhere locally represented,
\be\label{eq:adminA}
n\in\sB(\mod q), \ \forall q\ge1.
\ee
It cannot be the case that $\sA=\sB$, since, for example, $n=15$ is admissible, but a circle of radius $1/15$ does not appear in %the 
our
gasket. Nevertheless, as in Zaremba's problem, we have the following

\phantomsection
\begin{conja}
%[{Graham {\it et al} 2003 \cite[p. 37]{GrahamLagarias2003}}]
\label{conj:A}
Every sufficiently large admissible number
% appears in Figure \ref{fig:Apoll}.\\
is the curvature of some circle in $\sG$.
\end{conja}

This conjecture is stated by Graham-Lagarias-Mallows-Wilks-Yan \cite[p. 37]{GrahamLagarias2003}, 
%who call it the ``Strong Density Conjecture''
in the first of a lovely series of papers % on this topics.
Apollonian gaskets and generalizations.
They observe empirically that congruence obstructions for any integral gasket seem to be to the modulus $24$, 
and this %wa
is completely clarified (as we explain below) by Fuchs \cite{Fuchs2011}
in her thesis. 
Further convincing numerical evidence towards
 %it 
 the conjecture
 is %presented
given in Fuchs-Sanden \cite{FuchsSanden2011}.
Here is some recent progress.

\phantomsection
\begin{thma}[Bourgain-K. 2012 
\cite{BourgainKontorovich2012}]\label{thm:A}
Almost every %sufficiently large 
admissible number 
%appears in Figure \ref{fig:Apoll}. 
is the curvature of some circle in $\sG$. %\\
\end{thma}

Again, ``almost every'' is in the sense of density, that
\be\label{eq:BtoA}
{
\#(\sB\cap[1,N])\over
\#(\sA\cap[1,N])
}
\to
1,
\ee
as $N\to\infty$. %Of course
It follows 
 from 
 \eqref{eq:cong}
%implies 
that for $N$ large,  $\#(\sA\cap[1,N])$ is %the integer part of 
about
$N/3$ (there are $8$ admissible residue classes mod $24$), so \eqref{eq:BtoA} is equivalent to 
$$
\#(\sB\cap[1,N])\sim \frac N3.
$$
%The

%Here is s
Some  history 
on %of
%previous progress on
 this problem:
%is as follows. 
Graham {\it et al} \cite{GrahamLagarias2003} already made the first 
progress,
%dent
 proving that
\be\label{eq:GLbnd}
\#(\sB\cap[1,N])
\gg
N^{1/2}
.
\ee
Then Sarnak \cite{SarnakToLagarias} showed
\be\label{eq:SarBnd}
\#(\sB\cap[1,N])
\gg
{N\over \sqrt{\log N}}
,
\ee
before Bourgain-Fuchs \cite{BourgainFuchs2011} settled the so-called Positive Density Conjecture, that
\be\label{eq:BFbnd}
\#(\sB\cap[1,N])
\gg
N
.
\ee

A key observation in the proof of Theorem  \hyperref[thm:A]{A} is that the problem is nearly identical to Zaremba's, in the following sense. Recall from \eqref{eq:cOAp} that the orbit $\cO=\G\cdot\bv_{0}$ of the root quadruple $\bv_{0}$ under the Apollonian group $\G$ contains all quadruples of curvatures, and in particular its entries consist of all curvatures in $\sG$. Hence the set  $\sB$ of all curvatures is simply the finite union of sets of the form
\be\label{eq:BvGv}
\<\bw_{0}, \cO\>
=
\<\bw_{0}, \G\cdot\bv_{0}\>
,
\ee
%where 
as
$\bw_{0}$ 
%is one of
ranges through
 the standard basis vectors $\bbe_{1}=(1,0,0,0)^{t},\dots$, $\bbe_{4}=(0,0,0,1)^{t}$, 
 each
 picking off
  %the 
  one
  entry %ies 
  of $\cO$.
A heuristic analogy between Zaremba and the Apollonian problem is actually already given in 
\cite[p. 37]{GrahamLagarias2003}, but it is crucial for us that both problems are exactly of the form \eqref{eq:BvGv}; compare to \eqref{eq:DAvGv}.
That is, $n$ is represented
 if 
 and only if
 there is a $\g$ in the Apollonian group $\G$ and 
 some
 $\bw_{0}\in\{\bbe_{1},\dots,\bbe_{4}\}$ so that 
\be\label{eq:nRepA}
n=\<\bw_{0},\g\cdot\bv_{0}\>.
\ee
%
%\newpage
%
Before saying more about the proof of Theorem \hyperref[thm:A]{A}, we first discuss admissibility in greater detail.
%\\

\subsubsection{Local obstructions}\

Through \eqref{eq:BvGv}, the admissibility condition \eqref{eq:adminA} is again reduced 
%just 
to the study of the reduction of  $\G$ modulo  $q$.
%
%
%\newpage
%
An
%other 
important feature 
here
is that,
like in the Zaremba case, 
the group
%of 
$\G$ 
%it 
is  %``algebraically large'', being 
{Zariski dense} in $O_{Q}$. 
Recall % again
 that this  means:
% that 
if $P(\g)$ is a polynomial in the entries of a $4\times4$ matrix $\g$ which vanishes for every $
\g%\{g_{i,j}\}_{1\le i,j\le 4}
\in\G$, then $P$ also vanishes on all complex matrices in $O_{Q}$. % with entries in $\C$.

%\subsubsection{Spin representation from $\cA$ to $\G$}\
We would like again to exploit strong approximation, but neither $O_{Q}$ nor its  the orientation preserving subgroup $\SO_{Q}:=O_{Q}\cap\SL_{4}$
%do not 
have this property (%the latter is
being
 not even %simply 
 connected).
But there is a %well-known 
standard
method of applying strong approximation anyway, by first passing to
% the so-called {\it 
%spin group}, 
a certain cover,
as we now describe.

%\newpage

%As
%It is well-known in
From the theory of rational quadratic forms \cite{Cassels1978}, % that 
special orthogonal groups 
are
%double
 covered by so-called 
%have 
{\it spin 
%double covers
groups}, and it is a pleasant accident
% isomorphism that $\SO(3,1)$ is double-covered by $\SL_2(\C)$.
that,
since $Q$ has signature $(3,1)$,
 the spin group of $\SO_{Q}(\R)$ is isomorphic to $\SL_{2}(\C)$; let us %recall 
explain
this covering map.
%For simplicity, 
The formulae are nicer if
we first change variables
 (over $\Q$) from the quadratic form $Q$ to the equivalent form 
$$
\tilde Q(x,y,z,w):=
xw+y^{2}+z^{2}
.
$$
Observe that the matrix
$$
M%_{\tilde Q}
:=
\mattwo{-x}{y+iz}{y-iz}{w}
$$
has determinant equal to $-\tilde Q$ and is Hermitian, that is, fixed under transpose-conjugation. The group $\SL_2(\C)$,
consisting of $2\times2$ complex matrices of determinant one, acts on $M%_{\tilde Q}
$ by
$$
\SL_{2}(\C)\ni g:M
%_{\tilde Q}
\mapsto g\cdot M%_{\tilde Q}
\cdot \bar g^{t}=:M%_{\tilde Q}
'
=
\mattwo{-x'}{y'+iz'}{y'-iz'}{w'}
,
$$
with 
$
%g(M_{\tilde Q})
M%_{\tilde Q}
'
$ 
also Hermitian and of determinant $-\tilde Q$. Then it is easy to see that $(x',y',z',w')^{t}$ is a linear change of variables from $(x,y,z,w)^{t}$, via left multiplication by a matrix whose entries are quadratic in the entries of $g$. Explicitly, if
\be\label{eq:g}
g=
\mattwo{a+\ga i}{b+\gb i}{c+\g i}{d+\gd i} %\in\SL_2(\C)
,
\ee
then
the change of variables matrix is
\be\label{eq:rhoG}%\be%a
%&&
%\hskip-.5in
%\tilde
%\iota(g)
%\mapsto
%\\
%\nonumber
%&&
%:=
{1\over |\det(g)|^{2}}
\left(
\begin{array}{cccc}
 a^2+\alpha ^2 & 2 (a c+\alpha  \gamma ) & 2 (c \alpha - a \gamma)  & -c^2-\gamma ^2 \\
 a b+\alpha  \beta  & b c+a d+\beta  \gamma +\alpha  \delta  & d \alpha +c \beta -b \gamma -a \delta  & -c d-\gamma  \delta 
   \\
 a \beta -b \alpha  & -d \alpha +c \beta -b \gamma +a \delta  & -b c+a d-\beta  \gamma +\alpha  \delta  & d \gamma -c \delta
    \\
 -b^2-\beta ^2 & -2 (b d+\beta  \delta ) & 2( b \delta -d \beta ) & d^2+\delta ^2
\end{array}
\right)
.
\ee%\ee%a
Let $\tilde\rho$ be the (rational) map from $\SL_{2}(\C)$ to $\GL_{4}(\R)$, sending \eqref{eq:g} to \eqref{eq:rhoG}; then by construction (again one can  verify directly) the image is in $\SO_{\tilde Q}(\R)$. 
(Some minor technical points: Being quadratic in the entries, $\tilde\rho$ is a double cover, with $\pm I$ having the same image. Moreover, $\SL_2(\C)$ is connected while $\SO_{\tilde Q}(\R)$ has two connected components, so $\tilde\rho$ only maps onto
% $\SO_{\tilde Q}^{\circ}(\R)$,
 the identity component.) %; note in particular that the top left entry of \eqref{eq:rhoG} is non-negative.)
Then changing variables from $\tilde Q$ back to the Descartes form $Q$ by a conjugation, one gets the desired map 
$$
\rho:\SL_{2}(\C)\to\SO_{Q}(\R).
$$ %, as desired.

It is straightforward then to compute the pullback of $\G\cap\SO_{Q}$ under $\rho$ (see \cite{GrahamLagariasMallowsWilksYanI, Fuchs2011}), the answer being the following
\begin{lem}\label{lem:genSL2C}
There is\footnote{%
And
% we could 
one can easily
write it down explicitly: it is a conjugate of \eqref{eq:rhoG}, but much messier and 
%completely un-
not particularly
enlightening.  We spare the reader.%
} a homomorphism $\rho:\SL_{2}(\C)\to\SO_{Q}(\R)$ so that 
the group $\tilde\G:=\rho^{-1}(\G\cap\SO_{Q})$ sits in $\SL_{2}(\Z[i])$ and is generated by 
\be\label{eq:genSL2C}
\tilde\G=\<
\pm
\mattwo1{2}{0}1,\
\pm
\mattwo{1}{0}{2}{1},\ %\text{ and }\quad
\pm
\mattwo{1+2i}{-2}{-2}{1-2i}
\>.
\ee
Moreover, recalling the generators $S_{j}$ for $\G$ in \eqref{eq:SjAp}, one can arrange $\rho$ so that 
$
\rho:
\bigl( \begin{smallmatrix}
1&2\\ 0&1
\end{smallmatrix} \bigr)
\mapsto S_{2}S_{3}
,$
and
$
\rho:
\bigl( \begin{smallmatrix}
1&0\\ 2&1
\end{smallmatrix} \bigr)
\mapsto S_{4}S_{3}
.$
\end{lem}

%\newpage

In fact, we have just %explicitly 
realized a conjugate of 
the group
$\cA$ % of M\"obius transformations 
(or rather its index-two orientation preserving subgroup)  explicitly in terms of matrices in $\PSL(2,\C)\cong\Isom^{+}(\bH^{3})$.
%, as will be useful later.

From here, 
one follows the strategy outlined in \S\ref{sec:locZ}.
Using strong approximation for 
$\SL_{2}(\Z[i])$ (one considers reduction mod principal ideals $(q)$), Goursat's Lemma, some finite group theory, and 
%some 
other ingredients, Fuchs \cite{Fuchs2011} was able to determine completely the reduction of $\G$ modulo any $q$, and hence explain all local obstructions.
The answer is that all primes other than $2$ and $3$ are {\it unramified}, meaning, as in
 %(recall the definition from
  \S\ref{sec:locZ},
that for $(q,6)=1$, 
$$
\G\cap\SO_{Q}\ (\mod q)\quad=\quad \SO_{Q}(\Z/q\Z).
$$ 
Recall again that the %latter 
right hand side above
is a well-understood algebraic group. 
  And moreover, 
the prime $2$ stabilizes (with the same meaning as \S\ref{sec:locZ}) at the power $e_{0}(2)=3$, that is at  $8$, and the prime $3$ stabilizes immediately at $e_{0}(3)=1$. Then $\G(\mod 24)$ is some  explicit 
finite group, 
 %which one can write down explicitly, 
 and looking at all the values of \eqref{eq:BvGv} for the given root quadruple $\bv_{0}(\sG)$, one 
  immediately
  sees all admissible residue classes. 
\\

\subsubsection{Partial Progress}\label{sec:progA}\

%Before leaving our discussion of the Apollonian problem, we illustrate how 
Lemma \ref{lem:genSL2C} %gives even more information:
can already be quite useful, in particular, it easily implies 
\eqref{eq:GLbnd}
and 
\eqref{eq:SarBnd}, as follows.

The Apollonian group $\G$
%of
%\eqref{eq:GamAp}
 contains the matrix $S_{4}S_{3}$, which
 by Lemma \ref{lem:genSL2C}
  is the image under $\rho$ of $
\bigl( \begin{smallmatrix}
1&0\\ 2&1
\end{smallmatrix} \bigr)
$. 
The latter (and hence the former) is a {\it unipotent} matrix, meaning that all its eigenvalues are equal to $1$. These
have the important property that they
 grow only polynomially under exponentiation; in particular,
 $
 \bigl( \begin{smallmatrix}
1&0\\ 2&1
\end{smallmatrix} \bigr)^{k}
=
\bigl( \begin{smallmatrix}
1&0\\ 2k&1
\end{smallmatrix} \bigr)
,
$
and
one can check directly from \eqref{eq:SjAp} that
$$
(S_{4}S_{3})^{k}
=
\left(
\begin{array}{cccc}
 1 & 0 & 0 & 0 \\
 0 & 1 & 0 & 0 \\
 4 k^2-2 k & 4 k^2-2 k & 1-2 k & 2 k \\
 4 k^2+2 k & 4 k^2+2 k & -2 k & 2 k+1
\end{array}
\right)
.
$$
Put
%ting 
the above matrix into  \eqref{eq:nRepA} with the root quadruple $\bv_{0}$ for our fixed gasket 
%=(-10,18,23,27)^{t}$ 
from \eqref{eq:bv0G},
and take $\bw_{0}=\bbe_{4}$, say.
%this means that 
%gives that, 
Then
for any $k\in\Z$, the number
\be\label{eq:quadK}
\<\bbe_{4}\ , \ (S_{4}S_{3})^{k}\cdot\bv_{0}\>=
32 k^2
+ 24 k 
+ 
27 
\ee
is represented. 
That is, the set of represented numbers contains the values of this quadratic polynomial.
From this observation, made in \cite{GrahamLagarias2003}, %so 
it is immediate that
\eqref{eq:GLbnd} %is immediate.
holds.
Geometrically, these curvatures correspond to circles in the packing tangent to
 $C_{1}$ and $C_{2}$, since these are  untouched by the corresponding reflections through $\tilde C_{4}$ and $\tilde C_{3}$.
For example, the values $k=-2,-1,0,1,2$
in \eqref{eq:quadK}
give curvatures $107, 35, 27, 83, 203$, respectively.
These are  visible in Figure \ref{fig:Apoll};
%between curvatures $18$, $23$ and $27$.)
%near the middle, and 
they are all
tangent to the circles of curvature $-10$ (the bounding circle) and $18$, skipping every other such circle.
Using $\bw_{0}=\bbe_{3}$ instead of $\bbe_{4}$ in \eqref{eq:quadK} gives the polynomial 
$
32 k^2
- 8 k + 
23 
,
$
the values of which
correspond to the skipped circles.
\\

To prove \eqref{eq:SarBnd}, we make the following observation, due to Sarnak \cite{SarnakToLagarias}. 
It is well known that the matrices
$
\pm
 \bigl( \begin{smallmatrix}
1&0\\ 2&1
\end{smallmatrix} \bigr)
$
and
$
\pm
 \bigl( \begin{smallmatrix}
1&2\\ 0&1
\end{smallmatrix} \bigr)
$
(which map under $\rho$ to $S_{4}S_{3}$ and $S_{2}S_{3}$, respectively)
 generate the
 %following 
 group
$$
\gL(2):=
%\<
%\pm
%\bp
%1&2\\
%0&1
%\ep
%,
%\pm
%\bp
%1&1\\
%2&1
%\ep
%\>
%=
\left\{
\mattwo abcd\in\SL_{2}(\Z)
:
\begin{array}{c}
a\equiv d\equiv 1(\mod 2)
\\
 b\equiv c\equiv 0(\mod 2)
 \end{array}
\right\}
.
$$ 
This is the
so-called level-2 {\it principal  congruence} subgroup of $\SL_{2}(\Z)$.
%,
%, and we will denote it by $\gL(2)$.
Hence by Lemma \ref{lem:genSL2C},
the group $\G$ contains
\be\label{eq:XiDef}
\Xi:=\<\ S_{2}S_{3}\ , \ S_{4}S_{3}\ \>%<\G
=
\rho\left(\gL(2)\right)
.
\ee
The point is that $\gL(2)$
 is arithmetic: %again algebraic;
 %, being defined by polynomial equations (namely $ad-bc=1$ and the congruences, $a=2x+1$, etc.).
%
%Hence
%in particular,
 for
 any
 % even $k$ and odd $\ell$ with
  %$(2k,\ell)=1$
integer
  $\ell$ coprime to $2k$, there is a matrix 
$
 \bigl( \begin{smallmatrix}
*&*\\ 2k&\ell
\end{smallmatrix} \bigr)
$
in $\gL(2)$.
One can work out, with 
the same
$\bv_{0}$ and $\bw_{0}$ as above, that
\be\label{eq:quadKL}
\<\bbe_{4}\ ,\
\rho
%\left(
\bp
*&*\\
2k&\ell
\ep
%\right)
\cdot
\bv_{0}
\>
=
32 k^2 + 24 k \ell + 17 \ell^2
+10 
.
\ee
For example, the choices $(2k,\ell)=(4,-3),(2,-1),(4,-1),$ and $(6,-1)$ give curvatures $147, 35, 107,$ and  $243$, respectively, visible %down 
up
the %right 
left
side of Figure \ref{fig:Apoll}, all tangent to the bounding circle
(since $\Xi$ in \eqref{eq:XiDef} fixes $C_{1}$).
Observe also that setting $\ell=1$ in \eqref{eq:quadKL} recovers \eqref{eq:quadK}.
So, as observed by Sarnak \cite{SarnakToLagarias}, the set $\sB$ of represented numbers contains all {\it primitive}  (meaning with $2k$ and $\ell$ coprime) 
values
of the  shifted binary quadratic form in \eqref{eq:quadKL}. 
Note that the quadratic form has discriminant 
$24^2 - 4\cdot 32\cdot 17=-1600$, and so \eqref{eq:quadKL} is definite, taking only positive values.
The number of distinct primitive values of \eqref{eq:quadKL} up to $N$ was determined by  Landau \cite{Landau1908}:
it is asymptotic to a constant times $N/\sqrt{\log N}$, thereby proving \eqref{eq:SarBnd}.
A much more
delicate and
 clever but still ``elementary'' (no automorphic forms are harmed) argument goes into the proof of
\eqref{eq:BFbnd}, using an ensemble of such shifted binary quadratic forms.
For Theorem  \hyperref[thm:A]{A}, one %invokes
needs
 the theory of automorphic representations for the full Apollonian group, as 
hinted to at the end of \S\ref{sec:ApCount}.
\\

We now leave the discussion of the Apollonian problem, returning to it again in \S\ref{sec:C}.

 \newpage

\section{The Thin Pythagorean Problem}\label{sec:P}

 \begin{figure}
        \begin{subfigure}[t]{3in}
                \centering
                %\reflectbox
                {\includegraphics[height=1.5in]{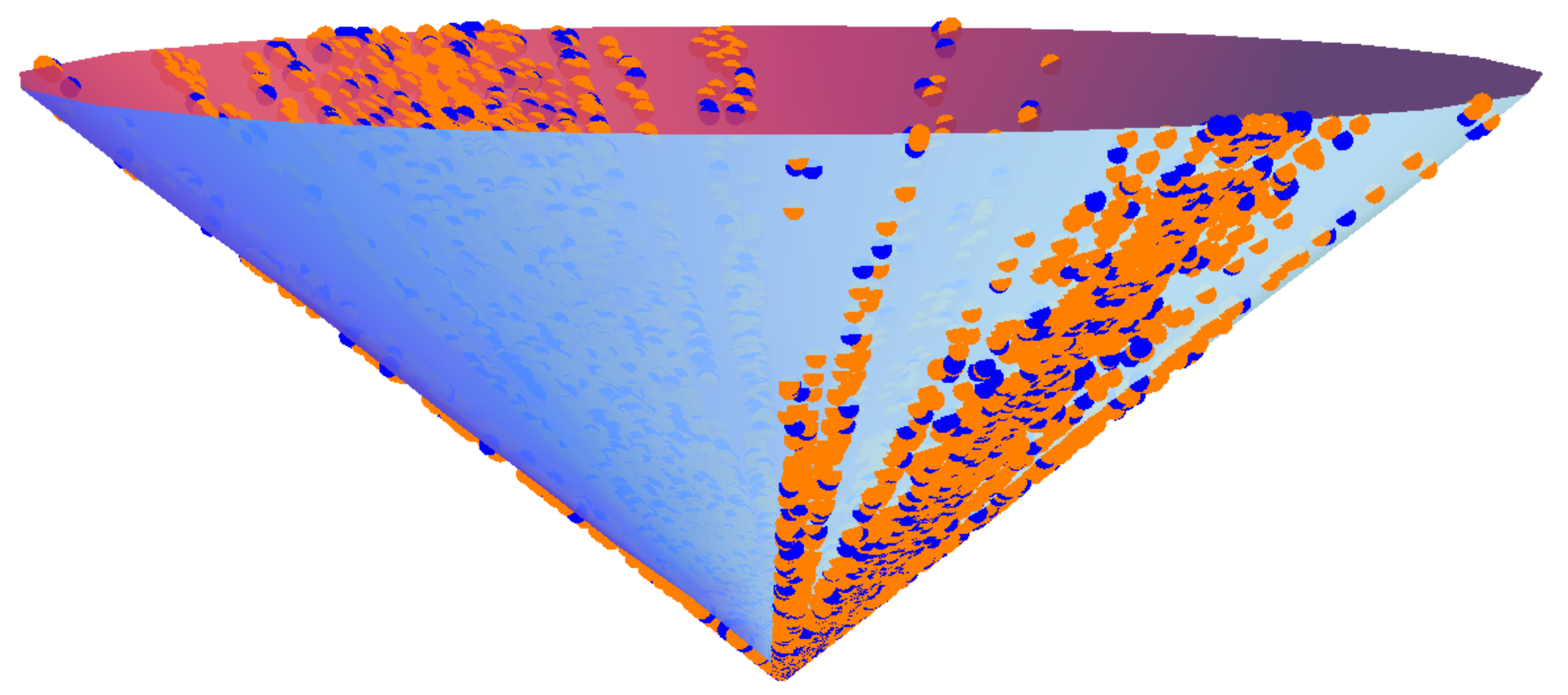}}
                \caption{View from the side}
                \label{fig:PythA}
        \end{subfigure}%
\qquad
        \begin{subfigure}[t]{1.5in}
                \centering
                \includegraphics[height=1.5in]{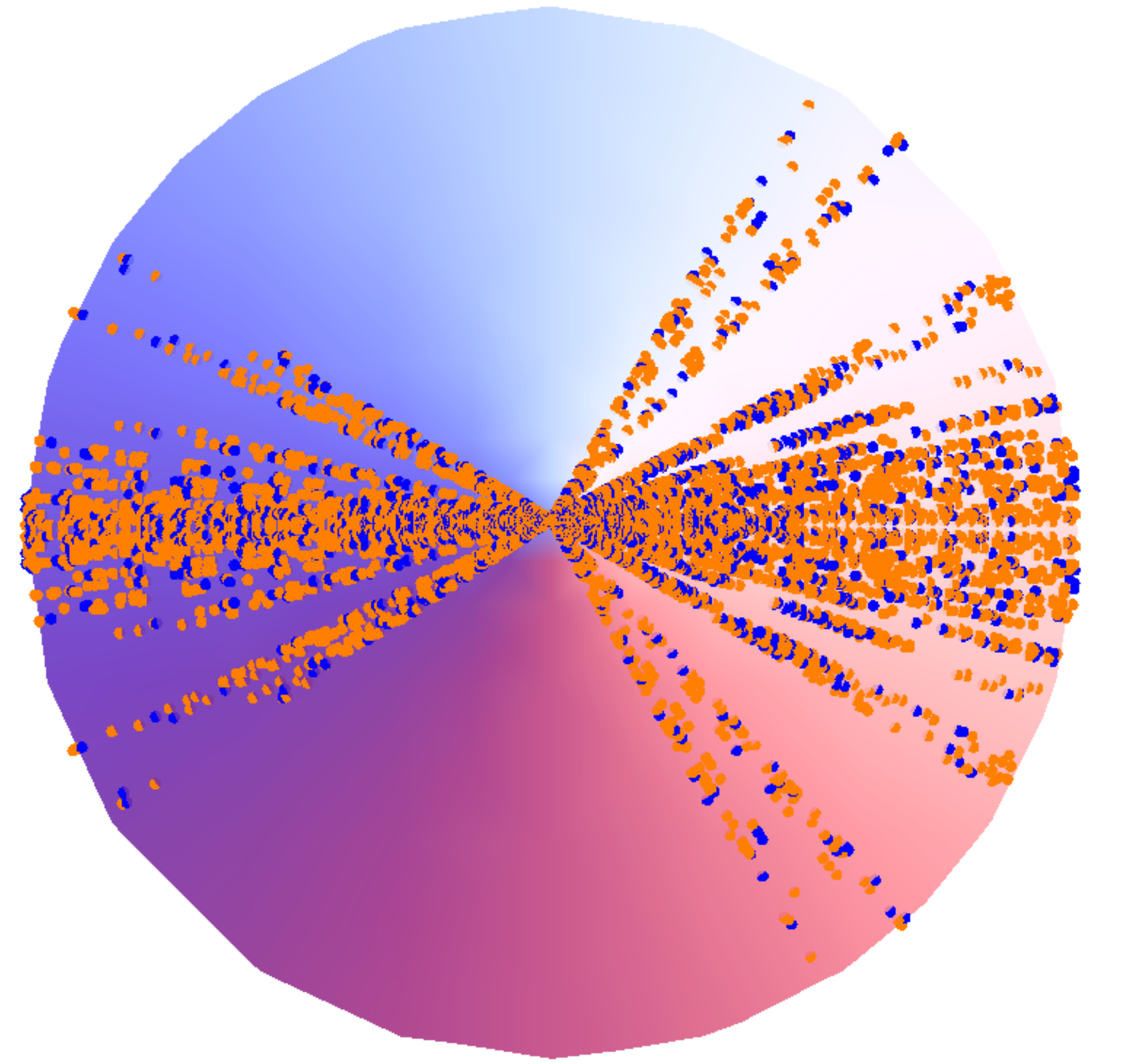}
                \caption{View from below}
                \label{fig:PythB}
        \end{subfigure}

\caption{The thin
Pythagorean 
orbit $\cO$ %of Pythagorean triples 
in \eqref{eq:cOIsP}. Points are marked according to whether the
 %sum of the 
 hypotenuse % and the even side 
 is % the square of 
 %a  
 prime (\protect\includegraphics[width=.1in]{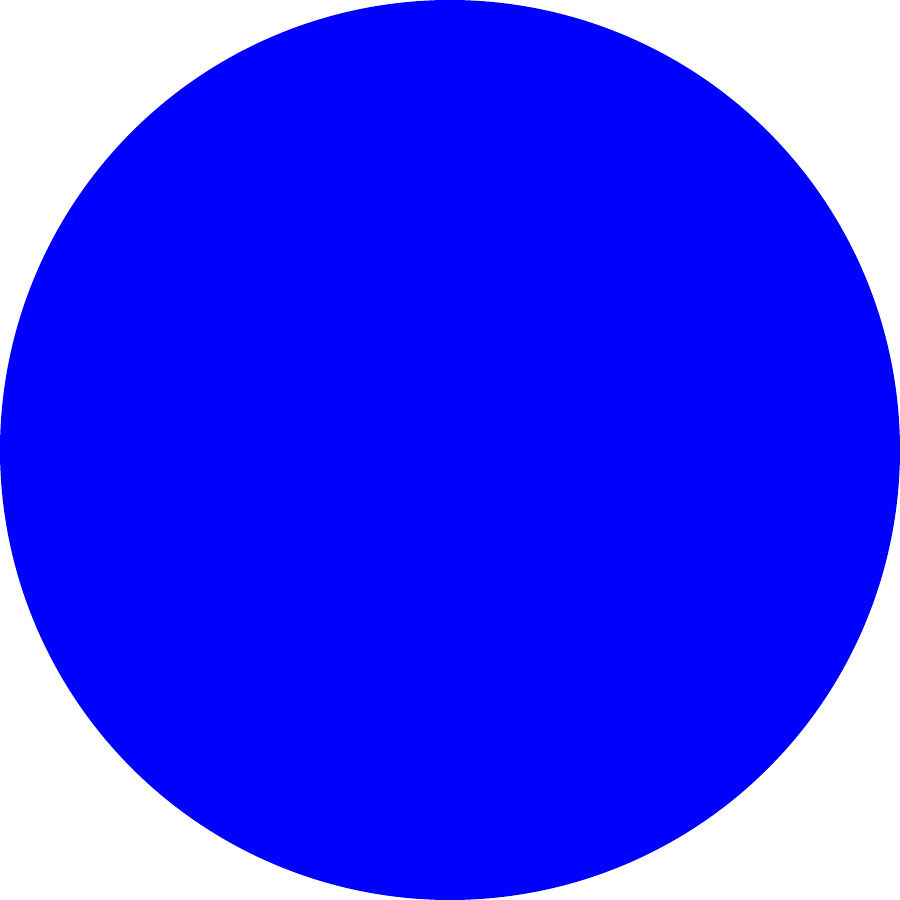}) or %a 
  composite (\protect\includegraphics[width=.1in]{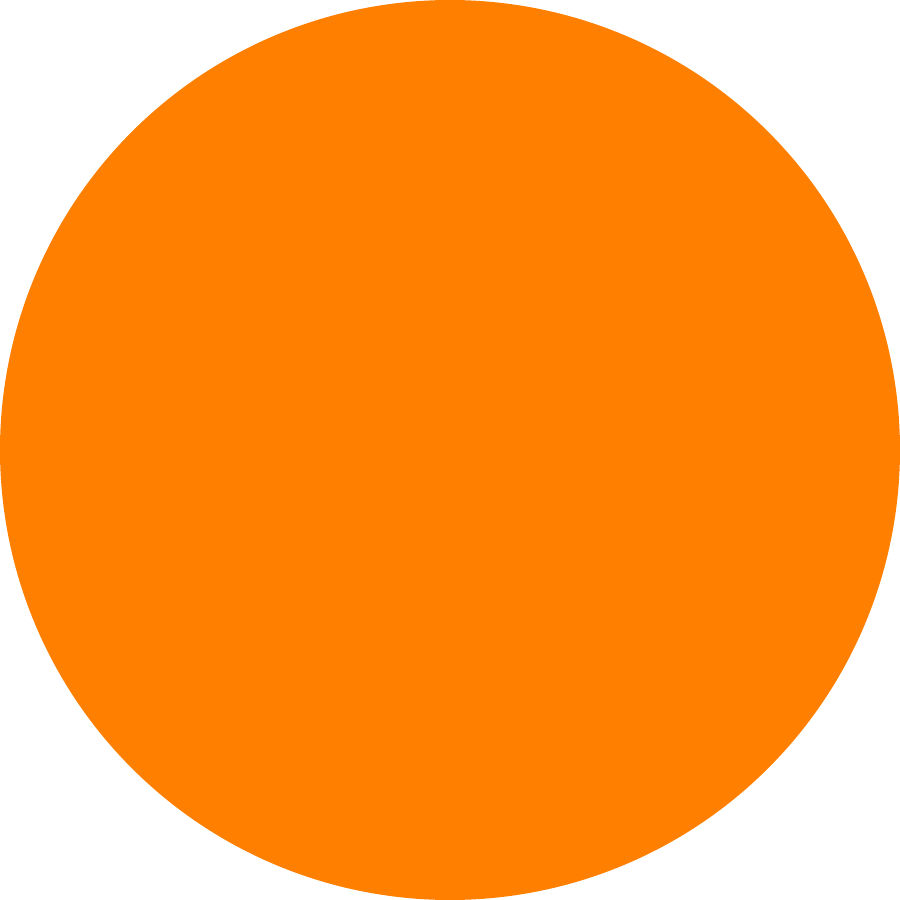})% number
  .}
\label{fig:Pyth}
\end{figure}

A {\it Pythagorean triple} $\bx=(x,y,z)^{t}$ is a %n integral 
point on the cone 
\be\label{eq:conePyth}
Q(\bx)=0,
\ee 
where $Q$ is the ``Pythagorean quadratic form''
$$
Q(\bx):=x^{2}+y^{2}-z^{2}.
$$
Throughout we consider only {\it integral} triples, $\bx\in\Z^{3}$, and assume that $x,y,$ and $z$ are coprime;
%, that is, the 
such a 
triple is called {\it primitive}. 
%From e
Elementary considerations %, this 
then
force the hypotenuse $z$ to be odd, and $x$ and $y$ to be of opposite parity; we assume %throughout 
henceforth
that $x$ is odd and $y$ is even. The cone has a singularity at the origin, so we only consider its top half, %
%\footnote{It is easier to visualize the bottom half, so we show the  in Figure \ref{fig:Pyth}.}
 assuming subsequently that the hypotenuse is positive, $z>0$.

Diophantus (and likely the Babylonians \cite{Plimpton322}, who preceded him by about as much as he precedes us) knew how to parametrize Pythagorean 
triples: % If % $\bx$ is %is an integral, primitive Pythagorean triple with %,
%, with positive entries and
% say, 
%$x$ and $z$ odd and $y$ even
%$\bx\equiv(1,0,1)(\mod 2)$, then t
Given $\bx$, there
is a pair $\bv=(u,v)$ of
 %are 
 coprime integers %$u$, $v$ 
 of opposite parity % (one even, one odd), 
so that
\be\label{eq:PythParam}
\left\{
\begin{array}{ccc}
x&=&u^{2}-v^{2}\\
y&=&2uv\\
z&=&u^{2}+v^{2}.
\end{array}
\right.
\ee
That the converse is true is elementary algebra: any such pair $\bv$ %=(u,v)$ %plugged 
inserted
into \eqref{eq:PythParam} gives rise to a triple $\bx$ satisfying \eqref{eq:conePyth}. For example, it is easy to see that the
 %base 
 triple 
 \be\label{eq:bx0}
 \bx_{0}=
 %(1,0,1)^{t}
(3,4,5)^{t}
 \ee 
 %arises from the 
 corresponds %ing base
 to the pair  
\be\label{eq:bv0}
\bv_{0}=
%(1,0)^{t}
(2,1)^{t}
.
\ee
%with $\bv=(u,v)$. (

\subsection{%An % Thin 
Orbits and the Spin Representation}\

%Just like the Descartes form,
As in the Apollonian case,
 the Pythagorean form $Q$ has
a {\it special} (determinant one) orthogonal group preserving it:
\be\label{eq:SOQ}
\SO_{Q}:=\{g\in\SL_{3}:Q(g\cdot \bx)=Q(\bx)\}
.
\ee
%As $Q$ is diagonal of signature $(2,1)$, this group is commonly written $\SO(2,1)$.
%and $\G$ is a subgroup of
And as before, this group is also better understood by passing to its spin cover. 
 %
 %It will turn out that,
 %s
 Since
 % that
 %for 
 the Pythagorean form 
 $Q$ has signature $(2,1)$, there is an accidental isomorphism between its spin group and $\SL_{2}(\R)$, %; 
 %we now present
%  this isomorphism is
  given explicitly as follows.
%The signature of $Q$ is clearly $(2,1)$, and there is an accidental isomorphism between the spin group of $\SO(2,1)$ and $\SL_{2}(\R)$.

Observe that $%\tilde G:=
\SL_{2}
%(\R)
$ acts on % the base 
a pair $\bv%_{0}
$
by 
% in the usual way (
left multiplication; %). 
via \eqref{eq:PythParam}, this action then extends to a linear 
%one
action 
on $\bx$. In coordinates, %one can 
it is an elementary computation
%compute %in this way 
that the action of 
$%\be\label{eq:gP}$$
%\tilde g=
%\mattwo abcd%\in\SL_{2}%(\R)
\bigl( \begin{smallmatrix}
a&b\\ c&d
\end{smallmatrix} \bigr)
$ %\ee %$$ 
on $\bv$ corresponds to left multiplication on $\bx$ by %that of
\be\label{eq:gIs}
%g:=
{1\over ad-bc}
\left(
\begin{array}{ccc}
 \frac{1}{2} \left(a^2-b^2-c^2+d^2\right) & a c-b d & \frac{1}{2}
   \left(a^2-b^2+c^2-d^2\right) \\
 a b-c d & b c+a d & a b+c d \\
 \frac{1}{2} \left(a^2+b^2-c^2-d^2\right) & a c+b d & \frac{1}{2}
   \left(a^2+b^2+c^2+d^2\right)
\end{array}
\right)
.
\ee
%on $\bx$. 
One can check %explicitly 
directly
from the definition \eqref{eq:SOQ} that \eqref{eq:gIs} is an element of $\SO_{Q}$, and hence we have explicitly constructed 
the spin homomorphism 
$$
\rho:\SL_{2}(\R)\to\SO_{Q}(\R):
\mattwo abcd\mapsto\eqref{eq:gIs}.
$$ 
%taking $\mattwo abcd$ to \eqref{eq:gIs}. 
%This is the analogue
%  of Lemma \ref{lem:genSL2C} for this setting. %(The %accidental isomorphism for 
%  pleasant accident 
%  the Descartes form of signature $(3,1)$ is the spin group of  $\SL_{2}(\C)$, while for the Pythagorean form of signature $(2,1)$, it is $\SL_{2}(\R)$). %) . %, we have thus realized a 

Given %an integer 
a Pythagorean
triple $\bx_{0}$,
 %satisfying \eqref{eq:conePyth}, 
 such as that in \eqref{eq:bx0}, 
the group $\G:=
\SO_{Q}(\Z)$ of
all
  {\it integer} matrices in $\SO_{Q}$ acts by left multiplication, giving the 
%a 
full
orbit $\cO=\G\cdot\bx_{0}$ of  all
Pythagorean triples (with our convention that $z>0$, $x$ is even, and $y$ is odd). 
  %Hence $\cO$ in \eqref{eq:cOIsP} is indeed an orbit of Pythagorean triples: e
%  Every $\bx$
  %\in\cO$ 
%  in % such an 
%  the
%  orbit
%  $\cO$
%  satisfies \eqref{eq:conePyth} because $\bx_{0}$ does, and the action of $\SO_{Q}$ does not change the value of $Q$.

Via \eqref{eq:PythParam} again, this $\SO_{Q}$ action on $\bx$
 %of $\SO_{Q}$ 
 is equivalent to the $\SL_{2}$ action on $\bv$.
%The above action
%on
%an 
For a %n 
 %integral 
 primitive $\bv\in\Z^{2}$, %, %To 
%one 
%we
% can
 %be made to 
% preserve
both
 the integrality and primitivity 
 are preserved
by restricting
%the domain of $\rho$
 %under 
the %above 
action 
%of $\bv$, 
%one can 
%restricting it % the action 
to just
%Restricting $\rho$ to 
the integral matrices $\SL_2(\Z)$. % preserves the integrality of $\bv$,
Moreover, one %can 
should
preserve the parity condition on $\bv$ by
%and 
restricting %the action 
further to only the 
%
% p%reserving integrality. 
%Recall the 
principal $2$-congruence subgroup 
$$
\gL(2)=
\bigg\{\g\in\SL_2(\Z):\g\equiv I(\mod 2)
\bigg\}
\
=
\
\<
\pm
%\bigl( \begin{smallmatrix}
\bp
1&2\\ 0&1
%\end{smallmatrix} \bigr)
\ep
,
%$ 
%and
\pm
\bp%\bigl( \begin{smallmatrix}
1&0\\ 2&1
\ep%\end{smallmatrix} \bigr)
\>
,
$$
which already appeared in %the  
\S\ref{sec:progA}.
%Recall it is generated by the matrices 
%$ 
One can check directly that the image \eqref{eq:gIs} of any $\g\in\gL(2)$ is an integral matrix, that is,  in $\SO_{Q}(\Z)$.
For $\bv_{0}$ corresponding to $\bx_{0}$, 
the orbit $\tilde\cO:=\tilde\G\cdot\bv_{0}$ %of 
under the full group $\tilde\G:=\gL(2)$
consists of all coprime $(u,v)$ with $u$ even and $v$ odd.
\\

Prompted by the Affine Sieve\footnote{%
We % do  not 
have insufficient room
to survey %these 
this
beautiful
 %developments
 theory, for which the reader is %invited 
directed to any number of %lovely 
excellent
surveys, e.g. \cite{Salehi2012}.}
\cite{BourgainGamburdSarnak2006, BourgainGamburdSarnak2010, SalehiSarnak2011}
one
 may 
 %ask 
% pose
% a variety of Diophantine %questions 
% problems
 %about
 wish to study
  {\it thin} orbits $\cO$ of Pythagorean triples. 
Here one replaces the full group $\SO_{Q}(\Z)$ by some finitely generated subgroup $\G$ of infinite index.
Equivalently one can consider an orbit $\tilde\cO$ of $\bv_{0}$ under an infinite index subgroup $\tilde\G$ of $\gL(2)$.
%For concreteness, we give 
We illustrate the general theory via
the following %explicit 
concrete
example.

We first give a sample $\tilde\cO$ orbit:
in  % analogy 
comparison
with the generators of $\gL(2)$, 
 let $\tilde\G$ be the group  generated % by these,
%explicitly 
by the following two 
matrices
\be\label{eq:tilGamIsP}
\tilde\G:=\<\pm\mattwo1201,\quad\pm\mattwo1041\>
%\qquad<\qquad\gL(2)
.
\ee
%It
This group clearly sits inside $\gL(2)$ but it is not immediately obvious whether it %has
is of
 finite or infinite index;
as we will show later, the index is infinite. %it is the latter. 
Taking the %same 
base pair
$\bv_{0}$
in \eqref{eq:bv0}, we form the orbit
\be\label{eq:cOtilP}
\tilde\cO:=\tilde\G\cdot\bv_{0}
.
\ee
Correspondingly, we can take the base triple $\bx_{0}$ in \eqref{eq:bx0}, and form the orbit 
\be\label{eq:cOIsP}
\cO:= \G\cdot \bx_{0}
\ee
of 
%the base triple 
$\bx_{0}$ % in \eqref{eq:bx0}, 
under the group
%$\G$ 
\be\label{eq:GammaIsP}
\G:=\<
M_{1},M_{2}\>
,
\ee
where $M_{1}$ and $M_{2}$ are the images under $\rho$ of the matrices generating $\tilde\G$;
%generated by the following two explicit $3\times3$ matrices:
one can elementarily compute from \eqref{eq:tilGamIsP} and \eqref{eq:gIs} that
%where
\be\label{eq:MsP}
M_{1}:=
\left(
\begin{array}{ccc}
 -1 & -2 & -2 \\
 2 & 1 & 2 \\
 2 & 2 & 3
\end{array}
\right)
,\qquad
M_{2}:=
\left(
\begin{array}{ccc}
 -7 & 4 & 8 \\
 -4 & 1 & 4 \\
 -8 & 4 & 9
\end{array}
\right)
.
\ee
Figure  \ref{fig:Pyth} illustrates this orbit $\cO$. We can visually verify that the orbit looks thin, and in the next subsection we confirm this rigorously.

\subsection{The Orbit is Thin}\

The group $\SL_2(\R)$ also acts on 
 the hyperbolic upper half-plane 
 $$
 \bH:=\{z=x+iy:x\in\R,y>0\}
 $$
  by fractional linear transformations,
  \be\label{eq:fracLin}
\mattwo abcd:z\mapsto{az+b\over cz+d}.
  \ee
The action of
 %this 
our group $\tilde \G$ in \eqref{eq:tilGamIsP} on
$\bH$
 has a fundamental domain (the definition is similar to \eqref{eq:fundDom}) given  by 
$$
\{z
%=x+iy
\in\bH:
|
\Re(z)
%x
|<1
%2
,\
|z-1/4|>1/4,\
|z+1/4|>1/4
\},
$$ 
where the distances above are Euclidean;
see Figure \ref{fig:HorbA}. The  hyperbolic measure is $y^{-2}\,{dx\,dy%\over y^{2}
}$, and hence this region again has {\it infinite} hyperbolic area. Equivalently, the index of $\tilde\G$ in $\gL(2)$ is infinite, as claimed.
% desired.

 \begin{figure}
        \begin{subfigure}[t]{0.5\textwidth}
                \centering
		\includegraphics[height=1.75in]{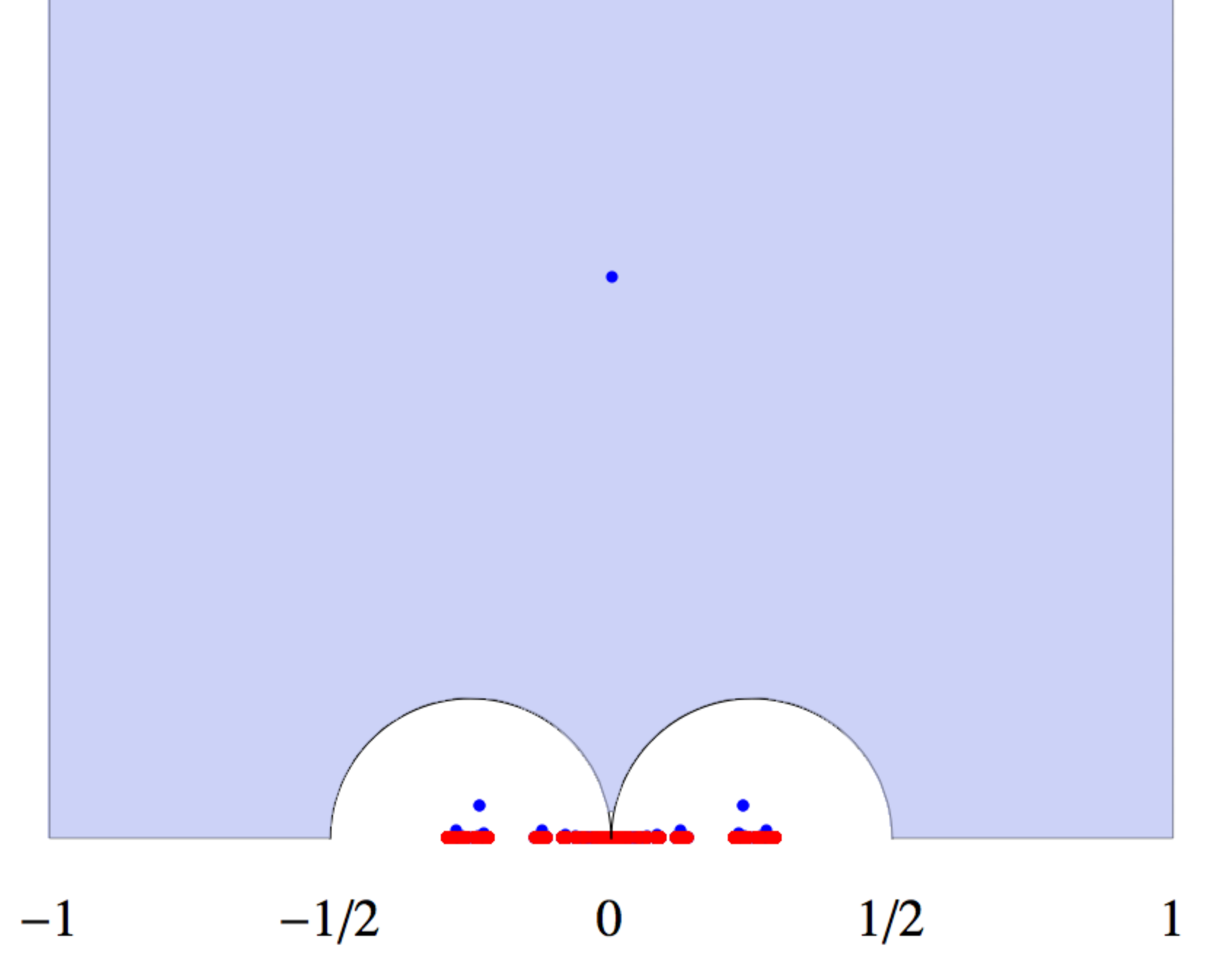}
                \caption{Upper half plane model}
                \label{fig:HorbA}
        \end{subfigure}%
\qquad\qquad
        \begin{subfigure}[t]{1.5in}
                \centering
		\includegraphics[height=1.75in]{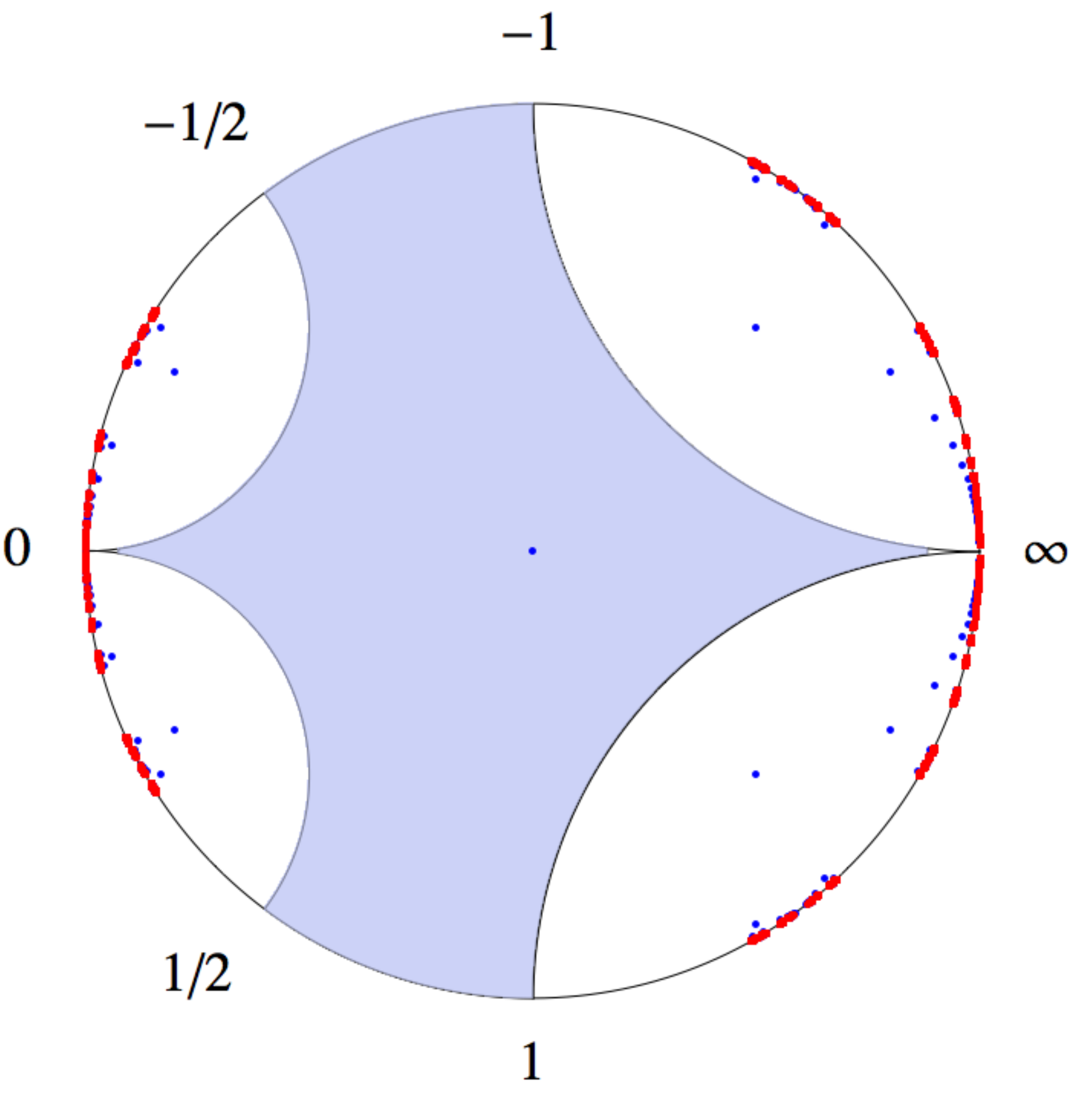}
                \caption{Disk model}
                \label{fig:HorbB}
        \end{subfigure}

\caption{The orbit of $i\in\bH$ under $\tilde\G$.}
\label{fig:Horb}
\end{figure}

%We %can 
%now justify calling $\cO$ a {\it thin} orbit. First we need s
%Some more notation.
%The group $\tilde \G$  in \eqref{eq:tilGamIsP} acts discontinuously on $\bH$, and hence has no limit points in $\bH$. But it does have 
Any orbit of a fixed base point in $\bH$ under $\tilde\G$ has some
{\it limit set} % points 
$\sC=\sC(\tilde\G)$
in the boundary 
$
\dd\bH. %=\R\cup\{\infty\}.
$ 
%The set of such is called the limit set, and is some
A piece of this 
Cantor-like 
%sub
set 
%of $\dd\bH$. A piece of it 
can already be seen in Figure \ref{fig:HorbA}. But to see it %more clearly
fully, we show in Figure \ref{fig:HorbB} the same $\tilde \G$-orbit in the disk model 
$$
\bD=\{z\in\C:|z|<1\},
$$ 
by % applying 
composing the action of $\tilde\G$ with
the map 
$$
\bH\to\bD:z\mapsto{z-i\over z+i}
$$ 
(which encodes the observation that points in the upper half plane are closer to $i$ than they are to $-i$). 
In
% fact, %after a conjugation,
the disk model, one more clearly sees
the limit set as 
 %this is exactly 
 the set of ``directions'' %in the cone 
 in which the orbit $\cO$ can grow -- juxtapose Figure \ref{fig:PythB} %and 
 on
 Figure \ref{fig:HorbB}.
This limit set 
$\sC$
has some Hausdorff dimension $\gd=\gd(\tilde \G)\in[0,1]$; one can estimate 
\be\label{eq:gdApprox}
\gd\approx0.59\dots
\ee 
This
 %number $\gd$ 
 dimension (also called the ``critical exponent of $\G$'')
 is again an % very 
 important geometric invariant, measuring the ``thinness'' of $\G$, as illustrated in the
  following
counting statement
%\begin{thm}[K. 2007 
\cite{MyThesis, Kontorovich2009, KontorovichOh2012}.
%]
Let $\|\bx\|$ be the Euclidean norm. There is some $c>0$ so that 
\be\label{eq:ConeCount}
\#\{\bx\in\cO:\|\bx\|<N\}
\sim c\, N^{\gd},
\qquad\text{
as $N\to\infty$.
}
\ee
%\end{thm}

\begin{comment}
Just a word about the proof: like in the Apollonian counting problem \S\ref{sec:ApCount},%
\footnote{Historically it is the other way around; \eqref{eq:ConeCount} preceded and inspired \eqref{eq:KO}.}
% the proof
 %
 %of \eqref{eq:ConeCount} 
one
 uses infinite-volume spectral theory (Lax-Phillips and Patterson-Sullivan), 
 but unlike %earlier,
 the Apollonian proof,
 the continuous spectrum is handled by some abstract operator theory, which we
  %will not
  have insufficient room to discuss. 
 %The proof uses the fact that $\tilde\G$ has unipotent elements, namely its generators \eqref{eq:tilGamIsP}
A key ingredient in the proof of the above is
that $\tilde\G$ 
%is finitely generated and 
has unipotent elements;
one can remove this restriction,
as is done in 
%, see 
\cite{KontorovichOh2012},
%.
 assuming only that  the dimension $\gd>1/2$ (which is needed to apply Patterson-Sullivan spectral theory; this condition is automatic if there are unipotents).
\end{comment}

%A
Once
again, \eqref{eq:ConeCount} should be compared with the orbit of $\bx_{0}$ under
 the full ambient group,  $\SO_{Q}(\Z)$. %, in which $\G$ is contained. 
 %Standard
 Elementary methods 
 %in automorphic forms or ergodic theory
 %\cite{DukeRudnickSarnak1993, EskinMcMullen1993}
 show
 that
$$
\#\{\bx\in\SO_{Q}(\Z)\cdot\bx_{0}:\|\bx\|<N\}
\sim c\, N
.
%,\qquad\text{
%as $N\to\infty$.
$$
So in passing from the full orbit to
 %our thin orbit 
 $\cO$,
the asymptotic 
%so the magnitude has 
drops from $N$ to $N^{\gd}$, with $\gd<1$. Thus  the orbit $\cO$ is {\it thin}.
\\

%Also, % it follows from
% t
 The fact
  that $\rho$
is a quadratic map in the entries
(see \eqref{eq:gIs})
 implies that the count \eqref{eq:ConeCount} on triples $\bx\in\cO$ is equivalent to the following asymptotic for the pairs $\bv\in\tilde\cO$:
 \be\label{eq:tilOCount}
\#\{
\bv\in\tilde\cO
:
\|\bv\|<N
\}
\sim
c'\cdot N^{2\gd}
,
 \ee
as $N\to\infty$. 
Note that the power of $N$ %has gone from $\gd$ to 
is now
$2\gd$.
This can also be seen immediately from %the fact 
%that 
%observe
 %using
  \eqref{eq:conePyth} and \eqref{eq:PythParam} that 
\be\label{eq:xTov}
\|\bx\|=\sqrt{x^{2}+y^{2}+z^{2}}=\sqrt 2 z
=
\sqrt 2(u^{2}+v^{2})
=
\sqrt 2\|\bv\|^{2}.
\ee
%
%
%So
(Geometrically, the cone \eqref{eq:conePyth} intersects the sphere of radius $N$ at a circle of radius $N/\sqrt2$.)
Observe that
 \eqref{eq:ConeCount} looks like the Apollonian asymptotic \eqref{eq:KO}, while \eqref{eq:tilOCount} 
is more similar to %Zaremba's (or rather, Hensley's) 
Hensley's
estimate \eqref{eq:N2del} in Zaremba's problem. This is just a consequence of choosing between working in the orthogonal group or its spin cover.
\\

\subsection{Diophantine
% Questions
Problems}\

One can now pose a variety of Diophantine
 %problems 
 questions
 about
 the
 values of %some 
% an integral
various
 functions % $f$ 
 on
 such %a 
 thin orbits.
  %$\cO$. 
Given
an orbit $\cO=\G\cdot\bx_{0}$ and 
 a function $f:\cO\to\Z$,
  %which is integral on $\cO$, 
  call
\be\label{eq:sPisP}
\sP:=%\<\bw_{0},\tilde\G\cdot\bv_{0}\>
f(\cO)
\qquad \subset\qquad \Z
\ee
the set of {\it represented} numbers. That is, $n$ is represented by the pair $(\cO,f)$ if there is some $\g\in\G$
   so that
$
n=f(\g\cdot\bx_{0})
.
$ 
And as before, we say $n$ is {\it admissible} if $n\in\sP(\mod q)$ for all $q$. 
For example, if $f$ is the ``hypotenuse'' function, $f(\bx)=z$, one can ask whether $(\cO,f)$ represents infinitely many admissible primes.
Evidence to the affirmative is illustrated in 
Figure \ref{fig:Pyth}, where a triple is highlighted if
 %triples whose 
its hypotenuse is prime. 
 Unfortunately 
this problem on thin orbits\footnote{%
For the full orbit of all Pythagorean triples,
infinitely many hypotenuses are prime. This follows from  \eqref{eq:PythParam} that $z=u^{2}+v^{2}$  and
 Fermat's theorem that all primes $\equiv1(\mod 4)$ are sums of two squares.}
%is 
seems
out of reach of current technology.
 %: are there infinitely many such?

But 
for a restricted class $\cF$ of functions $f$,
and orbits $\cO$ which are ``not too thin,''
recent progress
 %can be
has been  made 
towards the
local-global problem in $\sP$.
 % following 
%problem:
% For a
%restricted class $\cF$ of 
  %given a 
  %functions $f$ which
   %is 
   %integral on $\cO$, one may ask again
%    which numbers appear in $f(\cO)$?
    %
%The above restriction on $%f\in
Let
$\cF$ 
%is
 %simply that
 be the set of functions
  $f$ which are %be
   a linear, %function 
   not on the triples $\bx$, but on the corresponding
   %on the 
   pairs $\bv$.     
For example, it is not particularly well-know
 that
 in a Pythagorean triple,
  the sum of the hypotenuse $z$ and the even side $y$ is always a perfect square.
This follows immediately 
%Observe 
from \eqref{eq:PythParam}; in particular, $y+z=(u+v)^{2}$. So % the values of 
the function 
\be\label{eq:fbxP}
f(\bx)=\sqrt{y+z}=u+v
\ee
is
 integer-valued on $\cO$ and
 linear\footnote{%
Really we want the values of $|u+v|$, which
 within the positive integers
 %is
 are
  the union 
  of the values of $u+v$ and $-u-v$. Alternatively, we can assume that $-I\in\tilde\G$, as is the case for \eqref{eq:tilGamIsP}.%)
}
 in $\bv$.

Another way of saying this is to pass to the corresponding orbit $\tilde\cO=\tilde\G\cdot\bv_{0}$. 
Any such linear function on $\bv$ is of the form
 \be\label{eq:fbvP}
 f(\bv)=\<\bw_{0},\bv\>
,
 \ee
 for some fixed $\bw_{0}\in\Z^{2}$.
%For %our 
In
the
 example
  %above of $f(\bv)=u+v$
\eqref{eq:fbxP}, take 
$%\be\label{eq:bw0P}
\bw_{0}=(1,1)^{t}.
$ %\ee
Then $\cF$ consists of all functions on $\cO$ which, pulled back to $\tilde\cO$, are of the form \eqref{eq:fbvP}.

\phantomsection
\begin{thmp}[Bourgain-K. 2010 \cite{BourgainKontorovich2010}]\label{thm:P}
%Let $\bx_{0}$
%be a primitive, integral Pythagorean triple, and l
%Let
 %or $\bv_{0}$  be as in \eqref{eq:bv0},
%$\cO=\G\cdot\bx_{0}$ 
%and $\tilde\cO$ 
%be % corresponding 
%an
%the
%orbit
%of $\bx_{0}$
%of Pythagorean triples
%under a free, finitely generated %(but possibly infinite index) 
%subgroup $\G<\SO_{Q}(\Z)$
%which contains no unipotent elements.
% which, while allowed to be thin, is not too thin, in
%that
%and let
%$\tilde\G<
%\SL_{2}(\Z)
%\gL(2)$ be a finitely generated group with
Fix any $f\in\cF$
 %as above 
 and let $\sP$ be the set of represented numbers as in \eqref{eq:sPisP}. 
Assume that the %Pythagorean
 orbit $\cO=\G\cdot\bx_{0}$ is not too thin,
 in % the sense 
 that the
 exponent of $\G$ is sufficiently large
 %with Hausdorff dimension 
\be\label{eq:gd0P}
\gd>\gd_{0},
\ee
for some $\gd_{0}<1$. (The value
 $\gd_{0}=
 0.99995$ suffices.)
Then
 almost every admissible number is represented.
%
%Moreover, the
%``exceptional set'' of admissible $|n|<N$ which are not represented has cardinality at most $CN^{1-\vep_{0}}$, for some $%\vep_{0}>0$.
\end{thmp}

We are finally in position to relate this Pythagorean problem to the Apollonian and Zaremba's.
Indeed, passing to the corresponding orbit $\tilde\cO=\tilde\G\cdot\bv_{0}$ and fixing the function $f(\bv)=\<\bw_{0},\bv\>$, we have that 
$n$ is represented 
if there is a $\g\in\tilde\G$ so that
\be\label{eq:nRepP}
n=
\<\bw_{0},\g\cdot\bv_{0}\>
.
\ee
That is,
\be\label{eq:sPis}
\sP=\<\bw_{0},\tilde\G\cdot\bv_{0}\>
,
\ee
which is of the same form as \eqref{eq:DAvGv} and \eqref{eq:BvGv}.
The condition of admissibility is analyzed again given the generators of $\tilde\G$ by strong approximation, Goursat's Lemma, and finite group theory, as in \S\ref{sec:locZ}.

Note that in light of \eqref{eq:tilOCount}, the minimal dimension $\gd_{0}$ in \eqref{eq:gd0P} cannot go below $1/2$: the numbers in $\sP$ up to $N$ (counted {\it with} multiplicity) have cardinality roughly $N^{2\gd}$, so if $\gd$ is less than $1/2$, then certainly a local-global principle fails miserably.
(Such a phenomenon appeared already in the context of % the similarity between this setting and
Hensley's conjecture \eqref{eq:gdcA} in Zaremba's problem.)

%We now leave this specific problem, proceeding in the next section to sketch some of the ingredients going into the proofs of Theorems \hyperref[thm:A]{A}, \hyperref[thm:P]{P}, and \hyperref[thm:Z]{Z}.

\newpage

\section{The Circle Method: Tools and Proofs}\label{sec:C}

We %now 
briefly
review the previous three sections, unifying the (re)formulations of the problems.
The Apollonian, Pythagorean, and Zaremba Theorems will henceforth be referred to as
Theorem $X$, where
$$
X=
\hyperref[thm:A]{A}, 
\hyperref[thm:P]{P}, \text{ or }
\hyperref[thm:Z]{Z},
$$
respectively.
Theorem $X$ concerns 
the set $\sS$ of numbers of the form
\be\label{eq:sSis}
\sS=
\<\bw_{0},\G\cdot\bv_{0}\>
.
\ee
Here 
$$
\sS
=
\threecase
{\text{the set $\sB$ of curvatures }\eqref{eq:BvGv}}
{if $X=\hyperref[thm:A]{A}$,}
{%
\begin{array}{l}
\text{the set $\sP$ of square-roots of sums of 
%hypotenuses and even sides
}\\
\text{%the set $\sP$ of square-roots of sums of 
\quad hypotenuses and even sides }\eqref{eq:sPis}
\end{array}
}
{if $X=\hyperref[thm:P]{P}$,}
{\text{the set $\sD_{\cA}$ of denominators }\eqref{eq:DAvGv} }
{if $X=\hyperref[thm:Z]{Z}$,}
$$
$$
\G
=
\threecase
{\text{the Apollonian group $\G$}}
{if $X=\hyperref[thm:A]{A}$,}
{\text{an infinite index subgroup  $\tilde\G<\gL(2)$}}
{if $X=\hyperref[thm:P]{P}$,}
{\text{the semigroup $\G_{\cA}$}}
{if $X=\hyperref[thm:Z]{Z}$,}
$$
$$
\bv_{0}
=
\threecase
{\text{the root quadruple}}
{if $X=\hyperref[thm:A]{A}$,}
{%(2,1)^{t}
\text{%, or 
any coprime pair of opposite parity}}
{if $X=\hyperref[thm:P]{P}$,}
{(0,1)^{t}}
{if $X=\hyperref[thm:Z]{Z}$,}
$$
and
$$
\bw_{0}
=
\threecase
{\text{a standard basis vector $\bbe_{j}$}}
{if $X=\hyperref[thm:A]{A}$,}
{%(1,1)^{t}
\text{%, or 
any fixed pair}}
{if $X=\hyperref[thm:P]{P}$,}
{(0,1)^{t}}
{if $X=\hyperref[thm:Z]{Z}$.}
$$
But now we can forget the individual %situations
%settings
problems
 and just focus on the general setting \eqref{eq:sSis}; %, so 
 one need not keep the above taxonomy in one's head throughout. 

To study the local-global problem for $\sS$, we introduce the representation function
\be\label{eq:cRNdef}
\cR_{N}(n)
:=
\sum_{\g\in\gW_{N}}
\bo_{\{n=\<\bw_{0},\g\cdot\bv_{0}\>\}}
.
\ee
Here $N$ is a growing parameter, and $\gW_{N}$ is a certain subset of the radius $N$ ball in $\G$,
$$
\gW_{N}
\subset
\{
\g\in\G
:
\|\g\|<N
\}
,
$$
which we will describe in more detail later.
For now, one can just think of %it 
$\gW_{N}$
as the whole %$\G$-ball of 
radius $N$ ball. To get our bearings, 
let us recall roughly the size of $\gW_{N}$.
It is convenient to
%us 
introduce the parameter $\ga$, defined by
$$
\ga:=
\threecase
{\gd, \text{ the %Hausdorff 
dimension of an Apollonian packing}}
{if $X=\hyperref[thm:A]{A}$, see \eqref{eq:gdIs}}
{2\gd,\text{ where $\gd$ is the dimension of %the limit set of 
$\sC(\tilde\G)$}}
{if $X=\hyperref[thm:P]{P}$, see \eqref{eq:gdApprox}}
{2\gd_{\cA},\text{ where $\gd_{\cA}$ is the dimension of %the Cantor set 
$\sC_{\cA}$}}
{if $X=\hyperref[thm:Z]{Z}$, see \eqref{eq:gdZ}.}
$$
In each case $\ga$ satisfies
\be\label{eq:ga}
1<\ga<2.
\ee
Then the cardinality of such a ball $\gW_{N}$ is roughly
$$
\#\{
\g\in\G
:
\|\g\|<N
\}
\asymp
\threecase
{N^{\gd},}%\text{ $\gd=$H}}
{if $X=\hyperref[thm:A]{A}$, see \eqref{eq:KO}}
{N^{2\gd},}%(1,1)^{t}\text{, or any base pair}}
{if $X=\hyperref[thm:P]{P}$, see \eqref{eq:tilOCount}}
{N^{2\gd_{\cA}},}%(0,1)^{t}}
{if $X=\hyperref[thm:Z]{Z}$, see \eqref{eq:N2del}.}
.
$$
(Technically the quoted results are about counting in the corresponding orbits $\cO$ and not in the groups $\G$; 
but the order of magnitude is the same for both.)
We can write this uniformly by
 %saying that
giving  the cardinality of $\gW_{N}$ % is of order $N^{\ga}$. 
as
\be\label{eq:gWNsize}
|\gW_{N}|\asymp N^{\ga}.
\ee

Returning to \eqref{eq:cRNdef}, we see by construction that $\cR_{N}$ is nonnegative. Moreover observe that 
\be\label{eq:RNn0}
\text{if }
\cR_{N}(n)>0,
\text{
then certainly  $n$ is represented. %Note furthermore 
}
\ee
Also record that 
\be\label{eq:RNsupp}
\text{%
$\cR_{N}$ is supported on $n$ of size 
}
|n|\ll
%<c\cdot 
N.
\ee 

Recalling the notation $e(x)=e^{2\pi ix}$, the Fourier transform
\bea
\nonumber
\cS_{N}(\gt)
&:=&
\widehat{
\cR_{N}}(\gt)
=
\sum_{n\in\Z}
\cR_{N}(n)
e(n\gt)
\\
\label{eq:cSNdef}
&=&
\sum_{\g\in\gW_{N}}
e(\gt\<\bw_{0},\g\cdot\bv_{0}\>)
\eea
is a wildly oscillating exponential sum on the circle $\R/\Z=[0,1)$, whose graph looks something like Figure \ref{fig:Exp}.
%Of course o
One recovers $\cR_{N}$ through elementary Fourier inversion,
\be\label{eq:RNSN}
\cR_{N}(n)
=
\int_{\R/\Z}
\cS_{N}(\gt)
e(-n\gt)
d\gt
,
\ee
%but in a sense, we are just back to where we started, namely, nowhere.
but without further ingredients, one is going around in circles (no pun intended).

\begin{figure}
\includegraphics[width=\textwidth]{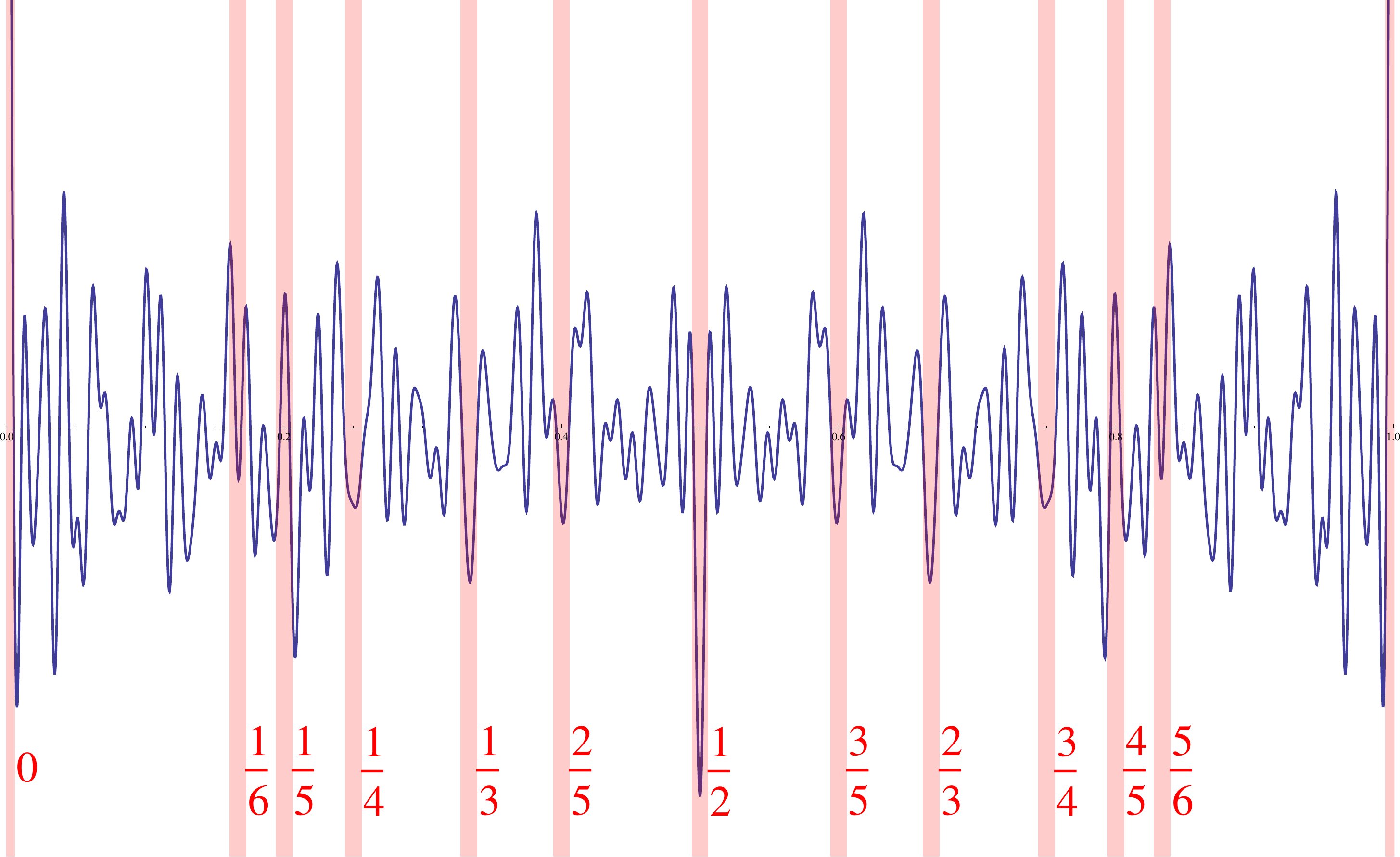}
\caption{The real part of an exponential sum of the form \eqref{eq:cSNdef}}
\label{fig:Exp}
\end{figure}

Building on the Hardy-Ramanujan technique for asymptotics of the partition function, Hardy and Littlewood
%advanced the following procedure for 
 had the 
%wild 
idea that
the bulk of the integral \eqref{eq:RNSN} %can 
could
be captured
just by integrating over
 frequencies $\gt$ that are
 very close to rational numbers $a/q$, $(a,q)=1$, with very small denominators $q$; some of these
 intervals are
 %is region is 
 shaded in Figure \ref{fig:Exp}. 
 These are now called the {\it major arcs} $\fM$; 
 %though they
 %Indeed,
%of course 
the name refers
not to their total length
(% despite comprising 
they comprise
an ever-shrinking fraction of the circle $\R/\Z$)
but
 to the fact that 
 they are supposed to account for a preponderance of % the integral \eqref{eq:RNSN}.
 $\cR_{N}(n)$.
Accordingly, we decompose \eqref{eq:RNSN} as
 $$
 \cR_{N}(n)=\cM_{N}(n)+\cE_{N}(n)
 ,
 $$
where the major arc contribution
\be\label{eq:cMNis}
\cM_{N}(n):=
\int_{\fM}
\cS_{N}(\gt)
e(-n\gt)
d\gt
\ee
is supposed to give the ``main'' term, and 
\be\label{eq:cENis}
\cE_{N}(n):=
\int_{\fm}
\cS_{N}(\gt)
e(-n\gt)
d\gt
\ee
should be the ``error''. Here $\fm:=[0,1)\setminus\fM$ are the complemenary so-called {\it minor arcs}.
If $\cM_{N}(n)$ is positive and
bigger than %satisfies $\cM_{N}(n)>
 $|\cE_{N}(n)|$, then certainly 
 \be\label{eq:RNge0}
 \cR_{N}(n)\ge \cM_{N}(n)-|\cE_{N}(n)|>0
,
 \ee 
 so again, $n$ is represented. 
In practice, one 
typically % proves 
tries to prove
an asymptotic formula (or at least a lower bound) for $\cM_{N}$, 
and % then tries to
then give 
an upper bound for $|\cE_{N}|$.

% , since
%e
The reason for this decomposition is that exponential sums such as $\cS_{N}$ should
be mostly  supported on $\fM,$
 having their biggest peaks and valleys
at (or very near) these frequencies (some of this phenomenon is visible in Figure \ref{fig:Exp}). 
 %Certainly 
 Indeed,
 the value $\gt=0$ is as big as $\cS_{N}$ will ever get,
 \be\label{eq:SNtriv}
| \cS_{N}(\gt)|\le \cS_{N}(0)=|\gW_{N}|,
 \ee
 %since 
%simply %because 
which follows trivially (and is thus called the {\it trivial bound})
from the triangle inequality:
every summand in \eqref{eq:cSNdef} is a complex number of absolute value $1$.
%But a
Also for other $\gt\in\fM,\ \gt\approx a/q$, the summands %of $\cS_{N}(\gt)$ 
should %also 
%somehow 
all point in a limited number of directions, colluding to
give a large contribution to $\cS_{N}$.
As we will see later, 
%there points 
at these frequencies, one is 
in %some
a  sense
measuring
  %how
  the distribution of
  $\sS$
  (or equivalently
   $\gW_{N}$) %is distributed
%  or rather, our set $\sS$,
   along 
   certain
   arithmetic progressions. %, % of modulus $q$, 
%   as we will see later.
   This strategy of coaxing out the (conjectural) main term for $\cR_{N}$
    %in this way 
    works in surprisingly great generality,
but can also give
false predictions (even for the Prime Number Theorem, see e.g. \cite{Granville1995}).
%misleading answers in some contexts).    
    % giving sharp (conjectural) asymptotics
%for everything from the Twin Prime Conjecture to 
%, Goldbach,
%and infinitely  many other  simil
%or your favorite problem.    
% 
% As we will see more clearly below, 
\\

%Also f
Having % agreed on 
made
this decomposition, %what do we expect for the 
we should determine what we expect for the main term.
From \eqref{eq:cSNdef}, we have that
$$
\sum_{n}\cR_{N}(n)=\cS_{N}(0)=|\gW_{N}|
,
$$
so 
recalling the support
%since $\cR_{N}$ is supported on %$|n|\ll N$
\eqref{eq:cRNdef}
of $\cR_{N}$, one might expect that 
an admissible number of size about $n\asymp N$ is represented roughly $|\gW_{N}|/N$ times. 
In particular, since every admissible number is expected to be represented, one would like to show, say, for $N/2\le n<N$, that
\be\label{eq:cMNbnd}
\cM_{N}(n)\gg \fS(n) {|\gW_{N}|\over N}.
\ee
Here $\fS(n)\ge0$ is a  certain product of local densities
 %arithmetic number 
 called the {\it singular series}, and will be discussed at greater length later.
 It alone is responsible for the notion of admissibility, vanishing on non-admissible $n$.
 For admissible $n$, it typically does not fluctuate too much; crudely one can show in many contexts 
the lower bound
 %that
 %, which measures the local density, and does not fluctuate too much. 
 %In particular, $\fS(n)$ %is zero 
 %vanishes
 %if $n$ is not admissible, and  otherwise is
  $\gg N^{-\vep}$ for any $\vep>0$. 
 %The singular series is well-understood, and f
 For ease of exposition, %we pretend henceforth 
 let us just pretend for now 
% that
 that every $n$ is admissible and remove the role of the singular series,
 %, which shall be clarified later, 
 %assuming 
 allowing ourselves to assume %just 
 that
 \be\label{eq:fSn}
 \fS(n)=1
.% ,
 \ee
%say.
%\\
Observe also that, in light of \eqref{eq:gWNsize} and \eqref{eq:ga}, the lower bound 
in \eqref{eq:cMNbnd} is of the order $N^{\ga-1}$, with $\ga>1$. That is, there should be %very many 
quite a lot of
representations
of an admissible $n\asymp N$
 %for $N$ 
 large, giving further %evidence
 indication
 that
  %for 
  every sufficiently large admissible number %being 
%  should
may be  represented.
\\

One is then left with the problem of estimating away the remainder %of the integral
term $\cE_{N}$,
and this is why (as Peter Sarnak likes to say) the circle method is a ``method'' and not a ``theorem'':
establishing such estimates is much more of an art than a science.
%
% for which 
The 
 Hardy-Littlewood procedure suggests somehow 
 exploiting the fact that on the minor arc frequencies, $\gt\in\fm$, the exponential sum $\cS_{N}$ in \eqref{eq:cSNdef}
 should itself already be quite small, 
 being a sum of canceling phases.
% %
% from the wildly spinning exp
 %
% establishing that, for frequencies $\gt$
  %not in the major arcs (the compliment of which is appropriately called the ``minor arcs'')
% in the minor arcs $\fm$, there 
%is so much cancellation from %the wildly irregular 
 %
% individual summands in \eqref{eq:cSNdef}
% actually 
 %For the minor arcs integral
%  for the former and upper bounds for the latter. 
%We first discuss how to approach
If one
could indeed prove 
 at the level of individual $n$
 an upper bound for the error term
  $\cE_{N}$
%one wishes to give an upper bound,
which is
 asymptotically smaller than the lower bound 
 %given 
\eqref{eq:cMNbnd}
 for $\cM_{N}$%
 %is able to do this 
, then one could 
immediately
conclude  that every sufficiently large admissible $n$ is represented.
Unfortunately, at present we do not know how to give such strong upper bounds on the minor arcs.

Instead, we settle for an ``almost'' local-global statement, by proving a sharp bound not for individual $n$, but 
for $n$ in an average sense, %in particular in $L^{2}$
as follows.
Parseval's theorem states that the $L^{2}$ norm of a function is %the same as 
equal to
that of its Fourier transform, that is, the Fourier transform is a unitary operator on these Hilbert spaces.
Using the definition
\eqref{eq:cENis},
Parseval's theorem then % tells us that
gives
\be\label{eq:Pars}
\sum_{n}|\cE_{N}(n)|^{2} = \int_{\fm}|\cS_{N}(\gt)|^{2}d\gt
.
\ee
Inserting our trivial bound \eqref{eq:SNtriv} for $\cS_{N}%(\gt)
$
into the above
 %into \eqref{eq:cSNdef} 
%is just $\le|\gW_{N}|$, giving
%gives 
yields
a trivial bound for \eqref{eq:Pars} of  
\be\label{eq:ParsTriv}
\int_{\fm}|\cS_{N}(\gt)|^{2}d\gt
\le
|\gW_{N}|^{2}%$ for \eqref{eq:Pars}
.  
\ee
We claim that it suffices for our applications to
establish
a bound of the form
\be\label{eq:E2}
\int_{\fm}|\cS_{N}(\gt)|^{2}d\gt
=o\left(
{|\gW_{N}|^{2}\over N}
\right)
.
\ee
That is, the above saves a little more than $\sqrt N$ on average over $\fm$ off of each term $\cS_{N}$ 
relative to the trivial bound \eqref{eq:ParsTriv}.
%i%n
 %\eqref{eq:Pars}. 
 %Let us 
We first explain why this % would be enough.
suffices.
 %the minor arcs, 
% that is, we seek

Let $\fE(N)$ be the set of exceptional $n$
 %up to $N$, that is, 
 (those that are admissible but not represented) in the range $N/2\le n<N$. 
%Certainly the $n$th Fourier coefficient \eqref{eq:SNhat} is nonzero if $\cM_{N}(n)>|\cE_{N}(n)|$, so 
Recalling \eqref{eq:RNge0},
the number of exceptions is bounded by
$$
\#\fE(N)
\le 
\sum_{N/2<|n|<N\atop n \text{ is admissible}}
\bo_{\{|\cE_{N}(n)|\ge \cM_{N}(n)\}}
.
$$
For  admissible $n$, we have the major arc lower bound \eqref{eq:cMNbnd} and recall our simplifying assumption \eqref{eq:fSn}; thus
\be\label{eq:fEn}
\#\fE(N)
\le
\sum_{n}%|n|<N\atop n \text{ is admissible}}
\bo_{\{|\cE_{N}(n)|\gg |\gW_{N}|/N\}}
.
\ee
Here is a pleasant (standard) trick: for those $n$ contributing a $1$ rather than $0$ to \eqref{eq:fEn}, we have
$$
1\ll{|\cE_{N}(n)|\over |\gW_{N}|/N},
$$ 
both sides of which may be squared. Hence \eqref{eq:fEn} implies that
$$
\#\fE(N)
\ll
{N^{2}\over |\gW_{N}|^{2}}\cdot
\sum_{n}|\cE_{N}(n)|^{2}
.
$$
Now we  apply Parseval \eqref{eq:Pars} and the bound \eqref{eq:E2} which we had claimed would suffice. This gives
$$
\#\fE(N)
=
o\left(
{N^{2}\over |\gW_{N}|^{2}}\cdot
{|\gW_{N}|^{2}\over N}
\right)
=
o(N)
,
$$
and thus
% almost every 
$100\%$ of the
admissible numbers in the range $[N/2,N)$ %is 
are
represented.
%, as claimed. 
%But of course 
Combining %these 
such
dyadic %pieces
intervals, % gives the claim.
we conclude that almost every admissible number is represented.
\\

Now ``all'' that is left is
to establish the major arcs bound \eqref{eq:cMNbnd} and the error bound \eqref{eq:E2}.
In the next two subsections, we focus individually on the tools needed to prove these claims.
\\

\subsection{The Major Arcs}\

%We now % sketch an argument 
%give some of the ingredients going into establishing
%leading to 
%\eqref{eq:cMNbnd}. 
Recall that $\cM_{N}$ in \eqref{eq:cMNis} is an integral over the major arcs $\gt\in\fM$, 
%and we %may 
%pretend 
%assume
%for simplicity that
with $\gt$
 %is just equal to 
 very close to
 a fraction $a/q$, with $q$ %very
  ``small'' (the meaning of which is explained below). %relative to $N$, as will be explained below. 
Also let us pretend for now that $\gW_{N}$ is just the whole $\G$-ball, 
$$
\gW_{N}=\{\g\in\G:\|\g\|<N\}.
$$ 
%Then % the analysis reduces  to 
%evaluating
%Then let us %first 
We
begin by trying to 
evaluate
 \eqref{eq:cSNdef} at $\gt=a/q$: %. P
%Then
%gives
$$
\cS_{N}\left(\frac aq\right)
=
\sum_{\g\in\G\atop\|\g\|<N}
e\left(
\frac aq
\<\bw_{0},
\g\cdot
\bv_{0}
\>
\right)
.
$$
An important observation in the above is that
%Observe that in the
 %exponent 
% above,
 the summation may be grouped according to the residue class mod $q$ of
the integer 
$\<\bw_{0},
\g\cdot
\bv_{0}
\>
.
$
%only needs to be defined $\mod q$.
%It is a pleasant
 %observation then 
%Then 
%Hence 
%Thus
%we have the pleasant fact
% that 
Or what is essentially the same,
we can decompose the sum
 %above
  according to the  
%
 %only depends on 
 the residue class of $\g(\mod q)$.
To this end,
 let $\G_{q}=\G(\mod q)$ be the set of such residue classes (which we have already studied %extensively 
 in the context of admissibility and strong approximation). Then we split the sum as
\be\label{eq:SNbreak}
\cS_{N}\left(\frac aq\right)
=
|\gW_{N}|
\sum_{\g_{0}\in\G_{q}}
e\left(
\frac aq
\<\bv_{0}\cdot\g_{0},
\bw_{0}
\>
\right)
\cdot
\left[
\frac1{|\gW_{N}|}
\sum_{\g\in\G\atop\|\g\|<N}
\bo_{\{\g\equiv\g_{0}(\mod q)\}}
\right]
,
\ee
where we have artificially multiplied and divided by the cardinality of $\gW_{N}$.
Now for $\g_{0}$ fixed,  the bracketed term is %clearly 
then
measuring the ``probability'' that $\g\equiv\g_{0}(\mod q)$.
%, and to understand how 
%We claim 
As one 
%the reader
may 
%have already guessed, % be expected, 
suspect,
our groups do not have particular preferences for certain residue classes over others; that is,
this probability 
%that 
%this is becoming 
becomes
equidistributed
as $N$ grows, %as long as $q$ does not grow too 
with $q$ also allowed to grow, but at a much %smaller 
slower
rate. (In fact, this is exactly what
 %determines the size of
we mean by the denominator $q$ being ``small''  -- relative to $N$ -- in the major arcs $\fM$.)
To explain %what this means
how this happens, we %need to 
briefly
discuss the notion of an {\it expander}. % graphs}.
\\
%,and for this one needs the concept of an expander
%To explain 

\begin{figure}
%$
%\begin{tabular}{cc}
\hskip-2.5in
\includegraphics[width=.5\textwidth]{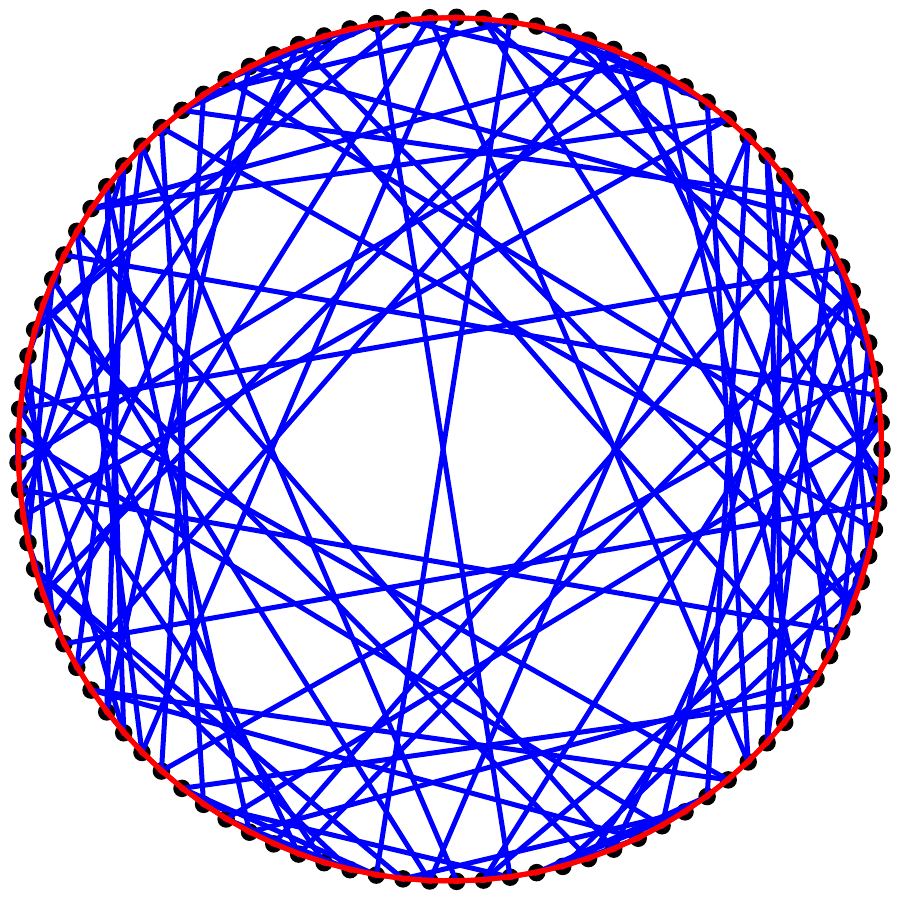}\vskip-1.5in
$
\begin{array}{l}
%\beann
\hskip2.75in{\color{red}\mattwo1{\pm2}01:}{\color{red}x\mapsto x\pm2}
%$$
\\
%$$
\hskip2.75in{\color{blue}\mattwo10{\pm4}1:}{\color{blue}x\mapsto x(\pm4x+1)^{-1}%(\mod p)
}
%$$
\end{array}
$
%\eeann
\vskip.75in
%\end{tabular}
%$
\caption{%Illustrating a
An expander; % graph
%, for 
shown with $q=101$% relative to the generators \eqref{eq:tilGamIsP}
}
\label{fig:Expand}
\end{figure}

Rather than going into the general theory (for which we refer the reader to the %thorough 
beautiful
survey \cite{Lubotzky2012}; see also \cite{Sarnak2004}), we content ourselves with but one illustrative example of expansion.
Figure \ref{fig:Expand} shows the following graph.
 %$\cG$. 
 %Let 
 For $q=101$, say,
  %and 
  take the vertices to be the elements of $\Z/q\Z$, %plotted 
 organized
% in the obvious way 
 around the unit circle by placing $x\in\Z/q\Z$ at $e(x/q)$. 
 %The edges are given as follows: 
For the edges,
 connect each 
 \be\label{eq:edges}
 x\text{ to }x\pm2\text{, and also to }x(\pm4x+1)^{-1},
 \ee 
 when %the element is 
 %invertible
 inversion $(\mod q)$ is possible. 
%Of course t
This is nothing %but 
more than
the fractional linear action (see \eqref{eq:fracLin}) of the generating matrices in \eqref{eq:tilGamIsP}
(and their inverses)
on $\Z/q\Z$.
%
%Expanders are %certain 
%graphs (or rather families of such)
%which, while being ``sparse,'' are simultaneously ``highly connected'', as we now explain.
%%
%
%First w
We
first
 claim
that
 %this 
our
graph
on $q$ vertices
 is ``sparse''. Indeed, %it has $q$ vertices, and
  the complete graph (connecting any vertex to any other) has 
  on the order of $q^{2}$ 
%  $\left({q\atop2}\right)\asymp q^{2}$
  edges, whereas our graph has only 
  on the order of 
  $
  %\asymp 
  q$ edges (since \eqref{eq:edges} implies that
  %by 
   any vertex is connected to at most %$4$ 
  four
  others). So we have
square-root the total number of possible edges, and our
%. So the 
graph is indeed quite sparse. 
% the graph is highly connected. 

%We also 
%
%Next we % claim that
%claim that, d
Despite having few edges, it is a fact that this graph is 
nevertheless
highly connected, 
%
%The condition of being is 
%roughly 
%meaning 
in the sense 
that a random walk on it % this graph 
is rapidly %uniformly 
mixing.
Moreover, 
this rate of mixing, properly normalized, is independent of the choice of $q$ above. That is, 
by varying $q$,
we in fact have a whole family
of such sparse but highly connected graphs, 
and
with a uniform mixing rate; this is exactly what characterizes an expander.
%\\
%were we to take a family of such graphs, with $q$ %=101$ replaced by arbitrarily l
%running 
%moving
%through the integers, this %uniform 
%mixing rate would be {\it uniform} over $q$. 

%So this roughly characterizes an expander: it is a
%family of 
%a sparse but uniformly highly connected graphs. %, with uniform mixing rate for a random walk.
%The
Proofs of % the
% mixing statement above
  %is proved
expansion   use, %ing 
among other things,   tools from additive combinatorics, %namely 
in particular, the
so-called sum-product \cite{BourgainKatzTao2004, Bourgain2008} and triple-product \cite{Helfgott2008, BreuillardGreenTao2010, PyberSzabo2010} estimates, and %a lot of fter 
quite a lot of other work (see e.g. \cite{
SarnakXue1991, 
Gamburd2002, 
BourgainGamburd2008, BourgainGamburdSarnak2010, Varju2010, BourgainVarju2011, SalehiVarju2011}).
Once one proves uniform expansion for
 %these
 such finite graphs, 
 %it is converted into
 the statements %are
 must be
  converted into
  the
archimedean 
 form needed 
 for the bracketed term in
% in
%the major arcs 
\eqref{eq:SNbreak}% integral
.
To handle such counting statements, one
 uses
 $$
 \threecase
 {
 \begin{array}{l}
 \text{infinite volume spectral and representation theory}\\
\text{ \`a la  \S\ref{sec:ApCount}, specifically 
 Vinogradov's thesis 
  \cite{Vinogradov2012},
}
\\
\
\end{array}
}
 {if $X=\hyperref[thm:A]{A}$,}
 {
 \begin{array}{l}
 \text{% 
 similar techniques %using 
 developed by}\\
 \text{
 Bourgain-K.-Sarnak \cite{BourgainKontorovichSarnak2010},}
 \end{array}
 }
 {if $X=\hyperref[thm:P]{P}$,}
 {
 \begin{array}{l}
 \
 \\
 \text{%
 the thermodynamic formalism,
  analytically continuing 
 }
 \\
\text{ certain  
Ruelle transfer operators \cite{Lalley1989, Dolgopyat1998, Naud2005} 
}
\\
\text{ 
and their 
``congruence''  extensions% in the semigroup
 ; see
\cite{BourgainGamburdSarnak2011},
}\end{array}
}
 {if $X=\hyperref[thm:Z]{Z}$.}
 $$
\begin{comment}
 in the group case 
in the Apollonian and Pythagorean cases (%
$X=\hyperref[thm:A]{A}$ or $\hyperref[thm:P]{P}$; that is, when $\G$ is a group),
and
for Zaremba's semigroup case $X=\hyperref[thm:Z]{Z}$,
one uses 
 the thermodynamic formalism,
  analytically continuing certain Ruelle transfer operators \cite{Lalley1989, Dolgopyat1998, Naud2005} and their 
``congruence''  extensions% in the semigroup
 ; see
\cite{BourgainGamburdSarnak2011}.
\end{comment}

Without going into details, the upshot is that, up to 
%up to 
acceptable error, 
%implying that
 the bracketed term in \eqref{eq:SNbreak} is
just
 $%$
%\left[
{1/%\over
 |\G_{q}|}
%\cdot
%\sum_{\g\in\G\atop\|\g\|<N}
%1
%\right]
%.
,
 $ %$
confirming the desired equidistribution.
Inserting this %analysis 
estimation
%of $\cS_{N}$ 
into $\cM_{N}$ in \eqref{eq:cMNis},
one uses these 
techniques and some more standard circle method analysis 
%lead 
to 
eventually %to
conclude
 \eqref{eq:cMNbnd}.
%\\

\subsection{The Minor Arcs}\

%\newpage

We use different strategies to prove \eqref{eq:E2} for  the 
Pythagorean and Zaremba settings
$X=\hyperref[thm:P]{P}$ or
$\hyperref[thm:Z]{Z}$,
% setting of \S\S\ref{sec:GAFA}-\ref{sec:Zaremba} %and
versus
 the Apollonian setting
  %\S\ref{sec:Apoll}, 
  $X=\hyperref[thm:A]{A}$,
  so we present them % separately.
individually.

\subsubsection{%The 
Pythagorean and Zaremba settings%
 %setting of \S\S\ref{sec:GAFA}-\ref{sec:Zaremba}
}\

To handle the minor arcs here, we make the observation   that the ensemble $\gW_{N}$ in the definition of  $\cS_{N}$ from \eqref{eq:cSNdef} 
need not be a full $\G$-ball. 
In fact, the definition of $\cS_{N}$
can be changed to, say,
\be\label{eq:cSNdef2}
\cS_{N}(\gt)
:=
\sum_{\g_{1}\in\G\atop\|\g_{1}\|<\sqrt N}
\sum_{\g_{2}\in\G\atop\|\g_{2}\|<\sqrt N}
e(\gt\<\bv_{0}\g_{1}\g_{2},\bw_{0}\>)
,
\ee
without
irreparably 
 damaging the major arcs analysis.
This new sum encodes much more of the (semi)group structure of $\G$, while %still having 
preserving
the property
\eqref{eq:RNn0}, where
% that 
%if we redefine
 $\cR_{N}$ is redefined by \eqref{eq:RNSN}.
 %, 
%a non-vanishing $n$th Fourier coefficient % implies that 
%if its $n$th Fourier coefficient does not vanish, then 
%forces 
%implies that
%$n$  %is represented by $\G$. 
%to be 
%is represented.
(In reality, we use %an 
even % a slightly 
more complicated exponential sums.) The advantage of \eqref{eq:cSNdef2} is that
we can now
 exploit this %e 
%bilinear
 structure
%of \eqref{eq:cSNdef2}
 \`a la 
 Vinogradov's
  %approach
  method 
  \cite{Vinogradov1937}
  for estimating bilinear forms.
Just one such maneuver is the following.

%Keeping
% we
%Fixing one of the variables, we can 
%a
Apply the Cauchy-Schwarz inequality 
to \eqref{eq:cSNdef2} 
in
the
$\g_{1}$ %fixed
variable,
% %
% the other, 
%as follows.
%
 estimating
$$
|\cS_{N}(\gt)|
\le
\left(
\sum_{\g_{1}\in\G\atop\|\g_{1}\|<\sqrt N}
1
\right)^{1/2}
\left(
\sum_{\g_{1}\in\SL_{2}(\Z)\atop\|\g_{1}\|<\sqrt N}
\left|
\sum_{\g_{2}\in\G\atop\|\g_{2}\|<\sqrt N}
e(\gt\<\bv_{0}\g_{1}\g_{2},\bw_{0}\>)
\right|^{2}
\right)^{1/2}
.
$$
Notice in the second appearance of $\g_{1}$, we have replaced 
%
%and begin to
%
%  Here is just one example 
%  
   %In so doing, we pass at some point from 
   the thin and mysterious group $\G$ (or semigroup $\G_{A}$) %to 
by   the full ambient group $\SL_{2}(\Z)$. 
On one hand, this %lets
allows
 us
 to
now
 use
  more classical tools to get the requisite cancellation in the minor arcs % sum.   
  integral.
On the other hand,   
   this type of perturbation argument only succeeds when $\gd$ is %very 
 near $1$, 
 %thereby 
 explaining the restrictions \eqref{eq:gdMin} and \eqref{eq:gd0P}. 
%There are many more ingredients

%\newpage

 \subsubsection{The Apollonian case% \S\ref{sec:Apoll}
 }\
 
%One cannot apply t
The above strategy fails for the Apollonian problem, because the Hausdorff dimension \eqref{eq:gdIs} is a fixed invariant which
% cannot be adjusted.
 refuses to be adjusted to suit our needs.
Instead, we 
recall that the Apollonian group $\G$ contains the special (arithmetic) subgroup $\Xi$ from \eqref{eq:XiDef}. 
Then, like \eqref{eq:cSNdef2}, we change the definition of the exponential sum to something of the form
\be\label{eq:SN3}
\cS_{N}(\gt):=
\sum_{\xi\in\Xi\atop\|\xi\|<X}
\sum_{\g\in\G\atop\|\g\|<T}
e(\gt\<\bv_{0}\cdot\xi\,\g,\bw_{0}\>)
,
\ee 
for certain parameters $X$ and $T$ chosen optimally in relation  to $N$.
One uses the full sum over the group $\G$ to capture the major arcs and % the %full
 admissibility conditions.
For the minor arcs bound,  one keeps $\g$ fixed and uses the
 %fact that
 classical %algebraic 
 arithmetic
 group 
  $\Xi$ 
  %is a 
  %
to get sufficient cancellation % in \eqref{eq:SN3} and \eqref{eq:Pars} (a type of Kloosterman refinement is used)
 to prove the desired bound \eqref{eq:E2}.
% This completes our sketch of the argument.
%\newpage

\newpage

%\subsection{Thermodynamical Formalism}\label{sec:therm}

%\newpage

%\newpage
\bibliographystyle{%
%plain}
alpha}

\bibliography{../../AKbibliog}

\

\end{document}